\documentclass[11pt]{amsart}
\usepackage{amsmath,amssymb,amsthm}
\usepackage[latin1]{inputenc}
\usepackage{version,tabularx,multicol}
\usepackage{graphicx,float}

\newcommand{\reff}[1]{(\ref{#1})}

\theoremstyle{plain}
\newtheorem{theo}{Theorem}[section]

\newtheorem{cor}[theo]{Corollary}
\newtheorem{prop}[theo]{Proposition}

\newtheorem{lem}[theo]{Lemma}
\newtheorem{defi}[theo]{Definition}
\theoremstyle{remark}
\newtheorem{rem}[theo]{Remark}

\newcommand{\cb}{{\mathcal B}}

\newcommand{\ce}{{\mathcal E}}

\newcommand{\cf}{{\mathcal F}}

\newcommand{\ch}{{\mathcal H}}
\newcommand{\cj}{{\mathcal J}}

\newcommand{\cl}{{\mathcal L}}
\newcommand{\cn}{{\mathcal N}}
\newcommand{\cm}{{\mathcal M}}
\newcommand{\cp}{{\mathcal P}}
\newcommand{\cq}{{\mathcal Q}}
\newcommand{\cs}{{\mathcal S}}

\newcommand{\ct}{{\mathcal T}}
\newcommand{\cx}{{\mathcal X}}

\newcommand{\E}{{\mathbb E}}

\newcommand{\N}{{\mathbb N}}
\renewcommand{\P}{{\mathbb P}}

\newcommand{\R}{{\mathbb R}}
\renewcommand{\S}{{\mathbb S}}

\newcommand{\rN}{{\rm N}}
\newcommand{\rP}{{\rm P}}

\newcommand{\ind}{{\bf 1}}

\newcommand{\Supp}{{\rm Supp}\;}
\newcommand{\Card}{{\rm Card}\;}

\newcommand{\val}[1]{\mathop{\left| #1 \right|}\nolimits}
\newcommand{\inv}[1]{\mathop{\frac{1}{ #1}}\nolimits}
\newcommand{\expp}[1]{\mathop {\mathrm{e}^{ #1}}}

\begin{document}

\includeversion{commentaries}

\title[Fragmentation associated to Lévy processes]{Fragmentation
  associated to Lévy processes using snake} 

\date{\today}
  
\author{Romain Abraham} 

\address{
MAPMO, Université d'Orléans,
B.P. 6759,
45067 Orléans cedex 2
FRANCE
}
  
\email{romain.abraham@univ-orleans.fr} 

\author{Jean-François Delmas}

\address{ENPC-CERMICS, 6-8 av. Blaise Pascal,
  Champs-sur-Marne, 77455 Marne La Vallée, France.}

\email{delmas@cermics.enpc.fr}

\begin{abstract}
  We  consider the height  process of  a Lévy  process with  no negative
  jumps, and  its associated continuous tree  representation.  Using Lévy
  snake  tools  developed  by  Duquesne and Le~Gall,  with  an  underlying
  Poisson process,  we construct a  fragmentation process, which  in the
  stable case corresponds to the self-similar fragmentation described by
  Miermont.  For the general  fragmentation process we
  compute a  family of dislocation  measures as well  as the law  of the
  size  of a tagged  fragment. We  also give  a special  Markov property
  for the snake which is interesting in itself.
  \end{abstract}



\maketitle


\section{Introduction}

We  present a  fragmentation  process associated  to general  critical
or sub-critical continuous
random  trees  (CRT)  which were  introduced  by  Le  Gall and  Le  Jan
\cite{lglj:bplpep}  and   developed  later  by  Duquesne   and  Le  Gall
\cite{dlg:rtlpsbp}.    This   extends   previous   work   from   Miermont
\cite{m:sfdfstsn} on  stable CRT. Although the underlying  ideas are the
same  in  both  constructions,  the  arguments in  the  proofs  are  very
different. Following Abraham and  Serlet \cite{as:psf} who deal with the
particular case of  Brownian  CRT, our arguments rely on Lévy Poisson
Snake   processes.   Those   path    processes   are   Lévy   Snake,   see
\cite{dlg:rtlpsbp},  with underlying   Poisson process.   To
prove the fragmentation property, we need some results
on Lévy Snake which are interesting by themselves. Eventually we give
the dislocation measure of the fragmentation process. We think this
construction provides non trivial examples of non self-similar
fragmentations, and that the tools developed here could give further results
on the fragmentation associated to CRT.

The next three subsections give a brief presentation of the
mathematical 
objects and state the mains results. The last one describes
the organization of the paper. 

\subsection{Exploration process} 

The  coding of  a tree  by its  height process  is now  well-known.  For
instance,  the  height  process   of  Aldous'  CRT  \cite{a:crt3}  is  a
normalized Brownian excursion. In \cite{lglj:bplpep}, Le Gall and Le Jan
associated to a Lévy process with no negative jumps that does not drift to
infinity, $X=(X_t, t\geq  0)$, a continuous state  branching process
(CSBP) and  a Lévy CRT which keeps  track of the genealogy  of the CSBP.
Let $\psi$ denote  the Laplace exponent of $X$.  We shall assume there
is no Brownian part, so that
\[
\psi(\lambda)=\alpha_0\lambda+\int_{(0,+\infty)}\pi(d\ell)
\left[\expp{-\lambda\ell}-1+\lambda\ell\right],  
\]
with  $\alpha_0\ge  0$  and  the   Lévy  measure  $\pi$  is  a  positive
$\sigma$-finite measure  on $(0,+\infty)$ such  that $\int_{(0,+\infty)}
(\ell\wedge \ell^2)\pi(d\ell)<\infty$.  Following \cite{dlg:rtlpsbp}, we
shall also assume that $X$  is of infinite variation a.s.  which implies
that  $\int_{(0,1)}\ell\pi(d\ell)=\infty$.  Notice those  hypothesis are
fulfilled in the stable case: $\psi(\lambda)=\lambda^\alpha$, $\alpha\in
(1,2)$.

Informally for the  height process, $H=(H_t, t\geq 0)$,  $H_t$ gives the
distance (which can be understood  as the number of generations) between
the individual labeled $t$ and the root, 0, of the CRT.  This process is
a key  tool in this  construction but it  is not a Markov  process.  The
so-called exploration process $\rho=(\rho_t,t\ge 0)$ is a càd-làg Markov
process taking  values in  $\cm_f(\R_+)$, the set  of finite  measure on
$\R_+$,  endowed with  the  topology of  weak  convergence.  The  height
process  can  easily  be  recovered  from  the  exploration  process  as
$H_t=H(\rho_t)$,  where  $H(\mu)$ denotes  the  supremum  of the  closed
support of the measure $\mu$ (with the convention that H(0)=0).  In some
sense  $\rho_t(dv)$  records the  ``number''  of  brothers, with  labels
larger than $t$, of the  ancestor of $t$ at generation $v$.  Furthermore
the jumps of $\rho$ are given by
\[
\rho_t-\rho_{t-} =\Delta_t \delta_{H_t},
\]
where $\Delta_t$  is the jump  of the Lévy  process $X$ at time  $t$ and
$\delta_x$ is the Dirac  mass at $x$.  Intuitively $\Delta_t$ represents
the  ``size''  of   the  progeny  of  such  individual   $t$. And the  set
$\big\{s\geq t; \min \{H_u,  u\in [t,s]\}\geq H_t\big \}$ represents the
``size''  of the  total descendants  of the  individual $t$.  
Such  individual $t$
corresponds to  a node in  the CRT.  To  each jump of $X$  corresponds a
node in the CRT and vice-versa.  Definition and properties of the height
process    and   exploration   process    are   recalled    in   Section
\ref{sec:levysnake}.

\subsection{Fragmentation}

A  fragmentation process  is a  Markov  process which  describes how  an
object with  given total  mass  evolves  as it breaks  into several
fragments randomly as time passes.  Notice there may be loss of mass but
no  creation.  This kind  of processes  has been  widely studied  in the
recent  years, see Bertoin~\cite{b:rfcp} and  references therein.  To be
more precise, the  state space of a fragmentation process  is the set of
the non-increasing sequences of masses with finite total mass
\[
\mathcal{S}^{\downarrow}=\left\{s=(s_1,s_2,\ldots); \;    s_1\ge    s_2\ge
  \cdots\ge                                         0\quad\text{and}\quad
  \Sigma(s)=\sum_{k=1}^{+\infty}s_k<+\infty\right\}.
\]
If we  denote by  $P_s$ the law  of a  $\cs^{\downarrow}$-valued process
$\Lambda     =(\Lambda     ^\theta,\theta\ge     0)$     starting     at
$s=(s_1,s_2,\ldots)\in\cs^{\downarrow}$,  we say  that $\Lambda  $  is a
fragmentation process if  it is a Markov process such that
$\theta\mapsto \Sigma(\Lambda^\theta)$ is non-increasing 
 and  if it fulfills the
fragmentation  property: the  law  of $(\Lambda  ^\theta,\theta \ge  0)$
under $P_s$ is the non-increasing reordering of the fragments of independent
processes of  respective laws $P_{(s_1,0,\ldots)}$,$P_{(s_2,0,\ldots)}$,
\ldots  In   other  words,  each  fragment   after  dislocation  behaves
independently  of the  others, and  its  evolution depends  only on  its
initial  mass.    As  a  consequence,   to  describe  the  law   of  the
fragmentation process  with any initial condition, it  suffices to study
the laws  $P_r:=P_{(r,0,\ldots)}$ for  any $r\in (0,+\infty)$,  i.e. the
law of the fragmentation process starting with a single mass $r$.

A fragmentation process is said to be self-similar of index $\alpha$ if,
for  any $r>0$, the  law of  the process  $(\Lambda ^\theta,\theta\ge 0)$
under $P_r$  is the law of the  process $(r\Lambda ^{r^\alpha \theta},\theta\ge
0)$  under $P_1$.  Bertoin~\cite{b:ssf}  proved that  the law  of a
self-similar   fragmentation   is  characterized   by:   the  index   of
self-similarity $\alpha$,  an erosion coefficient which  corresponds to a
deterministic  rate of  mass loss,  and a  dislocation measure  $\nu$ on
$\cs^{\downarrow}$ which describes sudden  dislocations of a fragment of
mass~1.

Connections between fragmentation processes and random trees or Brownian
excursion have been pointed out  by several authors.  Let us mention the
work of  Bertoin~\cite{b:fpcbm} who constructed  a fragmentation process
by  looking at  the  lengths of  the  excursions above  level  $t$ of  a
Brownian excursion. Aldous  and Pitman \cite{ap:sac} constructed another
fragmentation  process, which is  related  to  the additive  coalescent
process, by  cutting Aldous' Brownian CRT. Their proofs
rely  on  projective  limits   on  trees.   Those  results  have  been
generalized, by Miermont  \cite{m:sfdfstsh,m:sfdfstsn} to CRT associated
to  stable   Lévy  processes,  using  path
transformations  of  the  Lévy  process. Concerning  the  Aldous-Pitman's
fragmentation  process,   Abraham  and  Serlet   \cite{as:psf}  give  an
alternative  construction    using  Poisson   snakes. Our presentation
follow their ideas. However, we give next   a more
intuitive presentation which is in fact equivalent (see Section
\ref{sec:orf}). 

We  consider an  excursion of  the Lévy  process $X$  out of  $0$, which
correspond  also to  an excursion  of the  exploration process  (and the
height process)  out of  $0$. Let $\sigma$  denote the common  length of
those excursions.  Intuitively, $\sigma$  represents the ``size'' of the
total progeny of  the root $0$.  Let $\cj=\{t  \in [0,\sigma]; X_{t}\neq
X_{t-}\}$  the set  of  jumping times  of $X$  or  nodes of  the CRT,  and
consider $(T_t;  t\in \cj)  $ a countable  family of  independent random
variable such that $T_t$ is  distributed according to an exponential law
with parameter $\Delta_t$. At time  $T_t$, the node corresponding to the
jump $\Delta_t$  is cut from the  CRT.  Two individuals,  say $u\leq v$,
belongs to  the same fragment at time  $\theta$ if no node  has been cut
before time $\theta$ between them  and their most recent common ancestor
which is defined as $u\curlywedge  v= \inf \big\{t\in [0,u]; \min\{ H_r,
r\in [u,v]\}  =\min\{ H_r,  r\in [t,u]\}\big\}$.  Let  $\Lambda ^\theta$
denote  the  family  of  decreasing  positive Lebesgue  measure  of  the
fragments  completed by  zeros if  necessary so  that $\Lambda^\theta\in
\mathcal{S}^{\downarrow}$.  See  Section  \ref{sec:orf}  for  a  precise
construction.

\begin{figure}[ht]
\includegraphics[width=10cm]{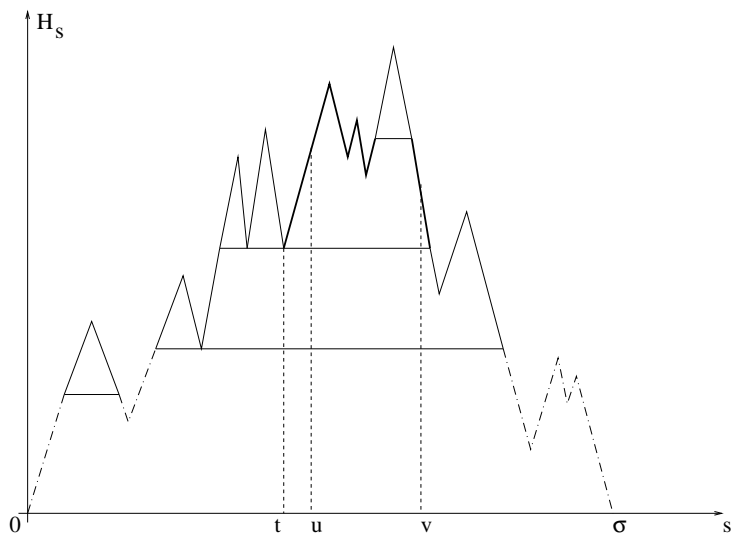}
\caption{Cutting at nodes: a modifier\label{fig:frag1}}
\end{figure}

Cutting nodes at time $\theta>0$  
may be viewed as adding horizontal lines under the epigraph of $H$ (see
figure \ref{fig:frag1}).

\begin{figure}[ht]
\includegraphics[width=10cm]{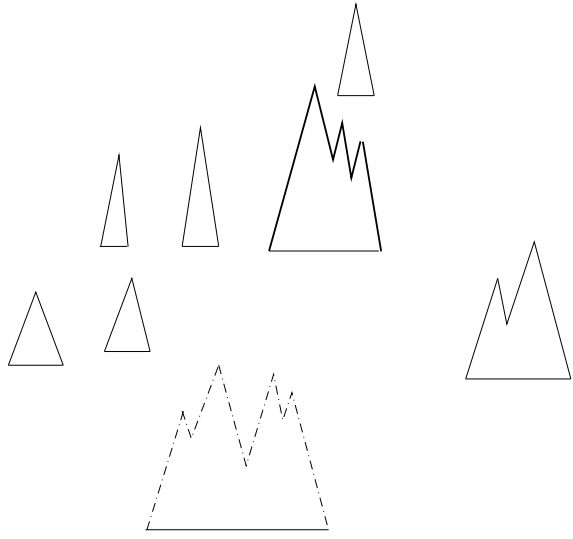}
\caption{Fragmentation of the excursion\label{fig:frag2}}
\end{figure}

We then consider the excursions obtained after cutting the initial
excursion along the horizontal lines and gluing together the
corresponding pieces of paths (for instance, the bold piece of the
path of $H$ in Figure \ref{fig:frag1} corresponds to the bold excursion
in Figure \ref{fig:frag2}). 
The lengths of these excursions, ranked in decreasing order, form the 
fragmentation process as $\theta$ increases. Of course, the figure are
caricatures as the process $H$ is very irregular and  the number of
fragments is  infinite.

Theorem \ref{theo:frag_property} asserts that the process
$(\Lambda ^\theta, \theta\geq 0)$ is a fragmentation process. There is  no loss
of mass thanks to Proposition  \ref{prop:tot-length}. 

In the stable
case, $\psi(\lambda)=\lambda^\alpha$ with  $\alpha\in (1,2)$, using
scaling properties, we get the fragmentation is self-similar with index
$1/\alpha$  and we
recover the results of Miermont \cite{m:sfdfstsn}, see Corollary
\ref{cor:alpha-frag}. In particular the dislocation measure is given by:
for any measurable non-negative
function $F$ on $\cs^\downarrow$,
\[
\int
F(x)\nu(dx)=\frac{\alpha(\alpha-1)\Gamma\bigl([\alpha-1]/\alpha\bigr)}{\Gamma(2-\alpha)}\E\left[S_1F(\Delta
  S_t/S_1,\ t\le 1)\right],
\]
where  $(S_t,t\ge 0)$  is a  stable subordinator  with  Laplace exponent
$\psi^{-1}(\lambda)=  \lambda^{1/\alpha}$, and $F(\Delta  S_t/S_1,\ t\le
1)$ has to  be understood as $F$ applied  to the decreasing reordering
of the sequence $(\Delta S_t/S_1,\ t\le 1)$.

In the general case, the fragmentation is not self-similar. However, if
$\ct=\{\theta\geq 0; \Lambda^\theta\neq \Lambda^{\theta-}\}$ denotes the
jumping times of the process $\Lambda $, we get as a direct consequence
of  Section
\ref{sec:comp-int} that 
\[
\sum_{\theta\in \ct} \delta _{\theta, \Lambda ^\theta} 
\]
is  a   point  process  with  intensity   $d\theta  \tilde  \nu_{\Lambda
  ^{\theta-}}     (ds)$,     where     $(\tilde     \nu_x,     x     \in
\mathcal{S}^{\downarrow})$  is a  family of  $\sigma$-finite  measures on
$\mathcal{S}^{\downarrow}$. There  exists  a  family  $(\nu_r,  r>0)$  of
$\sigma$-finite  measure on  $\mathcal{S}^{\downarrow}$,  which we  call
dislocation measures of the fragmentation  $\Lambda $, such that for any
$x=(x_1,x_2,  \ldots)\in \mathcal{S}^{\downarrow}$  and  any non-negative
measurable function, $F$, defined on $\cs^\downarrow$,
\[
\int F(s) \tilde \nu_x (ds)= \sum_{i\in \N^*; x_i>0} \int F(x^{i,s})
\nu_{x_i} (ds),
\]
where  $x^{i,s}$ is  the decreasing  reordering of  the  merging of
the sequences $s\in
\cs^\downarrow$ and  $x$, where $x_i$  has been removed of  the sequence
$x$.   This  means  that only  one  element  of  $x$ fragments  and  the
fragmentation depends  only on  the size of  this very  fragment.  The
dislocation 
measures  can  be  computed,  see  Theorem  \ref{th:calcul_nu_r}.   In
particular $\nu_r(dx)$-a.e. $\sum_{i\in \N^*} x_i=r$ assures there is no
loss  of mass  at the  dislocation.  The definition  of the  dislocation
measures  is more  involved than  in the stable  case. However,  it  can
still be
written  using the  law of  the jumps  of a  subordinator  with Laplace
exponent $\psi^{-1}$.

\subsection{The pruned exploration process}

In fact the dislocation measure is  computed by studying the evolution of a
tagged fragment, for example the one that contains the root of the CRT. 
Therefore, it is natural to consider first the exploration process of the
fragment containing the root at time $\theta$. 
The pruned exploration
process, $\tilde \rho=(\tilde \rho_t, t\geq 0)$, is defined by $\tilde
\rho_t=\rho_{C_t}$, where $C_t=\inf \{  r>0; 
A_r \geq t\}$ is the right continuous inverse of $A_t$, the Lebesgue
measure of the set of individuals prior to $t$ who belongs to the tagged
fragment at time $\theta$ (Section \ref{sec:pruned_process}). The
pruned process $\tilde \rho$ corresponds to the exploration process
associated to the dashed height process of Figures \ref{fig:frag1} and
\ref{fig:frag2}. To get
the law of the pruned exploration process (Section
\ref{sec:law_pruned}),    we use a Poisson Lévy snake
approach (Section \ref{sec:LPS})  and we prove a special Markov
property, Theorem \ref{th:SMP} in Section \ref{sec:Markov_special},
which is of independent interest. Notice this theorem differs from
Proposition 4.2.3 in \cite{dlg:rtlpsbp},  or Proposition 7 in
\cite{blglj:sbps}, where in both cases the exit measure is singular,
whereas here it is absolutely continuous w.r.t. to the Lebesgue measure.

Eventually,  using   martingales,  we   get  Theorem
\ref{thm:law_pruned}:  the  pruned  exploration
process  $\tilde \rho$ is  the exploration  process associated  to a
Lévy process, $X^{(\theta)}$,  with Laplace exponent $\psi^{(\theta)}$
defined by: for $\lambda\in \R_+$,
\[
\psi^{(\theta)}
(\lambda)=\psi(\lambda+\theta)-\psi(\theta).
\]
There exists other pruning procedure for Galton-Watson trees, see for
example \cite{dw:glt} and references therein.

Notice that conditionally on the length of the excursion, the excursions
of $X$ and $X^{(\theta)}$ out  of $0$ are equally distributed (see Lemma
\ref{lem:loi_cond}).  This  property,  as  well as  the  special  Markov
property  are essential to  prove the  fragmentation property.   We also
compute, see Proposition \ref{prop:s-sq}  the joint law of $\sigma$, the
initial mass of the fragment, and $\tilde \sigma$ the mass of the tagged
fragment at time $\theta$, under the excursion measure.

\subsection{Organization of the paper}

In Section  \ref{sec:levysnake}, we recall the construction  of the Lévy
CRT  and  give the  properties  we shall  use  in  this paper.   Section
\ref{sec:LPS} is  devoted to the  definition and some properties  of the
Lévy  Poisson snake.   {F}rom  this  Lévy Poisson  snake,  we define  in
Section  \ref{sec:pruned_process} the  pruned exploration  process which
corresponds to the tagged fragment  that contains 0.  Then, we introduce
in Section  \ref{sec:Markov_special} a  special Markov property  for the
Lévy Poisson snake:  Theorem \ref{th:SMP} and Corollary \ref{cor:super}.
We  compute  in  Section  \ref{sec:law_pruned}  the law  of  the  pruned
exploration   process,   see   Theorem  \ref{thm:law_pruned}.    Section
\ref{sec:prop_pruned} is then devoted to the study of some properties of
the pruned exploration process under the excursion measure.  Eventually,
we  construct  in  Section  \ref{sec:LS-F},  the  fragmentation  process
associated  to  our  Lévy  Poisson  snake and  prove  the  fragmentation
property, Theorem  \ref{theo:frag_property}, and check there  is no loss
of     mass,    Proposition    \ref{prop:tot-length}.      In    Section
\ref{sec:disloc}, we  identify completely  the law of  the fragmentation
process    by    computing    the    dislocation    measures,    Theorem
\ref{th:calcul_nu_r},   and   we   recover   the  result   of   Miermont
\cite{m:sfdfstsn} for the stable case in Corollary \ref{cor:alpha-frag}.


\section{Lévy snake: notations and properties}\label{sec:levysnake}

We  recall here  the construction  of the  Lévy continuous  random tree
(CRT) introduced in \cite{lglj:bplpep,lglj:bplplfss} and developed later
in \cite{dlg:rtlpsbp}. We  will emphasize on the height  process and the
exploration process which are the key tools to handle this tree. The
results of this section are mainly extract from \cite{dlg:rtlpsbp}.

\subsection{The underlying Lévy process}\label{subsec:levy}
We consider  a $\R$-valued Lévy process $(X_t,t\ge
0)$ with no negative jumps, starting from 0. Its law is characterized by
its Laplace transform: for  $\lambda\ge 0$ 
\[
\E\left[\expp{-\lambda X_t}\right]=\expp{t\psi(\lambda)},
\]
where its Laplace exponent, $\psi$, is given by 
\[
\psi(\lambda)=\alpha_0\lambda+\int_{(0,+\infty)}\pi(d\ell)
\left[\expp{-\lambda\ell}-1+\lambda\ell\right],  
\]
with $\alpha_0\ge 0$ and the Lévy measure $\pi$ is a positive
$\sigma$-finite measure on $(0,+\infty)$ such that
\begin{equation}
   \label{eq:cond_pi}
\int_{(0,+\infty)} (\ell\wedge
\ell^2)\pi(d\ell)<\infty\quad\mbox{and}
\quad\int_{(0,1)}\ell\pi(d\ell)=\infty.
\end{equation}
The first  assumption (with the  condition $\alpha_0\ge 0$)  implies the
process $X$ does not drift to infinity, while the second implies $X$ is
of infinite variation a.s.

For  $\lambda\geq   1/\varepsilon>0$,  we  have  $\expp{-\lambda\ell}-1+
\lambda\ell\geq     \frac{1}{2}\lambda     \ell    \ind_{\{     \ell\geq
  2\varepsilon\}}$, which implies  that $\lambda^{-1} \psi(\lambda) \geq
\alpha_0+  \int_{(2\varepsilon,\infty  )} \ell \; \pi(d\ell)$. We
deduce that
\begin{equation}
   \label{eq:psi/l}
\lim_{\lambda
\rightarrow\infty } \frac{\lambda}{\psi(\lambda)} =0.
\end{equation}

We introduce some processes related to $X$.  Let $\cj=\{s\geq 0; X_s\neq
X_{s-}\}$, the set of jumping times of $X$. For $s\in \cj$, we denote by
\[
\Delta_s= X_s- X_{s-}
\]
the  jump of $X$  at time  $s$ and  $\Delta_s=0$ otherwise.   The random
measure  $ \cx=\sum_{s\in\cj} \delta_{s,  \Delta_s}$  is a
Poisson point process with intensity $\pi(d\ell)$.  Let $I=(I_t,t\ge 0)$
be  the infimum  process of  $X$, $I_t=\inf_{0\le  s\le t}X_s$,  and let
$S=(S_t,t\ge 0)$  be the supremum process,  $S_t=\sup_{0\le s\le t}X_s$.
We will  also consider for every $0\le  s\le t$ the infimum  of $X$ over
$[s,t]$:
\[
I_t^s=\inf_{s\le r\le t}X_r.
\]

The point 0 is regular for the Markov 
process $X-I$, and $-I$ is the  local time of $X-I$ at 0 (see
\cite{b:pl}, chap. VII). Let $\N$ be the associated excursion measure of
the process $X-I$ out of 0, and $\sigma=\inf\{t>0; X_t-I_t=0\}$ the length
of the excursion of $X-I$ under $\N$. We will assume that under $\N$,
$X_0=I_0=0$. 

Since $X$ is of infinite variation, 0 is also regular for the Markov 
process $S-X$. The local time, $L=(L_t, t\geq 0)$, of $S-X$ at 0 will be
normalized so that 
\[
\E[\expp{-\beta S_{L^{-1}_t}}]= \expp{- t \psi(\beta)/\beta},
\]
where $L^{-1}_t=\inf\{ s\geq 0; L_s\geq t\}$ (see also \cite{b:pl}
Theorem VII.4 (ii)).

\subsection{The height process and the Lévy CRT}
For each $t\geq 0$, we consider the reversed process at time $t$,
$\hat X^{(t)}=(\hat X^{(t)}_s,0\le s\le t)$ by:
\[
\hat X^{(t)}_s=
X_t-X_{(t-s)-} \quad \mbox{if}\quad  0\le s<t,
\]
and $\hat X^{(t)}_t=X_t$. The two processes $(\hat X^{(t)}_s,0\le s\le t)$
and $(X_s,0\le s\le t)$ have the same law. Let $\hat S^{(t)}$ be the
supremum process of $\hat X^{(t)}$ and $\hat L^ {(t)}$ be the
local time at $0$ of $\hat S^{(t)} - \hat X^{(t)}$ with the same
normalization as $L$. 

\begin{defi}
   There exists a process $H=(H_t, t\geq 0)$, called the
   height process, such that for all $t\geq 0$, a.s. $H_t=\hat
   L^{(t)}_t$, and $H_0=0$. Furthermore $H$ is lower semi-continuous
   a.s. and a.s. for all $t'>t\geq 0$, the process $H$ takes all the
   values between $H_t$ and $H_{t'}$ on the time interval $[t,t']$. 
\end{defi}

The  height process  $(H_t, t\in  [0,\sigma])$  under $\N$  codes  a
continuous  genealogical structure,  the Lévy  CRT, via  the following
procedure. 

\begin{itemize}
   \item[(i)] To each $t\in [0,\sigma]$ corresponds a vertex at
   generation $H_t$. 
   \item[(ii)] Vertex $t$ is an ancestor of vertex $t'$ if
   $H_t=H_{t,t'}$, where 
\begin{equation}
   \label{eq:def_H}
H_{t,t'}=\inf\{H_u, u\in [t\wedge
t', t\vee t']\}.
\end{equation}
In general $H_{t,t'}$ is the generation of the last
   common ancestor to $t$ and $t'$. 
   \item[(iii)] We put $d(t,t')=H_t+H_{t'}- 2 H_{t,t'}$ and identify $t$
   and $t'$ ($t\sim t'$) if $d(t,t')=0$. 
\end{itemize}

The Lévy CRT coded by $H$ is then the quotient set $[0,\sigma]/ \sim$,
equipped with the distance $d$ and the genealogical relation specified
in (ii).

\subsection{The exploration process}
\label{sec:PLRT}
The height process is in general not Markov. But it is a very simple
function of a measure-valued Markov process, the so-called exploration
process.

If $E$  is a polish  space, let $\cb(E)$  (resp. $\cb_+(E)$) be the  set of
real-valued  measurable (resp. and non-negative) functions  defined on $E$
endowed    with   its   Borel    $\sigma$-field,   and    let   $\cm(E)$
(resp.  $\cm_f(E)$)  be  the  set  of  $\sigma$-finite  (resp.   finite)
measures  on  $E$, endowed  with  the  topology  of vague  (resp.  weak)
convergence.  For any measure $\mu\in\cm(E)$ and $f\in \cb_+(E)$, we write
$$\langle \mu,f\rangle =\int f(x)\,\mu(dx).$$

The     exploration    process     $\rho=(\rho_t,t\ge    0)$     is    a
$\cm_f(\R_+)$-valued process defined  as follows: for  every $f\in
\cb_+(\R_+) $,
$$\langle \rho_t,f\rangle =\int_{[0,t]} d_sI_t^sf(H_s),$$
or equivalently 
\begin{equation}\label{eq:def_rho}
\rho_t(dr)=\sum_{\stackrel{0<s\le t}
  {X_{s-}<I_t^s}}(I_t^s-X_{s-})\delta_{H_s}(dr). 
\end{equation}
In particular, the total mass of $\rho_t$ is
$\langle \rho_t,1\rangle =X_t-I_t$.

For $\mu\in \cm(\R_+)$, we set 
\begin{equation}
   \label{def:H}
H(\mu)=\sup\, \Supp \mu,
\end{equation}
where $ \Supp \mu$ is the closed support of $\mu$,  with the
convention $H(0)=0$. We have
\begin{prop}
Almost surely, for every $t>0$,
\begin{itemize}
\item $H(\rho_t)=H_t$,
\item $\rho_t=0$ if and only if $H_t=0$,
\item if $\rho_t\neq 0$, then $\Supp \rho_t=[0,H_t]$. 
\item $\rho_t= \rho_{t^-} + \Delta_t \delta_{H_t}$, where $\Delta_t=0$
  if $t\not\in \cj$. 
\end{itemize}
\end{prop}

In the definition  of the exploration process, as $X$  starts from 0, we
have  $\rho_0=0$ a.s.  To  get the  Markov property  of $\rho$,  we must
define  the  process  $\rho$  started  at any  initial  measure  $\mu\in
\cm_f(\R_+)$.

For $a\in  [0, \langle  \mu,1\rangle ] $,  we define the  erased measure
$k_a\mu$ by
\[
k_a\mu([0,r])=\mu([0,r])\wedge (\langle \mu,1\rangle -a), \quad \text{for $r\geq 0$}.
\]
If $a> \langle  \mu,1\rangle $, we set $k_a\mu=0$.   In other words, the
measure $k_a\mu$ is the measure $\mu$ erased by a mass $a$ backward from
$H(\mu)$.

For $\nu,\mu \in \cm_f(\R_+)$, and $\mu$ with compact support, we
  define the concatenation $[\mu,\nu]\in \cm_f(\R_+) $ of the
  two measures by: 
\[
\bigl\langle [\mu,\nu],f\bigr\rangle =\bigl\langle \mu,f\bigr\rangle +\bigl\langle \nu,f(H(\mu)+\cdot)\bigr\rangle ,
\quad f\in \cb_+(\R_+).
\]

Eventually,  we  set for  every  $\mu\in  \cm_f(\R_+)$  and every  $t>0$
$\rho_t^\mu=\bigl[k_{-I_t}\mu,\rho_t]$.  We say that $(\rho^\mu_t, t\geq
0)$  is  the  process  $\rho$  started at  $\rho_0^\mu=\mu$,  and  write
$\P_\mu$  for its  law. Unless  there is  an ambiguity,  we  shall write
$\rho_t$ for $\rho^\mu_t$.

\begin{prop}
The process $(\rho_t,t\ge 0)$ is a càd-làg strong Markov process in
$\cm_f(\R_+)$.
\end{prop}

\begin{rem}
  \label{rm:rho-L} 

{F}rom the construction of $\rho$, we get that a.s.
  $\rho_t=0$ if and only if $ -I_t\geq {\langle \rho_0,1\rangle }$ and $X_t-I_t=0$.  This implies that $0$ is also a regular point for $\rho$. Let $(\tau_s, s\geq 0)$ be the right continuous inverse of $-I$: $\tau_s=\inf \{ t>0; -I_t>s\}$. We get the local time at $0$ of $\rho^\mu$, $(L^0_t, t\geq 0)$, is given for $t\geq 0$, by
\[
L^0_t=-I_t+I_{t\wedge \tau_{\langle \mu,1\rangle }}.
\]
Notice that $\N$ is also the excursion measure of the process $\rho$ out
of  $0$, and that  $\sigma$, the  length of  the excursion,  is $\N$-a.e.
equal to $\inf\{ t>0; \rho_t=0\}$.
\end{rem}

\begin{rem}   
\label{rm:recons}
Recall  $(\Delta_s, s\geq 0)$  are the  jumps of  the process  $X$.  The
process $\rho$  is adapted to  the filtration generated by  the process
$X$, that is by the Poisson point process $\cx$, and $\rho_0$, completed
the usual way.   {F}rom the construction of $\rho$, we  get there exists a
measurable  function,  $\Gamma$,  defined on  $\cm(\R_+^2)\times
\cm_f( \R_+)$  (endowed with its Borel $\sigma$-field)
taking  values   in  $\cm_f(\R_+)$  (endowed  with   its  Borel
$\sigma$-field), such that
\[
\rho_t=\Gamma( \cx\ind_{[0,t]\times \R_+}, \rho_0).
\]
On the other hand, notice that a.s.  the jumping times of $\rho$ are also
the   jumping   times   of   $X$,   and   for   $s\in   \cj$,   we   have
$\rho_s(\{H_s\})=\Delta_s$.  We  deduce that $(\Delta_u,u\in  (s,t])$ is
measurable w.r.t.  the $\sigma$-field $\sigma(\rho_u,u\in [s,t])$.

\end{rem}

\subsection{The dual process and representation formula}
\label{sec:dual}

We  shall need the  $\cm_f(\R_+)$-valued process  $\eta=(\eta_t,t\ge 0)$
 defined by
\[
\eta_t(dr)=\sum_{\stackrel{0<s\le t}{X_{s-}<I_t^s}}(X_s-I_t^s)\delta
_{H_s}(dr).
\]
The process $\eta$ is the dual process of $\rho$ under $\N$ (see
Corollary 3.1.6 in \cite{dlg:rtlpsbp}). We write  (recall
$\Delta_s= X_s-X_{s-}$)
\begin{equation}
   \label{eq:def-kappa}
\kappa_t(dr)=\rho_t(dr)+\eta_t(dr)=\sum_{\stackrel{0<s\le
    t}{X_{s-}<I_s^t}}\Delta_s \delta _{H_s}(dr).
\end{equation}

We recall the Poisson representation of $(\rho,\eta)$ under $\N$. Let
$\mathcal{N}(dx\,   d\ell\,  du)$   be  a   Poisson  point   measure  on
$[0,+\infty)^3$ with intensity
$$dx\,\ell\pi(d\ell)\ind_{[0,1]}(u)du.$$
For every $a>0$, let us denote by $\mathbb{M}_a$ the law of the pair
$(\mu_a,\nu_a)$ of finite measures on $\R_+$ defined by:  for $f\in \cb_+(\R_+)$ 
\begin{align}
\label{def:mu_a}
\langle \mu_a,f\rangle  & =\int\mathcal{N}(dx\, d\ell\, du)\ind_{[0,a]}(x)u\ell f(x),\\
\label{def:nu_a}
\langle \nu_a,f\rangle  & =\int\mathcal{N}(dx\, d\ell\, du)\ind_{[0,a]}(x)\ell(1-u)f(x).
\end{align}
We eventually set $\mathbb{M}=\int_0^{+\infty}da\, \expp{-\alpha_0
  a}\mathbb{M}_a$. 

\begin{prop}\label{prop:poisson_representation1}
For every non-negative measurable function $F$ on $\cm_f(\R_+)^2$,
\[
\N\left[\int_0^\sigma F(\rho_t, \eta_t) \; dt \right]=\int\mathbb{M}(d\mu\,
    d\nu)F (\mu, \nu) ,
\]
where $\sigma=\inf\{s>0; \rho_s=0\}$ denotes the length of the
    excursion. 
\end{prop}
We recall Lemma 3.2.2 from \cite{dlg:rtlpsbp}, we shall use later.
\begin{prop}\label{prop:poisson_representation}
Let $\tau$ be an exponential variable of parameter $\lambda>0$
independent of $X$ defined under the measure $\N$. Then,
for every $F \in \cb_+(\cm_f(\R_+))$, we have 
$$\N\left(F(\rho_\tau)\ind_{\tau\le
    \sigma}\right)=\lambda\int\mathbb{M}(d\mu\,
    d\nu)F(\mu)\expp{-\psi^{-1}(\lambda)\langle \nu,1\rangle }.$$
\end{prop}
It is easy to deduce from this (see also the beginning of Section
3.2.2. \cite{dlg:rtlpsbp}) that for $\lambda>0$
\begin{equation}
   \label{eq:N_s}
\N\left[1 -\expp{-\lambda
  \sigma}\right] =\psi^{-1}(\lambda). 
\end{equation}

\section{The Lévy Poisson snake}
\label{sec:LPS}
As in  \cite{as:psf}, we want to  construct a Poisson snake  in order to
cut the Lévy CRT at its nodes.  For this, we will construct
a  consistent family  $(m^\theta=(m^\theta_t,t\ge 0),  \theta\ge  0)$ of
measure-valued processes. For fixed $\theta$ and $t$, $m^\theta$ will be
a point-measure whose  atoms mark the atoms of  the measure $\rho_t$ and
such that  the set  of atoms of  $m^{\theta+\theta'}$ contains  those of
$m^\theta$. To  achieve this, we  attach to each  jump of $X$  a Poisson
process indexed by $\theta$, with  intensity equal to this jump. In fact
only the first jump of the  Poisson processes will be necessary to build
the fragmentation process.

\subsection{Definition and properties}

Conditionally on
$\cx=\sum_{s>  0}  \delta_{s, \Delta_s}$,  we consider a  family
$(\sum_{u>0}\delta_{V_{s,u}}, s\in \cj)$ of independent 
Poisson  point measures  on
$\R_+$  with  respective intensity $ \Delta_s \;\ind_{\{u>0\}}  du$.   We   define  the  $\cm(\R_+^2)$-valued
process $M=(M_t,t\geq 0)$ by
\begin{equation}
   \label{def:mt}
M_t(dr,dv)=\sum_{\stackrel{0<s\le t}{X_{s-}<I_t^s}} \; (I^s_t-X_{s-})
(\sum_{u>0} \delta_{V_{s,u}}(dv))  \; \delta_{H_s}(dr).
\end{equation}

\begin{rem}
 The additional coefficient $I_t^s-X_{s-}$ is not very
 important   and  is   only   needed   for  the   process   $M$  to   be
 right-continuous.
\end{rem}

Let $\theta>0$. For $t\geq 0$, notice that 
\[
M_t(\R_+\times [0,\theta]) \leq  \sum_{0<s\leq t} \Delta_s \xi_s,
\]
with $\xi_s=\Card \{ u>0; V_{s,u} \leq \theta \}$. 
In particular, we
have for $T>0$, 
\begin{equation}
   \label{eq:majo_M}
\sup_{t\in [0,T]} M_t(\R_+\times [0,\theta]) \leq  \sum_{0<s\leq
  T} \Delta_s \xi_s.
\end{equation}
Notice the variable
$\xi_s$ are, conditionally on $\cx$,  independent and distributed as
Poisson random variables with parameter $\theta\Delta_s$. We have 
  $\E  [  \sum_{0<s\leq   T}  \Delta_s   \xi_s  |   \cx]  =\theta
  \sum_{0<s\leq  T} \Delta_s^2$. As  $\int_{(0,\infty )}  (\ell^2 \wedge
  \ell ) \pi(d\ell)$ is finite, this implies the quantity $\sum_{0<s\leq
  T} \Delta_s^2$ is finite a.s. In particular we have a.s. 
\[
\sup_{t\in [0,T]} M_t(\R_+\times [0,\theta]) <\infty ,
\]
and $M_t$ is a $\sigma$-finite measure on $ \R_+^2$. Notice that a.s. 
\begin{equation}
   \label{eq:decomp-M}
M_t(dr,dv)=\rho_t(dr) M_{t,r}(dv),
\end{equation}
where $M_{t,r}$ is a $\sigma$-finite counting measure on $\R_+$.

We call the process $\cs=((\rho_t,M_t), t\geq 0)$ the Lévy Poisson snake
started at  $\rho_0=0, M_0=0$.  To get the  Markov property of  the Lévy
Poisson snake, we  must define the process $\cs$  started at any initial
value $(\mu, \Pi)\in \S$, where $\S$ is the set of pair $(\mu,\Pi)$ such
that  $\mu\in \cm_f(\R_+)$   and $\Pi(dr,dv)=\mu(dr)
\Pi_r(dv)$, $\Pi_r$ being $\sigma$-finite  measures on $\R_+$, such that
$\Pi(\R_+\times  [0,\theta])<\infty $  for all  $\theta\geq 0$.   We set
$H^\mu_t=H(k_{-I_t}    \mu)$.     Then,    we   define    the    process
$M^{\mu,\Pi}=(M^{\mu,\Pi}_t,    t\geq    0)$    by:   for    $\varphi\in
\cb_+(\R_+^2)$,
\[
\langle M^{\mu,\Pi}_t, \varphi\rangle = \int_{(0,\infty )}  \varphi(r,v) k_{-I_t}
\mu(dr) \Pi_r(dv) 
+ \int_{(0,\infty )}  \varphi( r+ H^\mu_t, v) M_t(dr, dv).
\]
We shall write $M$ for $M^{\mu, \Pi}$. By construction 
and since $\rho$
is an homogeneous Markov process, the Lévy
Poisson snake $\cs=(\rho,M)$ is an homogeneous  Markov process. 

We  now denote by  $\P_{\mu,\Pi}$ the  law of  the Lévy  Poisson snake
starting at 0  from $(\mu,\Pi)$, and by $\P_{\mu,\Pi}^*$  the law of the
Lévy  Poisson snake  killed when  $\rho$ reaches  $0$. We  deduce from
\reff{eq:majo_M}, that a.s. 
\begin{equation}
   \label{eq:majo_M2}
\E_{\mu,\Pi} \left[\sup_{t\in [0,T]} M_t(\R_+\times [0,\theta])\Bigm|\cx\right]
  \leq  \theta \sum_{0<s\leq 
  T} \Delta_s^2 +\Pi(\R_+\times [0,\theta]) <\infty .
\end{equation}

Let $\cf=(\cf_t, t\geq 0)$ be the filtration generated by $\cs$
completed the usual way. Notice this
filtration is also generated by the processes $(\cx([0,t], \cdot), t\geq
0)$ and $(\sum_{s\in \cj, \, s\leq t} \sum_{  u\geq   0 }   \delta_{V_{s,u}},  t\geq
0)$. In particular the filtration $\cf$ is right continuous. And by
construction, we have that $\rho$ is Markovian with respect to $\cf$. 

\begin{prop}
The Lévy Poisson snake, $\cs$,  is  a càd-làg strong Markov process in
$\S\subset \cm_f(\R_+)\times \cm(\R_+^2)$.  
\end{prop}

\begin{proof}
   
We first check the process $M$ is right continuous. Recall
\reff{eq:decomp-M}. 
We have by construction a.s. for all $t'>t$, 
\[
M_{t'} (dr,  dv)
= k_{-I_{t'}^t} \rho_t(dr) M_{t,r}(dv)
+ \rho_{t'}(dr) \ind_{\{r> H_{t,t'}\}}
M_{t',r}(dv),
\]
where $H_{t,t'}$ is defined by \reff{eq:def_H}.
Thanks to \reff{eq:majo_M}, we have, for $\theta>0$, 
\[
\int_{\R_+} \rho_{t'}(dr)
\ind_{\{r> H_{t,t'}\}}   M_{t',r}([0,\theta])\leq
\sum_{t<s\leq  t'}  \Delta_s  \xi_s   .
\]
In particular this quantity decreases  to 0 as $t'\downarrow t$ a.s.  By
the  properties of  the  exploration  process, we  recall  that a.s.   $
k_{-I_{t'}^t}  \rho_t=  \rho   _{t"}$,  where  $t"=\inf\{  s\in  [t,t'];
I^t_s=I^t_{t'}\}$.
{F}rom the right continuity of $\rho$,
we deduce that a.s.  for the vague convergence
\[
 \lim_{t'\downarrow t} M_{t'} =M_t. 
\]
This implies the right continuity of the process $M$ for the vague
topology on $\cm(\R_+^2)$. 

Now, we check the process $M$ has left limits. Let $t<t'$. For $r\in [0, H_{t,t'}]$,
we  have  $k_{-I_{t'}^t}  \rho_t(dr)  M_{t,r}=\ind_{\{r\leq  H_{t,t'}\}}
\rho_{t'}(dr) M_{t',r}$, as well as
\[
M_{t} (dr,  dv)
= \ind_{\{r\leq  H_{t,t'}\}}
 \rho_{t'}(dr) M_{t',r}(dv)
+ [\rho_{t}(dr) - k_{-I_{t'}^t}
\rho_t(dr)] M_{t,r}(dv).
\]
If $\rho$ is continuous  at $t'$, then either $\rho_{t'}(\{H_t'\})=0$ or
$H_{t,t'}= H_{t'}$ for $t$ close  enough to $t'$. In particular, since
$\lim_{t\rightarrow t'} H_{t,t'}=H_{t'}$, we have
$\lim_{t\uparrow   t'}   \ind_{\{r\leq   H_{t,t'}\}}   \rho_{t'}(dr)   =
\rho_{t'}(dr)$.  If $\rho$ is not  continuous at $t'$, this implies that
$\rho_{t'} (dr)=\rho_{t'-} (dr)  + \Delta_{t'} \delta_{H_{t'}} (dr)$ and
for  $t$  close  enough  to  $t'$,  $  H_{t,t'}<H_{t'}$.  Then,  we  get
$\lim_{t\uparrow   t'}   \ind_{\{r\leq   H_{t,t'}\}}
\rho_{t'}(dr) =  \rho_{t'-}(dr)$.  In  any case, we  have a.s.   for the
vague convergence
\[
 \lim_{t\uparrow t'} \ind_{\{r\leq H_{t,t'}\}}
 \rho_{t'}(dr) M_{t',r}(dv) =  \rho_{t'-}(dr) M_{t',r}(dv).
\]
Now, we check that for the vague topology
\[
 \lim_{t\uparrow t'} [\rho_{t}(dr) - k_{-I_{t'}^t}
\rho_t(dr)] M_{t,r}(dv) =0 .
\]
For this purpose, we remark that
\begin{align*}
   \E_{\mu,\Pi} \left[\int_{\R_+}\!\! [\rho_{t}(dr) - k_{-I_{t'}^t}
\rho_t(dr)] M_{t,r}([0,\theta])| \cx\right]
& = \theta \int_{\R_+}\!\! [\rho_{t}(dr) - k_{-I_{t'}^t}
\rho_t(dr)] (\rho_t(\{r\})+ \eta_t(\{r\}) )\\ 
& \leq   \theta  (\langle \rho_t+ \eta_t,1\rangle  ) \int_{\R_+} [\rho_{t}(dr) -
k_{-I_{t'}^t} 
\rho_t(dr)]\\ 
&= \theta  (\langle \rho_t+ \eta_t,1\rangle  ) (-I_{t'}^t).
\end{align*}
As $\rho$ and $\eta$ are respectively càd-làg and càg-làd
process, they are bounded over any finite interval a.s.
Since $\lim_{t\uparrow t'}I_{t'}^t=0$, we deduce that 
\[
\lim_{t\uparrow t'} \E_{\mu,\Pi} \left[\int_{\R_+} [\rho_{t}(dr) -
k_{-I_{t'}^t} 
\rho_t(dr)] M_{t,r}([0,\theta])| \cx\right] =0.
\]
Thanks to \reff{eq:majo_M2} and Fatou's Lemma, we deduce that 
\[
\lim_{t\uparrow t'} \int_{\R_+} [\rho_{t}(dr) -
k_{-I_{t'}^t} 
\rho_t(dr)] M_{t,r}([0,\theta]) =0.
\]
Therefore, we conclude that for vague topology,
\[
\lim_{t\uparrow t'}  M_t =M_{t'-}.
\]

We deduce that for the vague topology on $\cm(\R_+^2)$, the process $M$
is a.s. càd-làg. This implies the process $\cs$ is a.s. càd-làg. 

We check the strong Markov property of $\cs$. 
Mimicking the proof of Proposition 1.2.3 in \cite{dlg:rtlpsbp}, and
using properties of Poisson point measure, one gets that, for any $\cf$-stopping time $T$, we have a.s. for every $t>0$,
\begin{align*}
\rho_{T+t} & =\left[k_{-I_t^{(T)}}\rho_T,\rho_t^{(T)}\right]\\
M_{T+t}(dr,dv) & =k_{-I_t^{(T)}}\rho_t^{(T)}(dr)M_{T,r}(dv)+M_t^{(T)}(dr+H(k_{-I_t^{(T)}}\rho_T),dv)
\end{align*}
where $I^{(T)}$,  $\rho^{(T)}$ and $M^{(T)}$  are the analogues  of $I$,
$\rho$   and   $M$   with   $X$   replaced  by   the   shifted   process
$X^{(T)}=(X_{T+t}-X_T,t\ge 0)$. This implies the strong Markov property.
\end{proof}

\subsection{Poisson representation of the snake}

Notice  that   a.s.  $(\rho_t,M_t)=(0,0)$   if  and   only  if
$\rho_t=0$. In  particular, $(0,0)$  is a regular  point for  the Lévy
Poisson snake, with associated local  time $(L^0_s, s\geq 0)$.  We still
write $\N$ for the excursion measure  of the Lévy Poisson snake out of
$(0,0)$, with the same normalization as in Section \ref{sec:dual}.

We  decompose the  path of  $\cs$ under  $\P^*_{\mu, \Pi}$  according to
excursions of  the total  mass of $\rho$  above its minimum,  see Section
4.2.3  in \cite{dlg:rtlpsbp}.  More  precisely let  $(\alpha_i,\beta_i),
i\in I$ be the excursion intervals of the process $\langle \rho,1\rangle
$  above its minimum  under $\P^*_{\mu,  \Pi}$. For  every $i\in  I$, we
define $h_i=H_{\alpha_i}$ and $\cs^i=(\rho^i, M^i)$ by the formulas
\begin{align*}
\langle \rho_t^i,f\rangle                                                             &
=\int_{(h_i,+\infty)}f(x-h_i)\rho_{(\alpha_i+t)\wedge
\beta_i}(dx)\\                        
\langle M_t^i,\varphi\rangle                         &
=\int_{(h_i,+\infty)\times[0,+\infty)}\varphi(x-h_i, v)M_{(\alpha_i+t)\wedge
\beta_i}(dx,dv). 
\end{align*}

It  is easy to  adapt Lemma  4.2.4. of \cite{dlg:rtlpsbp} to  get the
following Lemma.
\begin{lem}
\label{lem:dlg-decomp}
   Let $(\mu, \Pi)\in \cm_f(\R_+)\times \cm(\R_+^2)$. The point
   measure $\displaystyle \sum_{i\in I} \delta_{(h_i,\cs^i)}$ is under
    $\P^*_{\mu, \Pi}$ a Poisson point measure with intensity $\mu(dr) \N[d\cs]$.
\end{lem}

\subsection{The process $m^{(\theta)}$}
For   $\theta\geq   0$,   we define the  $\cm(\R_+)$-valued   process
$m^{(\theta)}=(m_t^{(\theta)},t\geq 0)$ by
\begin{equation}
   \label{eq:def_m}
m_t^{(\theta)}(dr)=M_t(dr,(0,\theta]). 
\end{equation}
We make two remarks. We  have for $s>0$, 
\begin{equation}
   \label{eq:m=0/x}
\P_{0,0}(m^{(\theta)}_s=0|\cx)= \expp{- \theta \sum_{{0<r\le s}, \;
  {X_{s-}<I_t^s}} \Delta_s}= \expp{- \theta \langle\kappa_s,1\rangle }.
\end{equation}
Notice that for $s\in \cj$, i.e.  $\Delta_s>0$, we have  $M_s(\{H_s\}, dv)=\Delta_s \sum_{u\geq
0} \delta_{V_{s,u}}(dv)$, where conditionally on $\cx$, $\sum_{u\geq 0}
\delta_{V_{s,u}}(dv)$ is a Poisson point measure with intensity $\Delta_sdu$. In
particular, we have
\[
\P_{\mu,\Pi}( m_s^{(\theta)}(\{H_s\})>0|\cx)=\P(M_s(\{H_s\}\times 
(0,\theta])>0|\cx)=1-\expp{-\theta\Delta_s}. 
\]
{F}rom Poisson point measure properties, we get the following Lemma.
\begin{lem}
\label{lem:PPPjumps}
   The  pruned random  measure $\displaystyle  \cx^\theta=\sum_{s\geq 0}
   \ind_{\{m^{(\theta)}_s(\{H_s\})>  0\}}  \delta_{s,  \Delta_s}$  is  a
   Poisson         point          process         with         intensity
\begin{equation}
   \label{eq:def-nq}
n^\theta(d\ell)=(1-\expp{-\theta \ell}) \pi(d\ell).
\end{equation}
\end{lem}

We   shall use later the following property, which is a
consequence of Poisson point measure properties.
\begin{prop}
\label{prop:M_trans}
Let $M^\theta=(M^\theta_t, t\geq 0)$ be the measure-valued process defined by
\[
M^\theta _t(dr, [0,a])=M_t(dr, (\theta,\theta+a]), \quad \text{for all
  $a\geq 0$}. 
\]
Then, given $\rho$, $M^\theta$ is independent of $M\ind_{\R_+\times
  [0,\theta]}$  and is distributed as $M$.
\end{prop}

Eventually, the next Lemma on time reversibility can easily be deduced
from Corollary 3.1.6 of  \cite{dlg:rtlpsbp} and the construction of
$M$. 

\begin{lem}
   \label{lem:reversib}
   Under    $\N$,    the    processes   $((\rho_s,    \eta_s,    \ind_{\{
     m^{(\theta)}_s=0\}}),      s     \in     [0,      \sigma])$     and
   $((\eta_{(\sigma-s)-}, \rho_{(\sigma-s)-},  $ $\ind_{\{ m^{(\theta)}_
     {(\sigma-s)-}=0\}}),   s   \in   [0,   \sigma])$  have   the   same
   distribution.
\end{lem}


\section{The pruned exploration process}\label{sec:pruned_process}

Let us fix $\theta>0$. We shall write $m$ for the process $m^{(\theta)}$
defined  in the previous  Section.  We  define the  following continuous
additive functional  of the process $((\rho_t,m_t),t\ge  0)$: for $t\geq
0$
\[
A_t=\int_0^t \ind_{\{m_s=0\}}\; ds,
\]
and $C_t=\inf \{  r>0; A_r > t\}$ its  right continuous inverse, with
the  convention  that  $\inf\emptyset=\infty   $.   Notice  $C_t$  is  a
$\cf$-stopping  time for  any $t\geq  0$ and is finite a.s. from Corollary \ref{cor:C_0=0} below.

  We   define  the   pruned  exploration   process   $\tilde  \rho=(\tilde
\rho_t=\rho_{C_t}, t\geq 0)$ and the pruned Lévy Poisson snake $\tilde
\cs= ( \tilde \rho, \tilde M)$, where $\tilde M=(M_{C_t}, t\geq 0)$.  In
particular the law of $\tilde M$ knowing $\tilde \rho$ is the law of $M$
knowing $\rho=\tilde \rho$.  Notice  the process $\tilde \rho$, and thus
the  process $\tilde  \cs$,  is  càd-làg.  We  also  set $\tilde  H_t=
H_{C_t}$.  Let $\tilde \cf=(\tilde \cf_t, t\geq 0)$ be the filtration generated by the
pruned  exploration process $\tilde  \cs$ completed  the usual  way.  In
particular   $\tilde \cf_t\subset   \cf_{C_t}$,   where   if   $\tau$   is   an
$\cf$-stopping time, then $\cf_\tau$ is the $\sigma$-field associated to
$\tau$.

We introduce the following Laplace exponent $\psi^{(\theta)}$  defined for
  $\lambda\geq 0$ by 
\begin{equation}
   \label{eq:def_psi-q}
{\psi^{(\theta)}
(\lambda)=\psi(\lambda+\theta)-\psi(\theta).}
\end{equation}

\begin{lem}
\label{lem:A_s=0}
   We have the following properties.
\begin{itemize}
   \item[(i)] For $\lambda>0$, $\N[1  -\expp{-\lambda
     A_\sigma}]={\psi^{(\theta)}}^{-1} (\lambda)$.
   \item[(ii)] $\N$-a.e. 0 and $\sigma$ are points of increase for
   $A$. More precisely, $\N$-a.e. for all $\varepsilon>0$, we have 
$A_\varepsilon>0$ and $A_\sigma-A_{(\sigma-\varepsilon) \vee 0}>0$.
   \item[(iii)] $\N$-a.e. the set $\{s; m_s\neq 0\}$ is dense in
   $[0,\sigma]$. 
\end{itemize}
\end{lem}

Before  going into  the proof  of this  Lemma, let  us state  two direct
consequences.     {F}rom    excursion    decomposition,     see    Lemma
\ref{lem:dlg-decomp}, the  second part of  Lemma \ref{lem:A_s=0} implies
the following corollary.
\begin{cor}
\label{cor:C_0=0}
For any initial measures $\mu,\Pi$,   $\P_{\mu,\Pi}$-a.s. the process $(C_t,t\geq 0)$ is finite and starts at 0 if $m_0=0$.
\end{cor}

We  define  $\tilde \sigma=\inf\{t>0;  \tilde  \rho_t=0\}$.  {F}rom  the
second  part of  Lemma \ref{lem:A_s=0},  we get  that $\sigma=\inf\{t>0;
\rho_t=0\}$  is   a  left  increasing   point  of  $A$   ($\N$-a.e.   or
$\P_{(\mu,\Pi)}$-a.s., $\mu\neq 0$). Therefore, we have $\lim_{r\uparrow
  A_\sigma} C_r=\sigma$.  As $\rho$ is  left continuous at  $\sigma$, we
get that $\lim_{r\uparrow A_\sigma}  \tilde \rho_r=0$ which implies that
$\tilde\sigma\le A_\sigma$.  Since $\tilde \sigma \geq A_\sigma$, we get
that $\N$-a.e.
\begin{equation}
   \label{eq:s-tilde}
\tilde \sigma= A_\sigma.
\end{equation}
This equality holds also $\P_{(\mu,\Pi)}$-a.s., for $\mu\neq 0$. 

\begin{proof}[Proof of Lemma \ref{lem:A_s=0}.]
We first prove (i). 
Let  $\lambda>0$. Before  computing $v=\N[1  -\exp{-\lambda A_\sigma}]$,
notice  that  $A_\sigma\leq  \sigma$  implies, thanks  to
\reff{eq:N_s}, that 
$v\leq \N[1 -\exp{-\lambda \sigma}] =\psi^{-1}(\lambda)<+\infty$.  We have
\[
   v= \lambda \N\left[\int_0^\sigma dA_t\, \expp{-\lambda \int_t^\sigma
       dA_u}\right] 
=\lambda \N\left[\int_0^\sigma dA_t\, \E^*_{\rho_t, 0}[ \expp{-\lambda
  A_\sigma}] \right],
\]
where  we  replaced $\expp{-\lambda  \int_t^\sigma  dA_u}$  in the  last
equality by $\E^*_{\rho_t,  0}[ \expp{-\lambda A_\sigma}]$, its optional
projection.   In order  to  compute  this last  expression,  we use  the
decomposition of  $\cs$ under $\P^*_{\mu, \Pi}$  according to excursions
of   the  total   mass  of   $\rho$   above  its   minimum,  see   Lemma
\ref{lem:dlg-decomp}.  Using the same notations as in this Lemma, notice
that under  $\P^*_{\mu, 0}$,  we have $A_\sigma=A_\infty  =\sum_{i\in I}
A^i_\infty $, with
\begin{equation}
   \label{eq:def-Ai}
A^i_T =\int_0^T \ind_{\{M^i_t(\R_+\times [0,\theta])=0\}} dt.
\end{equation}
By Lemma \ref{lem:dlg-decomp}, we get 
\[
\E^*_{\mu,  0}[ \expp{-\lambda A_\sigma}] = \expp{ - \langle \mu,1\rangle  \N[ 1-
  \exp{-\lambda A_\sigma}]}=\expp{-v\langle \mu,1\rangle }.
\]
Now, for fixed $t$, recall \reff{eq:m=0/x}. By conditioning with respect
to $\cx$ or to $\rho$ thanks to Remark \ref{rm:recons}, we have
\[
 v=\lambda  \N\Big[\int_0^\sigma dA_t \, \expp{ -v \langle \rho_t,1\rangle }\Big]
= \lambda \N\Big[\int_0^\sigma dt \,\ind_{\{m_t=0\}} \expp{ -v
    \langle \rho_t,1\rangle }\Big] 
= \lambda \N\Big[\int_0^\sigma dt  \expp{ -(v+\theta) \langle \rho_t,1\rangle  - \theta
  \langle \eta_t,1\rangle  }\Big].  
\]
Now we use Proposition \ref{prop:poisson_representation1} to get 
\begin{align}
\nonumber
   v
&=\lambda\int_0^{+\infty}da\, \expp{-\alpha_0   a}\mathbb{M}_a [\expp{
  -(v+\theta) \langle \mu,1\rangle  - \theta   \langle \nu,1\rangle  }] \\
\nonumber
&=\lambda\int_0^{+\infty}da\, \expp{-\alpha_0   a} \exp\Big\{- \int_0^a dx \int_0^1
  du \int _{(0,\infty )} \ell\pi(d\ell) \left[ 1- \expp{
  -(v+\theta) u\ell - \theta (1-u) \ell }\right]\Big\} \\
\label{eq:calc-v0}
&=\lambda\int_0^{+\infty}da\, \exp\Big\{- a \int_0^1  du
  \,\psi'(\theta+vu)\Big\}\\  
&=\lambda\frac{v}{\psi(\theta+v) -\psi(\theta)}.
\label{eq:calc-v}
\end{align}
where, for the third equality,  we used 
\begin{equation}
   \label{eq:def-psi1}
\psi'(\lambda)= \alpha_0+ \int_{(0,\infty )} \pi(d\ell)\; 
\ell(1-\expp{-\lambda \ell }).
\end{equation}
Notice   that    if   $v=0$,   then    \reff{eq:calc-v0}   implies   $v=
\lambda/\psi'(\theta)  $,  which  is  absurd. Therefore  we  have  $v\in
(0,\infty)$, and we can divide \reff{eq:calc-v} by $v$ to get
$\psi^{(\theta)}(v)=\lambda$. This proves (i). 

Now, we prove  (ii). If we let $\lambda \rightarrow \infty  $ in (i) and
use  that $\lim_{r\rightarrow \infty  } \psi^{(\theta)}(r)=+\infty  $, then  we get
that $\N[A_\sigma>0]=+\infty  $.  Notice that for  $(\mu,\Pi)\in \S$, we
have under  $\P^*_{\mu, \Pi}$, $A_\infty  = \sum_{i\in I}  A^i_\infty $,
with $A_i$ defined  by \reff{eq:def-Ai}. Thus Lemma \ref{lem:dlg-decomp}
imply that if $\mu\neq 0$,  then $\P^*_{\mu, \Pi}$-a.s.  $I$ is infinite
and $A_\infty >0$.   Using the Markov property at time  $t$ of the snake
under $\N$, we get that  for any $t>0$, $\N$-a.e.  on $\{\sigma>t\}$, we
have  $A_\sigma-A_t>0$.   This  implies  that  $\sigma$ is  a  point  of
increase   of   $A$  $\N$-a.e.   By   time   reversibility,  see   Lemma
\ref{lem:reversib}, we also  get that $0$ is a point  of increase of $A$
$\N$-a.e.

To prove  (iii), recall that $\int_{(0,1 )} \ell \pi(d\ell)=+\infty
  $ implies that $\cj=\{s\geq 0; \Delta_s>0\}$ is dense in $\R_+$
  a.s. Moreover, for every $t>r\geq 0$,
$$\sum_{r\le s\le t}\Delta_s=+\infty\ a.s.$$
Now, by the properties of Poisson point measures, we have
$$\P(\forall s\in[r,t],\ m_s=0)=\E\left[\expp{-\sum_{r\le s\le t}\Delta _s}\right]=0$$
which proves (iii).

\end{proof}


\section{A special Markov property}\label{sec:Markov_special}

Let us fix $\theta>0$ and  use the notations of the previous Section.

In order  to define  the excursion  of the Lévy  Poisson snake  out of
$\{s\ge  0;\ m_s=0\}$,  we define  $O$ as  the interior  of  $\{s\ge 0,\
m_s\ne 0\}$. 

\begin{lem}
\label{lem:O}
    $\N$-a.e. the open set $O$ is non empty. 
\end{lem}

\begin{proof}
Thanks to Lemma  \ref{lem:A_s=0}, (iii), $\{s\ge 0,\ m_s\ne  0\}$ is non
empty. For any element, $s$, of  this set, there exists $u\leq H_s$ such
that $m_s([0,u])\ne  0$ and $\rho_s(\{u\})>0$.  Then we consider
$\tau_s=\inf \{t>s,\rho_t(\{u\})=0\}$. 
By the right continuity of $\rho$, $\tau_s>s$ and clearly
$(s,\tau_s)\subset O$ $\N$-a.e. Therefore $O$ in non empty.
\end{proof}

We  write  $O=\bigcup_{i\in  I}(\alpha_i,\beta_i)$ and say that
$(\alpha_i,\beta_i)_{i\in I}$ are  the excursions
intervals  of the Lévy  Poisson snake  $\cs=(\rho, M)$  out of
$\{s\ge 0,\ m_s=0\}$.

Next  we prove a  special Markov  property out  of $\{s\ge  0,\ m_s=0\}$
under the excursion measure $\N$. Using the right continuity of $\rho$
and the definition of $M$, we get that for $i\in I$, $\alpha_i>0$,
$\alpha_i\in \cj$, that is $\rho_{\alpha_i}(\{ H_{\alpha_i}\})=
\Delta_{\alpha_i}$, and $M_{\alpha_i}([0, H_{\alpha_i}),
[0,\theta])=0$. For  every $i\in I$, let us define the
measure-valued  process  $\cs^i=(\rho^i,M^i)$  by:  for  every  $f\in
\cb_+(\R_+)$, $\varphi\in \cb_+(\R_+^2)$, $t\ge 0$,
\begin{align*}
\langle \rho_t^i,f\rangle                                                             &
=\int_{[H_{\alpha_i},+\infty)}f(x-H_{\alpha_i})\rho_{(\alpha_i+t)\wedge
\beta_i}(dx)\\                        
\langle M_t^i,\varphi\rangle                         &
=\int_{(H_{\alpha_i},+\infty)\times [0,+\infty)}\varphi(x-H_{\alpha_i}, v)M_{(\alpha_i+t)\wedge
\beta_i}(dx,dv). 
\end{align*}
Since   $\rho^i_0= \delta_ {\Delta_{\alpha_i}}$, with
$\Delta_{\alpha_i}>0$,  then   for every $t<
\beta_i-\alpha_i$, the  measure $\rho^i_t$ charges  0. As $M^i_0=0$
we have for every $t< \beta_i-\alpha_i$, $M^i_t(\{0\}\times \R_+)=0$.
We call $\Delta_{\alpha_i}$ the starting mass of $\cs^i$.

Recall  $\tilde \cf_\infty  $  is   the  $\sigma$-field  generated  by  $\tilde
\cs=((\rho_{C_t},    M_{C_t}),   t\geq    0)$    and   $\P_{\mu,    \Pi
}^{*}(d\cs)$  denotes  the  law   of  the  snake  $\cs$  started  at
$(\mu,\Pi)$    and    stopped    when    $\rho$    reaches    0.     For
$\ell\in[0,+\infty)$,     we      will     write     $\P_\ell^*$     for
$\P_{\delta_\ell,0}^*$. Recall \reff{eq:def-nq} and define the measure $\rN$ by
\begin{equation}
   \label{eq:def_rN}
\rN(d\cs)=\int_{(0,+\infty)}\pi(d\ell)
\left(1-\expp{-\theta\ell}\right)\P_\ell^*(d\cs)=\int_{(0,\infty )}
n^{(\theta)} (d\ell) \P_\ell^*(d\cs). 
\end{equation}

If  $Q$ is  a measure  on $\S$  and $\phi$  is a  non-negative measurable
function defined on a  measurable space $\R_+\times \Omega\times \S$, we
denote by
\[
Q[\phi(u,\omega,\cdot)]=\int _{\S}\phi(u,\omega,\cs)Q(d\cs).
\]
In other words, the integration concerns only the third component of the
function $\phi$.

Recall the definition of $\tilde \sigma$ given after Corollary
\ref{cor:C_0=0}.

\begin{theo}(Special Markov property)
\label{th:SMP}
   Let  $\phi$   be  a   non-negative  measurable  function   defined  on
   $\R_+\times \Omega\times \S$ such that $t\mapsto  \phi(t,\omega,\cs)$ is progressively $
   \tilde \cf$-measurable for any $\cs\in\S$. Then, we have $\N$-a.e. 
\begin{equation}
   \label{eq:MS}
\N\left[\exp\left(-\sum_{i\in
    I}\phi(A_{\alpha_i},\omega,\cs^i)\right)\biggm|\tilde \cf_\infty \right]
=\exp\left(-\int_0^{\tilde   \sigma}
    du\,\rN\left[1-\expp{-\phi(u,\omega,\cdot)}\right]\right).
\end{equation}
In  particular, the  law of  the excursion  process  $\displaystyle \sum
  _{i\in  I}\delta_{(A_{\alpha_i},\cs^i)}$,  given  $\tilde \cf_\infty$
  under  $\N$,  is the  law  of a  Poisson  point  measure of  intensity
  $\ind_{[0, \tilde \sigma]}(u) \;  du\; \rN(d\cs)$. 
\end{theo}

Before going into the proof of this Theorem, we give a corollary we
shall use later. 

\begin{cor} 
  \label{cor:super}  
The law of the  excursion process $\displaystyle \sum _{i\in I}\delta_{(
A_{\alpha_i},\rho_{\alpha_i-},  \cs^i)}$,  given $\tilde \cf_\infty
$,  is  the  law of  a  Poisson  point  measure  of  intensity $
\ind_{ [0, \tilde \sigma]}(u) du  \;
\delta_{ \tilde \rho_u} (d\mu)\; \rN(d\cs)$.
\end{cor}

\begin{proof}
   This is a direct consequence of Theorem \ref{th:SMP} and Lemma
   \ref{lem:cor=s_ex}. 
\end{proof}

The rest of this Section is devoted to the proof of the special Markov
property. 

\subsection{A remark and notations}
To begin with, let us remark  that to prove Theorem \ref{th:SMP}, we may
only  consider   function $\phi$ satisfying the hypothesis of Theorem
\ref{th:SMP}  and  those two conditions:
\begin{itemize}
   \item[($h_1$)] $\phi(s,\omega,\cs)=0$  if the starting mass of $\cs$ is less than
     $\eta$, that is $\langle \rho_0,1\rangle \le \eta$,   for a fixed
positive real number $\eta$.
   \item[($h_2$)] $t\mapsto \phi(t,\omega,\cs)$ is continuous for all
$\cs\in \S$ a.s. 
\end{itemize} 
Indeed if \reff{eq:MS}  holds for such functions then  by monotone class
Theorem and monotonicity it holds also for every function satisfying the
hypothesis  of  Theorem  \ref{th:SMP}.  {F}rom  now on,  but  for  Lemma
\ref{lem:cor=s_ex}, we  fix $\eta>0$, and we assume  the function $\phi$
satisfies the  hypothesis of Theorem  \ref{th:SMP} and ($h_1$).  We will
assume    ($h_2$)    only    for   Sections    \ref{sec:comp_lim}    and
\ref{sec:proof_th_MS}.

Let $\varepsilon<\eta$ and let us define by induction (under the
measure $\N$) the following stopping times:
$T_0^\varepsilon=0$
and, for every integer $k\ge 0$,
\begin{align*}
S_{k+1}^\varepsilon & =\inf\left\{s>T_k^\varepsilon,\ m_s(\{H_s\})>0,\
\rho_s(\{H_s\})>\varepsilon\right\}\\
T_{k+1}^\varepsilon & =\inf\left\{s>S_{k+1}^\varepsilon,\ m_s=0\right\}
\end{align*}
with the convention $\inf\emptyset=\sigma$.
Let us then denote
\begin{equation}
   \label{eq:N-e}
N_\varepsilon=\sup\{k\in\N,\ S_k^\varepsilon\ne \sigma\}.
\end{equation}
Notice $N_\varepsilon$ is finite $\N$-a.e. as there is a finite number
of jumps $\Delta_t>\varepsilon$. 

For  every $k\le  N_\varepsilon$, we  define the  measure-valued process
$\cs^{k,\varepsilon}=(\rho^{k,\varepsilon},M^{k,\varepsilon})$   in  the
same way  as the  processes $\rho^i$ and  $M^i$: for  every non-negative
continuous functions $f$ and $\varphi$, and $t\ge 0$,
\begin{align*}
\langle \rho_t^{k,\varepsilon},f\rangle  &
=\int_{[H_{S_k^\varepsilon},+\infty)}f(x-H_{S_k^\varepsilon})\rho_{(S_k^\varepsilon+t)\wedge
T_k^\varepsilon}(dx)\\
\langle M_t^{k,\varepsilon},\varphi\rangle  &
=\int_{(H_{S_k^\varepsilon},+\infty)
  \times[0,+\infty)}\varphi(x-H_{S_k^\varepsilon},v)
M_{(S_k^\varepsilon+t)\wedge   T_k^\varepsilon}(dx,dv).
\end{align*}
We call $\Delta_{S^{\varepsilon}_k}$ the starting mass of $\cs^{k,\varepsilon}$. Notice that
$\rho^{k,\varepsilon}_0=\delta_{\Delta_{S^{\varepsilon}_k}} $  and
$\Delta_{S^{\varepsilon}_k} \ge \varepsilon$ for $k\in \{1, \ldots,
N_\varepsilon\}$. 
Notice also that $\N$-a.e, 
\begin{equation}
   \label{eq:cv_st}
\lim_{\varepsilon\to
  0}\bigcup_{k\in\N}(S_k^\varepsilon,T_k^\varepsilon)=\bigcup_{i\in
  I}(\alpha_i,\beta_i).
\end{equation}

\subsection{Approximation of the excursion process}

\begin{lem}\label{lem:approx}
$\N$-a.e., we have for $\varepsilon>0$ small enough
\begin{equation}
   \label{eq:=cs}
\sum_{i\in I}\phi(A_{\alpha_i},\omega,\cs^i)= \sum_{k=2}^{N_\varepsilon}
  \phi(A_{S_k^\varepsilon},\omega,\cs^{k,\varepsilon}). 
\end{equation}
\end{lem}
\begin{proof}
Let $I_\eta$ be the set of indexes $i\in I$, such that the starting mass
of $\cs^i$ is larger than $\eta$. Because of ($h_1$), we have
$$\sum_{i\in I}\phi(A_{\alpha_i},\omega,\cs^i)=\sum_{i\in
  I_\eta}\phi(A_{\alpha_i},\omega,\cs^i).$$

Let $\varepsilon<\eta$. Then, for any $i\in I_\eta$,
there         exists        $k\in        \N^*$,         such        that
$\cs^{k,\varepsilon}=\cs^i$. 

Furthermore,  all  the  others  excursions
$\cs^{k,\varepsilon}$ which don't belong to $\{\cs^i,  i\in I_\eta\}$
either have a starting mass less that $\eta$ (and thus
$\phi(A_{S_k^\varepsilon},\omega,\cs^{k,\varepsilon})=0$), or have a
starting mass greater that $\eta$ but
$m_{S_k^\varepsilon}([0,H_{S_k^\varepsilon}))>0$. But, as the set
    $\{s\ge 0, \Delta_s>\eta\}$ is finite, there exists only a finite
    number of excursions $\cs^i$ which straddle a time $s$ such that
    $\Delta_s>\eta$. Therefore, the minimum over those excursions of
    their starting mass, say $\eta'$, is positive a.s. and, if we
    choose $\varepsilon<\eta'$, there are no excursions
    $\cs^{k,\varepsilon}$ with initial mass greater than $\eta$ which
    do not correspond to a $\cs ^i$.

Consequently, if we choose $\varepsilon<\eta\wedge\eta'$,
 we
  have $$\sum_{i\in I}\phi(A_{\alpha_i},\omega,\cs^i)=
  \sum_{k=1}^{N_\varepsilon} 
  \phi(A_{S_k^\varepsilon},\omega,\cs^{k,\varepsilon}).$$

Notice  also,   that  because   of  Lemma  \ref{lem:A_s=0}   (iii),  for
$\varepsilon>0$ small enough, the starting mass of $\cs^{1,\varepsilon}$
is  less than  $\eta$.  Therefore, we  deduce  that \reff{eq:=cs}  holds
$\N$-a.e. for $\varepsilon>0$ small enough.

\end{proof}

We  can now prove the next Lemma which gives   Corollary
\ref{cor:super}. 
\begin{lem}
   \label{lem:cor=s_ex}
  Let  $\psi$  be a  bounded  non-negative measurable  function  defined  on
$\R_+\times \cm_f(\R_+)\times \S$.  $\N$-a.e., we have
\begin{equation*}
\sum_{i\in I}\psi(A_{\alpha_i},\rho_{\alpha_i-},\cs^i)=
\sum_{i\in I}\phi(A_{\alpha_i},\omega,\cs^i), 
\end{equation*}
where $\phi(t,\omega,\cs)=\psi(t,\tilde \rho_t(\omega), \cs)$. 
\end{lem}

\begin{proof}
   First we assume that $\psi(t,\mu,\cs)=0$ if the starting mass of
   $\cs$ is less than
   $\eta$. The same arguments as those used to prove Lemma \ref{lem:approx}
   yields that $\N$-a.e. for $\varepsilon>0$ small enough, we have 
\[
\sum_{i\in I}\psi(A_{\alpha_i},\rho_{\alpha_i-}, \cs^i)=
  \sum_{k=2}^{N_\varepsilon} 
  \psi(A_{S_k^\varepsilon},\rho_{S_k^\varepsilon-},\cs^{k,\varepsilon}). 
\]
Notice     that    by    construction,     $\rho_{S_k^\varepsilon-}    =
\rho_{T_k^\varepsilon}$  and  that  $m_{T^\varepsilon_k}=0$.  Using  the
strong Markov property at time $T_k^\varepsilon$ and the second part of
Corollary \ref{cor:C_0=0}, we deduce that $\N$-a.e. for all $k\in
\N^*$, 
\begin{equation}
    \label{eq:cat=t}
C_{A_{T_k^\varepsilon}}=T_k^\varepsilon. 
\end{equation}
Therefore, as
$A_{S_k^\varepsilon}= A_{T_k^\varepsilon}$, 
we have
$\N$-a.e. 
\[
\tilde \rho_{A_{S_k^\varepsilon}}=
\tilde \rho_{A_{T_k^\varepsilon}}=
\rho_{T_k^\varepsilon}= \rho_{S_k^\varepsilon-}.
\]
Hence, we have that $\N$-a.e. for $\varepsilon>0$ small enough,
\[
\sum_{i\in I}\psi(A_{\alpha_i},\rho_{\alpha_i-}, \cs^i)=
  \sum_{k=2}^{N_\varepsilon} 
  \phi(A_{S_k^\varepsilon},\omega,\cs^{k,\varepsilon}), 
\]
with $\phi(t,\omega,\cs)=\psi(t,\tilde \rho_t(\omega), \cs)$. Now, we complete the proof
using Lemma \ref{lem:approx} and letting $\eta\downarrow 0$.

\end{proof}

\subsection{A measurability result}

We shall use later the next additive functional defined  for $s\geq 0$
by 
\begin{equation}
   \label{eq:def-Ae}
 A_s^\varepsilon=\sum_{k=0}^{N_\varepsilon}\int_0^s\ind_{
  [T_k^\varepsilon,S_{k+1}^\varepsilon]}(u)\;du.
\end{equation}

For $k\geq 1$, we  consider the $\sigma$-field
$\cf^{(\varepsilon), k}$ generated by the family of processes
\[
\left(\cs_{(T_l^\varepsilon+s)\wedge
S_{l+1}^\varepsilon   -},\   s>0\right)_{ l\in  \{0, \ldots,   k-1\}}.
\]
Notice that 
for $k\in \N^*$
\begin{equation}
   \label{eq:inclusion}
 \cf^{(\varepsilon),k} \subset \cf_{S^\varepsilon_k}.
\end{equation}

\begin{lem}
\label{lem:mesuraiblite_f}
For    any    $\varepsilon>0$,    $k\in    \N^*$,   the    function    $
\phi(A_{S^\varepsilon_k},\omega  ,  \cdot)$ is  
$\cf^{(\varepsilon), k}$-measurable.
\end{lem}

\begin{proof}
We set $C_s^\varepsilon$ the right continuous inverse of
$A_s^\varepsilon$ and we
define the filtration
$\tilde \cf^{(\varepsilon)}=(\tilde \cf^{(\varepsilon)}_t, t\geq 0)$ generated by
the process $(\cs_{C_s^\varepsilon},s\ge 0)$.

We consider the
counting process $(R_t, t\geq 0)$ defined by 
$\displaystyle R_t=\inf\{ k\geq 0;
S^\varepsilon_{k+1}>A^\varepsilon_t\}$.
Consider the filtration
$\cf^{(\varepsilon)}=(\cf^{(\varepsilon)}_t, t\geq 0)$, where $
\cf^{(\varepsilon)}_t =\tilde \cf^{(\varepsilon)}_t \vee \sigma(R_s, s\leq
t)$. In particular for $k\geq 1$,
$A^\varepsilon_{S^\varepsilon_k}=\inf\{t\geq 0; R_t=k\}$  is a
$\cf^{(\varepsilon)}$-stopping time. Notice then that 
$\cf^{(\varepsilon),
  k}=\cf^{(\varepsilon)}_{A^\varepsilon_{S^\varepsilon_k}}$. 

By  the monotone  class Theorem,  to prove  the Lemma,  it is  enough to
consider   simple   processes,   $\phi$,  defined   by   $\phi(t,\omega,
\cs)=g(\cs)Z\ind_{\{r\leq t\}}$, where $r\geq  0$, $Z\in \tilde \cf_r$, and $g$
is a real measurable function defined on $\S$.  For $k\in \N^*$, we have
$\phi(A_    {S^\varepsilon_k},\omega,     \cdot)=    gZ    \ind_{\{r\leq
  A_{S^\varepsilon_k}\}} $.  
Notice that 
\begin{align*}
   A^\varepsilon_{C_r}
&=\inf\{u>0; C_u^\varepsilon> C_r\}\\
&=\inf\{u>0; \int_0^{C_u^\varepsilon} \ind_{\{m_s=0\}} ds > r\}\\
&=\inf\{u>0; \int_0^{C_u^\varepsilon} \ind_{\{m_s=0\}} dA^\varepsilon_s > r\}\\
&=\inf\{u>0; \int_0^u \ind_{\{m_{C^\varepsilon_t}=0\}} dt > r\},
\end{align*}
where we  used that $A^\varepsilon$  is the right continuous  inverse of
$C^\varepsilon$  for the  first equality,  $C$ is  the  right continuous
inverse of  $A$ for the  second, $\{s;m_s=0\} \subset  \bigcup_{k\geq 0}
[T^\varepsilon _k,  S^\varepsilon_{k+1}]$ for the third,  and the change
of variable $t=A^\varepsilon_s$ for the last.
This   gives  that   $A^\varepsilon_{C_r}$  is   a  $\cf^{(\varepsilon)}
$-stopping time. By composition of random change time, we also have
$\tilde \cf_r\subset \cf^{(\varepsilon)}_{ A^\varepsilon_{C_r}}$. 
Eventually, we have 
\[
\{r\leq  A_{S^\varepsilon_k}\}
=\{r\leq  A_{T^\varepsilon_k}\}
=\{C_r\leq T^\varepsilon_k\}=\{A^\varepsilon_{C_r}\leq
A^\varepsilon_{T^\varepsilon_k} \}
=\{A^\varepsilon_{C_r}\leq
A^\varepsilon_{S^\varepsilon_k} \},
\]
where we  used $A_{S^\varepsilon_k} =A_{T^\varepsilon_k}$  for the first
equality, \reff{eq:cat=t} and the definition  of $C$ for the second, and
similar properties for $A^\varepsilon$ for the two last ones.  We deduce
then      that       $Z      \ind_{\{r\leq      A_{S^\varepsilon_k}\}}=Z
\ind_{\{A^\varepsilon_{C_r} \leq A^\varepsilon_{S^\varepsilon_k}\}} $ is
measurable            with            respect            to            $
\cf^{(\varepsilon)}_{A^\varepsilon_{S^\varepsilon_k}}=\cf^{(\varepsilon),
  k}$.  This ends the proof of the Lemma.
   
\end{proof}

\subsection{Computation of the conditional expectation of the approximation}
\begin{lem}
\label{lem:ccea}
For every $\tilde \cf_\infty$-measurable non-negative random variable
$Z$, we have 
$$
\N\left[Z\exp\left(-\sum_{k=2}^{N_\varepsilon}
    \phi(A_{S_k^\varepsilon},\omega,\cs^{k,\varepsilon})\right)\right]=  
\N\left[Z\prod_{k=2}^{N_\varepsilon}
\rN\left[\expp{-\phi(  A_{S_k^\varepsilon},\omega, \cdot)}\Bigm|\rho_0>\varepsilon\right] \right]. 
$$
\end{lem}

\begin{proof}

For every  integer $p\ge 2$,  we consider a non-negative  random variable
$Z$    of    the    form    $Z=Z_0   Z_1$,    where    $Z_0\in
\cf^{(\varepsilon), p} $   and  $Z_1\in  \sigma  (\cs_{(T_k^\varepsilon+s)\wedge
S_{k+1}^\varepsilon  -},\ s\ge0,\  k\geq p)$  are bounded  non-negative and
such that $\N[Z_0]<\infty $.

To compute $\displaystyle 
\N\left[Z\exp\left(-\sum_{k=2}^p\phi(A_{S_k^\varepsilon},\omega,
    \cs^{k,\varepsilon})\right)\right]$,   
we first apply the strong Markov property at time $T_p^\varepsilon$. We obtain
\[
\N\left[Z\exp\left(-\sum_{k=2}^p
  \phi(A_{S_k^\varepsilon},\omega,\cs^{k,\varepsilon}) 
  \right)\right]
=\N\left[Z_0\exp\left(-\sum_{k=2}^p
  \phi(A_{S_k^\varepsilon}, \omega,\cs^{k,\varepsilon})\right)
  \E_{\rho_{T_p^\varepsilon},0}^* 
\bigl[Z_1]\right].
\] 
Notice    that   $\rho_{T_p^\varepsilon}=\rho_{S_p^\varepsilon-}$,   and
consequently  $\rho_{T_p^\varepsilon}$  is  measurable with  respect  to
$\mathcal{F}_{S_p^\varepsilon}$.   So,  when we  use  the strong  Markov
property at time $S_p^\varepsilon$, we get thanks to Lemma
\ref{lem:mesuraiblite_f} and \reff{eq:inclusion}, 
\begin{multline*}
\N\left[Z\exp\left(-\sum_{k=2}^p\phi(A_{S_k^\varepsilon},\omega,
    \cs^{k,\varepsilon})\right)\right]\\  
=\N\left[Z_0\exp\left(-\sum_{k=2}^{p-1}\phi(A_{S_k^\varepsilon},\omega,
    \cs^{k,\varepsilon})\right)
 \E_{\rho_0^{p,\varepsilon},0}^*\left[\expp{-\phi(  A_{S_p^\varepsilon} ,\omega,\cdot)}\right]
   \E_{\rho_{T_p^\varepsilon},0}^*[Z_1]\right].
\end{multline*}

Recall  $p\geq  2$.   Conditionally on  $\cf_{T^\varepsilon_{p-1}}$,  on
$N_\varepsilon\geq p$,  the measure $\rho_0^{p,\varepsilon}$  is a Dirac
mass  and, by  the Poisson  representation of  Lemma \ref{lem:PPPjumps},
this mass  is the first atom  of the Poisson  point measure $\cx^\theta$
that  lies  in  $(\varepsilon,+\infty)$.   
Consequently,  the  mass  of
$\rho_0^{p,\varepsilon}$   is   distributed   according   to   the   law
$n^\theta(d\ell\,|\,  \ell>\varepsilon)$.
{F}rom  Poisson  point  measure
properties, notice  that $\rho_0^{p,\varepsilon}$ is  also independent of
$  \sigma(\cs_t,  t<  S^\varepsilon_{p})$
and thus of $\cf^{(\varepsilon),p}$.

Therefore,       conditionally      on       $N_\varepsilon\geq      p$,
$\rho_0^{p,\varepsilon}$       is       independent      of       $Z_0$,
$\rho_{T_p^\varepsilon}=  \rho_{S_p^\varepsilon-}$ and, thanks  to Lemma
\ref{lem:mesuraiblite_f} of  $\phi( A_{S_p^\varepsilon},\omega, \cdot)$.
So, by conditioning with respect to $\cf^{(\varepsilon),p}$, we get
\begin{multline}
\label{eq:eq-prec}
\N\left[Z\exp\left(-\sum_{k=2}^p
    \phi(A_{S_k^\varepsilon},\omega,\cs^{k,\varepsilon})\right)\right]\\  
=\N\left[Z_0\exp\left(-\sum_{k=2}^{p-1}\phi(A_{S_k^\varepsilon},\omega,
    \cs^{k,\varepsilon})\right)
\rN\left[\expp{-\phi(  A_{S_p^\varepsilon}, \omega,\cdot)}\Bigm|\rho_0>\varepsilon\right]
    \E_{\rho_{T_p^\varepsilon},0}^*[Z_1]\right]. 
\end{multline}

\begin{rem}
  \label{rem:e_p}  {F}rom point Poisson  measure property,  notice that,
  conditionally  on  $\cf_{T^\varepsilon_{p-1}}$ and  $N_\varepsilon\geq
  p\geq 2$, $e_p^\varepsilon= S_p^\varepsilon- T_{p-1}^ \varepsilon $ is
  an  exponential random  variable with  parameter 
\begin{equation}
   \label{eq:n_e}
n_\varepsilon=n^\theta(  \ell>\varepsilon)=
\int_{(\varepsilon,+\infty)}\pi(d\ell)\left(1-\expp{-\theta\ell}\right).
\end{equation}  
 And,  conditionally  on $N_\varepsilon$  and
  $N_\varepsilon\geq  2$, the  random variables  $(e_k^\varepsilon, k\in
  \{2,  \ldots,  N_\varepsilon\})$  are independent  exponential  random
  variables with parameter $n_\varepsilon$.
\end{rem}

Now,  using   one  more  time   the  strong  Markov  property   at  time
$T_p^\varepsilon$, we get from \reff{eq:eq-prec} 
\begin{multline*}
\N\left[Z\exp\left(-\sum_{k=2}^p\phi(A_{S_k^\varepsilon},\omega,
    \cs^{k,\varepsilon})\right)\right]\\ 
=\N\Bigg[Z
\rN\left[\expp{-\phi(  A_{S_p^\varepsilon},\omega, \cdot)}\Bigm|\rho_0>\varepsilon\right]\exp 
\left(-\sum_{k=2}^{p-1}\phi(A_{S_k^\varepsilon},\omega,
    \cs^{k,\varepsilon})\right)  \Bigg].
\end{multline*}

{{F}rom} monotone class  Theorem,  this equality holds  also for  any $Z\in
\cf^{(\varepsilon),\infty}    $    non-negative.     Thanks   to    Lemma
\ref{lem:mesuraiblite_f},     the     non-negative    random     variable
$Z'=Z\rN[\expp{-\phi(
A_{S_p^\varepsilon},\omega,\cdot)}|\rho_0>\varepsilon]$        is        measurable
w.r.t. $\cf^{(\varepsilon),\infty}  $.  So, we may  iterate the previous
argument and eventually get that  for any non-negative random variable $Z
\in \cf^{(\varepsilon),\infty} $, we have
\[
\N\left[Z\exp\left(-\sum_{k=2}^p\phi(A_{S_k^\varepsilon},\omega,
    \cs^{k,\varepsilon})\right)\right]=  
\N\left[Z\prod_{k=2}^p
\rN\left[\expp{-\phi(  A_{S_k^\varepsilon},\omega, \cdot)}\Bigm|\rho_0>\varepsilon\right] \right].
\]
Let $p\to +\infty$ and notice that $\tilde \cf_\infty  \subset
\cf^{(\varepsilon),\infty} $ to end the proof. 
\end{proof}

\subsection{An ancillary result}

Recall \reff{eq:N-e} and \reff{eq:n_e}
We prove the next result.

\begin{lem}
\label{lem:cv-Ne}
There exists a positive sequence $(\varepsilon_j, j\in \N^*)$ decreasing
to 0, such that  $\N$-a.e.:
\begin{itemize}
   \item[(i)]  $\displaystyle \lim_{j\rightarrow \infty }
   \frac{N_{\varepsilon_j}}{n_{\varepsilon_j}} = A_\sigma$. 
   \item[(ii)] For any $g\in \cb_+(\R_+)$ bounded continuous, we have 
\[
\lim_{j\rightarrow \infty }
  \inv{n_{\varepsilon_j}}\sum_{k=2}^{N_{\varepsilon_j}}
  g(A_{S^{\varepsilon_j}_k}) =  \int_0^{\tilde
  \sigma} g(u) \, du.
\]
\end{itemize}
\end{lem}

\begin{proof}

  Notice that  $\{s; m_s\neq 0\}  \subset O \cup \{s;  \Delta_s\neq 0\}$
  (see  proof  of Lemma  \ref{lem:O}).  As  $\{s;  \Delta_s\neq 0\}$  is
  discrete, we  have thanks to \reff{eq:cv_st}, that  $\N$-a.e.  for all
  $s\geq   0$,    $\lim_{\varepsilon\to   0}A_s^\varepsilon=A_s$   where
  $A^\varepsilon$ is defined by \reff{eq:def-Ae}. {F}rom Dini Theorem this
  convergence is  uniform on $[0,\sigma]$ $\N$-a.e.  In particular, (ii)
  will be  proved once we  proved (i) and  that $\N$-a.e. for  any $g\in
  \cb_+(\R_+)$ bounded continuous, we have
\begin{equation}
    \label{eq:cv-ke}
\lim_{j\rightarrow \infty }
  \inv{n_{\varepsilon_j}}\sum_{k=2}^{N_{\varepsilon_j}}
  g(A^{\varepsilon_j}_{S^{\varepsilon_j}_k}) =  \int_0^{\tilde
  \sigma} g(u) \, du.
\end{equation}

{{F}rom} Remark  \ref{rem:e_p}, we see  there exists a sequence  of random
variables  $(e_k^\varepsilon,  k\geq 2)$,  such  that conditionally  on
$N_\varepsilon$, they are independent exponential variables of parameter
$n_\varepsilon$ (see \reff{eq:n_e}) and 
\[
A_\sigma^\varepsilon=(S^\varepsilon_1 -
T^\varepsilon_0)+  \sum_{k=2}^{N_\varepsilon} e_k^\varepsilon +
(S^\varepsilon_{N_\varepsilon+1} - T^\varepsilon_{N_\varepsilon}).
\]
We set $e_0^\varepsilon=S^\varepsilon_{N_\varepsilon+1} -
T^\varepsilon_{N_\varepsilon}$ and $ e_1^\varepsilon=
S^\varepsilon_1 -
T^\varepsilon_0$, so that we have the compact notation
$A_\sigma^\varepsilon=\sum_{k=0}^{N_\varepsilon} e_k^\varepsilon $ and
$A^\varepsilon_{S^\varepsilon_k} =\sum_{l=1}^k e^\varepsilon_l$ for
$k\leq N_\varepsilon$.

Because of Lemma  \ref{lem:A_s=0} (ii) and (iii) we  have that $\N$-a.e.
$\lim_{\varepsilon  \downarrow   0}  e^\varepsilon_0  =\lim_{\varepsilon
\downarrow    0}    e^\varepsilon_1=0$.    We deduce that  $\N$-a.e.
 \[
\lim_{\varepsilon\downarrow
  0}\sum_{k=2}^{N_\varepsilon}e_k^\varepsilon=\lim_{\varepsilon\downarrow
  0}A_\sigma^\varepsilon=A_\sigma.
\]

Conditionally on $N_\varepsilon$,  the random variables $(n_\varepsilon
e_k^\varepsilon,  k\geq 2)  $ are  independent exponential  variables of
parameter 1. The previous equality and the law of large numbers implies
that 
$\N$-a.e. for any
positive deterministic sequence $(\varepsilon_j, j\in
\N )$ that decreases to 0, and we obtain (i).

To get \reff{eq:cv-ke}, we choose the  sequence $(\varepsilon_j, j\in
\N )$ so that
for some $\delta\in(0,1/3)$, we have 
\[
\sum_{j=1}^{+\infty}n_{\varepsilon_j}^{-(1-3\delta)/2}<+\infty.
\]
As a consequence of (i), there exists a (random) integer $J$ such that,
if $j\ge J$, 
\[
N_{\varepsilon_j}\le n_{\varepsilon_j}^{1+\delta}.
\]
Notice  that to  prove \reff{eq:cv-ke},  it  is enough  to consider  $g$
bounded and Lipschitz. We have for $j\geq J$,
\begin{align*}
\left|
 \inv{n_{\varepsilon_j}}\sum_{k=2}^{N_{\varepsilon_j}}
  g(A^{\varepsilon_j}_{S^{\varepsilon_j}_k})
-\frac{1}{n_{\varepsilon_j}}
  \sum_{k=2}^{N_{\varepsilon_j}}g\left(\frac{k}{
  n_{\varepsilon_j}}\right)\right|  
& \le
 C_g\frac{1}{n_{\varepsilon_j}} \sum_{k=2}^{N_{\varepsilon_j}}
 \left|\sum_{l=2}^ke_l^{\varepsilon_j}-\frac{k-1}{n_{\varepsilon_j}}\right|  
 + C_g \frac  {N_{\varepsilon_j}}{n_{\varepsilon_j}}
 (e_1^{\varepsilon_j} + \inv{n_{\varepsilon_j}}) \\
& \le C_gZ(\varepsilon_j)+ C_g \frac  {N_{\varepsilon_j}}{n_{\varepsilon_j}}
 (e_1^{\varepsilon_j} + \inv{n_{\varepsilon_j}})
\end{align*}
where $C_g$ is the Lipschitz constant of $g$ and 
\[
Z(\varepsilon)=
\frac{1}{n_\varepsilon}\sum_{k=2}^{n_\varepsilon^{1+\delta}}
\left|\sum_{l=2}^ke_l^{\varepsilon}-\frac{(k-1)}{n_{\varepsilon}}\right|. 
\]
In order to prove that $\lim_{j\rightarrow \infty } Z(\varepsilon_j)=0$,
we compute the expectation of $Z(\varepsilon)$:
\[
\E\bigl[Z(\varepsilon)\bigr] 
 =\frac{1}{n_\varepsilon}\sum_{k=2}^{n_\varepsilon^{1+\delta}} 
\E\left[\left|\sum_{l=1}^ke_l^\varepsilon-\frac{k-1}{n_\varepsilon}\right|\right]
 =\frac{1}{n_\varepsilon^2}\sum_{k=2}^{n_\varepsilon^{1+\delta}}\E\left[\left|\sum_{l=2}^kn_\varepsilon e_l^\varepsilon-(k-1)\right|\right].
\]
But, as the law $n_\varepsilon e_l^\varepsilon$ is the exponential law
with parameter 1, we have 
$$
\E\left[\left(\sum_{l=2}^kn_\varepsilon
    e_l^\varepsilon-(k-1)\right)^4\right] =6k(k-1). 
$$
Thus, the quantity  $\E\bigl[Z(\varepsilon)\bigr]$ is bounded from above
by 
\[
\frac{1}{n_\varepsilon^2}\sum_{k=2}^{n_\varepsilon^{1+\delta}}
\E\left[\left(\sum_{l=2}^kn_\varepsilon 
    e_l^\varepsilon-(k-1)\right)^4\right]^{1/4} 
 \le
2\frac{1}{n_\varepsilon^2}\sum_{k=2}^{n_\varepsilon^{1+\delta}}\sqrt{k  }\\ 
 \le 2n_\varepsilon^{3(1+\delta)/2-2}\\
 \le 2n_\varepsilon^{-(1-3\delta)/2}.
\]
In particular, the series $\sum_{j\geq 1} \E[Z(\varepsilon_j)]$ converges
  and  as  $Z(\varepsilon) $  is  non-negative,  this  implies the  series
  $\sum_{j\geq 1} Z(\varepsilon_j)$ converges a.s. and thus $\N$-a.e. we
  have
\[
\lim_{j\to+\infty}Z(\varepsilon_j)=0.
\]
The convergence of the Riemann's sums gives that $\N$-a.e. 
\[
\frac{1}{n_{\varepsilon_j}}\sum_{k=2}^{N_{\varepsilon_j}}g
\left(\frac{k}{n_{\varepsilon_j}}\right) 
 =\frac{N_{\varepsilon_j}}{n_{\varepsilon_j}}
\frac{1}{N_{\varepsilon_j}}\sum_{k=2}^{N_{\varepsilon_j}}g
\left(\frac{N_{\varepsilon_j}}{n_{\varepsilon_j}}
  \frac{k}{N_{\varepsilon_j}}\right)
 \underset{j\to+\infty}{\longrightarrow}A_\sigma\int_0^1g(uA_\sigma)\,
du=\int _0^{A_\sigma}g(u)\, du.
\]
Then we deduce \reff{eq:cv-ke}  from \reff{eq:s-tilde}, and this
finishes the proof. 
\end{proof}

\subsection{Computation of the limit}
\label{sec:comp_lim}

\begin{lem}
\label{lem:comp_lim}
We assume  ($h_2$), that  is $t\mapsto \phi(t,\omega,\cs)$  is continuous  for all
$\cs\in \S$.   We have, for  the sequence $(\varepsilon_j, j\in  \N^*) $
from Lemma \ref{lem:cv-Ne}, that $\N$-a.e.
\[
   \lim_{j\rightarrow \infty
   }\prod_{k=2}^{N_{\varepsilon_j}}\rN\left[\expp{-\phi(
   A_{S_k^{\varepsilon_j}},\omega, \cdot)} 
   \Bigm|\rho_0>\varepsilon_j\right]
= \exp\left(-\int_0^{\tilde   \sigma}
    du\,\rN\left[1-\expp{-\phi(
   u,\omega, \cdot)}\right]\right).
\]
\end{lem}

\begin{proof}
For any    sequence $(\varphi_k, k\in \N)$ of non-negative
measurable 
function on $\S$, such that $\varphi_k(\cs)=0$ if $\langle \rho_0,1\rangle \leq \eta$,
we have for 
$\varepsilon\in (0,\eta)$, 
\[
\prod_{k=2}^{N_\varepsilon} 
 \rN\left[\expp{-\varphi_k}\Bigm|\rho_0>\varepsilon\right]
 =\prod_{k=2}^{N_\varepsilon}\left(1-\frac
 {\rN\left[1-\expp{-\varphi_k}\right]}{\rN[\rho_0>\varepsilon]}\right). 
\]
Recall          \reff{eq:n_e},          and         notice          that
$\rN\left[1-\expp{-\varphi_k}\right]\leq  \rN[\rho_0\geq    \eta]
\leq    \rN[\rho_0>  \varepsilon]=
n_\varepsilon$      and       $\lim_{\varepsilon      \downarrow      0}
n_\varepsilon=+\infty $.
Since $\log (1-x)=-x+h(x)$, with $\val{h(x)}\leq x^2$ for $x\in
[0,1/2]$, we have  for $\varepsilon$ small enough such that
$\rN[\rho_0>\eta]/n_\varepsilon\leq 1/2$, 
\begin{align*}
\prod_{k=2}^{N_\varepsilon}\left(1-\frac{\rN\left[1-\expp{-\varphi_k}\right]}
  {n_\varepsilon}\right) 
& =\exp\left(\sum_{k=2}^{N_\varepsilon}
 \ln\left(1-\frac{\rN\left[1-\expp{-\varphi_k}\right]}
 {n_\varepsilon}\right)\right)\\  
&= \exp\left(-\frac{1}{n_\varepsilon}
   \sum_{k=2}^{N_\varepsilon}
   \rN\left[1-\expp{-\varphi_k}\right]\right) \exp\left(
   \sum_{k=2}^{N_\varepsilon}
   h(\rN[1-\expp{-\varphi_k}]/ n_\varepsilon) \right),
\end{align*}
and $ \sum_{k=2}^{N_\varepsilon}
   h(\rN\left[1-\expp{-\varphi_k}\right]/ n_\varepsilon) \leq
   \rN[\rho_0>\eta] ^2  N_\varepsilon/ n_\varepsilon^2$. {F}rom the
   hypothesis on $\phi$, we can take
   $\varphi_k=\phi(A_{S^{\varepsilon_j}_k},\omega, \cdot)$. Then,  we deduce
   from Lemma 
   \ref{lem:cv-Ne} (i),  that
   $\N$-a.e. 
\[
\lim_{j\rightarrow \infty } \sum_{k=2}^{N_{\varepsilon_j}}
   h\Big(\rN\Big[1-\expp{-\phi(
   A_{S_k^{\varepsilon_j}},\omega, \cdot)}\Big]   / n_{\varepsilon_j} \Big) =0.
\]
Since ($h_2$) is satisfied,  we deduce that $t\mapsto
\rN\left[1-\expp{-\phi(t,\omega, \cdot)}\right]$ is
continuous. We get from Lemma 
   \ref{lem:cv-Ne} (ii), that
   $\N$-a.e. 
\[
\lim_{j\rightarrow \infty }    \frac{1}{n_{\varepsilon_j}}
   \sum_{k=2}^{N_{\varepsilon_j}}
   \rN\left[1-\expp{-\phi(
   A_{S_k^{\varepsilon_j}},\omega, \cdot)}\right] =\int_0^{\tilde
   \sigma}\rN\left[1-\expp{-\phi(u,\omega, \cdot)}\right]
  \, du .
\]
This finishes the proof of the Lemma. 
\end{proof}

\subsection{Proof of Theorem \ref{th:SMP}}
\label{sec:proof_th_MS}

Let $Z\in \tilde \cf_\infty $ non-negative such that $\N[Z]<\infty $. Let $\phi$
satisfying hypothesis of Theorem  \ref{th:SMP}, ($h_1$) and ($h_2$). We
have 
\begin{align*}
   \N\left[ Z\exp\left(-\sum_{i\in
    I}\phi(A_{\alpha_i},\omega,\cs^i)\right) \right]
& =\lim_{j\rightarrow\infty }    \N\left[
    Z\exp\left(-\sum_{k=2}^{N_{\varepsilon_j}} 
  \phi(A_{S_k^{\varepsilon_j}},\omega,\cs^{k,{\varepsilon_j}})\right)
    \right]\\
& =\lim_{j\rightarrow\infty }
\N\left[Z\prod_{k=2}^{N_{\varepsilon_j}}
\rN\left[\expp{-\phi(  A_{S_k^{\varepsilon_j}},\omega, \cdot)}\Bigm|\rho_0>\varepsilon_j\right]\right]\\
& = \N\left[Z \exp\left(- \int_0^{\tilde
   \sigma}\rN\left[1-\expp{-\phi(u, \omega,\cdot)}\right]
  \, du \right)
    \right],
\end{align*}
where we  used Lemma \ref{lem:approx} and dominated  convergence for the
first  equality,  Lemma \ref{lem:ccea}  for  the  second equality,  Lemma
\ref{lem:comp_lim} and dominated convergence for the last equality. By
monotone class Theorem (resp. monotonicity), we can remove hypothesis
($h_2$) (resp. ($h_1$)). To ends the proof, it suffices to remark that $\exp\left(- \int_0^{\tilde
   \sigma}\rN\left[1-\expp{-\phi(u, \omega,\cdot)}\right]
  \, du \right)$ is $\tilde \cf_\infty$-measurable and so this is
$\N$-a.e. equal to the conditional expectation (i.e. the left hand side
term of \reff{eq:MS}).

\section{Law of the pruned exploration process}
\label{sec:law_pruned}
Recall    notations   of    Section    \ref{sec:LPS}   and    definition
\reff{eq:def_m}.    We  still   fix   $\theta>0$  and   write  $m$   for
$m^{(\theta)}$.
Notice  that  $\psi^{(\theta)}=\psi(\theta+\cdot)-\psi(\theta)$, defined
by \reff{eq:def_psi-q} is  the Laplace exponent of a  Lévy process, with
Lévy measure satisfying  \reff{eq:cond_pi}. The
exploration  process, $\rho^{(\theta)}$,  of this  Lévy process  is thus
well  defined.   

The  aim of this section is to prove the following Theorem.
\begin{theo}\label{thm:law_pruned}
For every finite measure $\mu$, the law of the pruned process $\tilde
\rho$ under $\P_{\mu,0}$ is the law of the exploration process
$\rho^{(\theta)}$ associated to a L\'evy process with Laplace exponent
$\psi^{(\theta)}$ under $\P_\mu$.
\end{theo}
The next Corollary is a direct consequence of this Theorem. 
\begin{cor}
\label{cor:Npruned}
The  excursion measure  of $\tilde  \rho$ outside  $0$ is  equal  to the
excursion measure of $\rho^{(\theta)}$ outside $0$.
\end{cor}

\subsection{A martingale problem for $\tilde \rho$}
In this  section, we  shall compute  the law of  the total  mass process
$(\langle \tilde\rho_{t\wedge \tilde  \sigma} ,1\rangle,\ t\ge 0)$ under
$\P_\mu=\P_{\mu,0}$, using  martingale problem characterization.  We will
first show how a martingale problem  for $\rho$ can be translated into a
martingale problem for $\tilde \rho$.  (In a forthcoming paper, we shall
compute   the  infinitesimal   generator  of   $\rho$   for  exponential
functionals.)   Unfortunately,   we  were  not  able   to  use  standard
techniques  of  random  time  change,  as developed  in  Chapter  6  of
\cite{ek:mp} and used for Poisson snake in \cite{as:psf}, mainly because
$\displaystyle  t^{-1}  \left[\E_{\mu}  [  f(\rho_t)\ind_{\{m_t=0\}}]  -
  f(\mu)\right]$ does not  have a limit as $t$ goes down  to 0, even for
exponential functionals.

Let  $F,  K\in  \cb(\cm_f(\R_+))$  bounded  such  that,  for  any  $\mu\in
\cm_f(R_+)$,  $\displaystyle  \E_\mu\left[\int_0^\sigma \val{K(\rho_s)}\;
  ds  \right]<\infty  $  and   $M_t=F(\rho_{t\wedge  \sigma})
-\int_0^{t\wedge   \sigma}  K(\rho_s)$,   for  $t\geq   0$,   define  an
$\cf$-martingale.    In particular, notice     that    $\E_\mu\left[\sup_{t\geq    0}
  \val{M_t}\right]<\infty $. Thus, we can define for $t\geq 0$,
\[
N_t=\E^*_\mu[M_{C_t}|\tilde \cf_t].
\]
\begin{prop}
\label{prop:mart-tr}
 The process  $N=(N_t, t\geq 0)$ is an $\tilde \cf$-martingale.
And we have for all $\mu\in \cm_f(\R_+)$, $\P_{\mu}$-a.s.
\[
 \int_0^{\tilde   \sigma} du
 \int_{(0,\infty )} \left(1-\expp{-\theta \ell}
\right)\;\pi(d\ell)\val{F([\tilde \rho_u, \ell\delta_0])- F(\tilde
  \rho_u)} <\infty ,
\]
and the representation formula for $N_t$:
\begin{equation}
   \label{eq:Mart-tr}
N_t=F(\tilde \rho_{t\wedge \tilde \sigma})  - 
\int_0^{t\wedge  \tilde \sigma} du\; \left( K(\tilde \rho_u) + 
\int_{(0,\infty )} \left(1-\expp{-\theta \ell}
\right)\;\pi(d\ell)\Big(F([\tilde \rho_u, \ell\delta_0])- F(\tilde
\rho_u)\Big)\right).
\end{equation}
\end{prop}

\begin{proof}
  Notice that $N=(N_t, t\geq 0)$ is an $\tilde \cf$-martingale. Indeed, we have
for $t,s\geq 0$, 
\begin{align*}
   \E_{\mu}[N_{t+s}|\tilde \cf_t]
&=\E_{\mu}[\E_{\mu}[M_{C_{t+s}}|\tilde \cf_{t+s}]|\tilde \cf_t]\\
&=\E_{\mu}[M_{C_{t+s}}|\tilde \cf_t]\\
&=\E_{\mu}[\E_{\mu}[M_{C_{t+s}}|\cf_{C_t}]|\tilde \cf_t]\\
&=\E_{\mu}[M_{C_{t}}|\tilde \cf_t],
\end{align*}
where we used the stopping time Theorem for the last equality.
To compute $\E_{\mu}[M_{C_t}|\tilde \cf_t]$, we write
$M_{C_t}= N'_t - 
M'_{C_t} $, where for $u\geq 0$, 
\[
M'_{u}=\int_0^{u\wedge \sigma} K(\rho_s)\ind_{\{m_s\neq 0\}} \; ds.
\]

Recall that $C_0=0$ $\P_{\mu}$-a.s. by Corollary \ref{cor:C_0=0}.
In particular, we get
\begin{align*}
   N'_t
&=F(\rho_{C_t\wedge \sigma})  - 
\int_0^{C_t\wedge \sigma} K(\rho_s) \ind_{\{m_s= 0\}}  \; ds \\
&=F(\tilde \rho_{t\wedge \tilde \sigma})  - 
\int_0^{C_t\wedge \sigma} K(\rho_s)  \; dA_s \\
&=F(\tilde \rho_{t\wedge \tilde \sigma})  - 
\int_0^{t\wedge  \tilde \sigma} K(\tilde \rho_u)  \; du ,
\end{align*}
where  we used  the  time change  $u=A_s$  for the  last  equality.   In
particular, as  $\tilde \sigma $ is  an $\tilde \cf$-stopping time,  we get that
the process $(N'_t, t\geq 0)$  is $\tilde\cf$-adapted.  Since  $N_t=N'_t
-  \E_{\mu}[M'_{C_t}|\tilde \cf_t]  $, we  are left  with  the  computation  of
$\E_{\mu}[M'_{C_t}|\tilde \cf_t] $.

In Section  \ref{sec:Markov_special}, the arguments are  given under the
excursion measure,  but they can readily be  extended under $\P_{\mu}$
or   $\P^*_{\mu,0}$.   In    particular,   the   result   of   Corollary
\ref{cor:super} holds also under $\P_{\mu}$ or $\P^*_{\mu,0}$.
We keep  the notations of Section  \ref{sec:Markov_special}. We consider
$(\rho^i,  m^i)$, $i\in  I$  the excursions  of  the process  $(\rho,m)$
outside  $\{s,m_s=0\}$ before  $\sigma$ and  let  $(\alpha_i, \beta_i)$,
$i\in I$ be the corresponding  interval excursions. In particular we can
write
\[
\int_0^{C_t\wedge \sigma} \val{K(\rho_s)} \ind_{\{m_s\neq 0\}}\;
ds  =\sum_{i\in    I}     \Phi( A_{\alpha_i},
\rho_{\alpha_i-},\rho_i) ,
\]
with
\[
\Phi(u,\mu,\rho)=\ind_{\{u<t\}} \int_0^{\sigma(\rho)}  \val{K([\mu, \rho_s])} \; ds,
\]
where $\sigma(\rho)=\inf\{v>0; \rho_v=0\}$. 
We deduce from   Corollary  \ref{cor:super},  that $\P_\mu$-a.s. 
\begin{equation}
    \label{eq:K-hat}
\E_{\mu}\left[\int_0^{C_t\wedge \sigma} \val{K(\rho_s)} \ind_{\{m_s\neq 0\}}\;
ds|\tilde \cf_\infty   \right]  =   \int_0^{\tilde   \sigma} 
\ind_{\{u<t\}} \hat K(\tilde \rho_u) \; du ,
\end{equation}
with, $\hat K$ defined for $\nu\in \cm_f(\R_+)$ by  
\[
\hat K(\nu) 
=
\int_{(0,\infty )} \left(1-\expp{-\theta \ell} \right)\;\pi(d\ell)
\;\E_\ell 
\left[\int_0^\sigma \val{K([\nu, \rho_s])} \; ds\right]. 
\]
Since $\E_{\mu}\left[\int_0^{ \sigma} \val{K(\rho_s)}  \; 
ds\right]$ is finite, we deduce that $\P_{\mu}$-a.s. $du$-a.e.
$\ind_{\{u<\tilde \sigma\}}\hat K(\tilde
\rho_u) $ is finite. 

We  define $\tilde  K\in \cb(\cm_f(\R_+))$  for $\nu\in  \cm_f(\R_+)$ by
\begin{equation}
   \label{eq:def-tK}
\tilde  K(\nu) 
=
\int_{(0,\infty )} \left(1-\expp{-\theta \ell} \right)\;\pi(d\ell)
\;\E_\ell 
\left[\int_0^\sigma {K([\nu, \rho_s])} \; ds\right],  
\end{equation}
if   $\hat  K(\nu)<\infty   $,  or   by  $\tilde   K(\nu)=0$   if  $\hat
K(\nu)=+\infty  $.  In  particular,  we have  $|\tilde K(\nu)|\leq  \hat
K(\nu)$  and $\P_{\mu}$-a.s.   $\int_0^{\tilde \sigma}  |\tilde K(\tilde
\rho_u)| \;  du$ is finite.  Using Corollary \ref{cor:super}  once again
(see \reff{eq:K-hat}), we get that $\P_\mu$-a.s.,
\begin{equation}
   \label{eq:EMC-tK}
\E_{\mu}\left[M'_{C_t} |\tilde \cf_\infty   \right]  
= \E_{\mu}\left[\int_0^{C_t\wedge \sigma} {K(\rho_s)} \ind_{\{m_s\neq 0\}}\;
ds|\tilde \cf_\infty   \right]  =   \int_0^{t\wedge\tilde   \sigma} 
 \tilde K(\tilde \rho_u) \; du. 
\end{equation}
To rewrite $\tilde K$, we notice that, for $\nu$ with compact support,
$\displaystyle  \E _\ell 
\left[\int_0^\sigma {K([\nu, \rho_s])} \; ds\right]$ is equal to 
$\E_{[\nu, \ell\delta_0]} 
\left[\int_0^{\tau_\ell}  {K(\rho_s)} \; ds\right]$,
where   $\tau_\ell=\inf\{s;  -I_s\geq   \ell\}$   is  an  $\cf$-stopping
time. Notice that  $\P_{[\nu, \ell \delta_0]}$-a.s.   $\tau_\ell\leq \sigma$
and $\rho_{\tau_\ell}=\nu$. We deduce  from the stopping time Theorem
that 
\begin{equation}
 \label{eq:Emul-K}
\E_{[\nu, \ell\delta_0]} 
\left[\int_0^{\tau_\ell}  {K(\rho_s)} \; ds\right]
=\E_{[\nu, \ell\delta_0]} 
\left[-M_{\tau_\ell} + F(\rho_{\tau_\ell})\right]
=-F([\nu, \ell\delta_0])+ F(\nu).      
\end{equation} 
Therefore, we get from \reff{eq:def-tK} and \reff{eq:EMC-tK}
\[
\E_{\mu}\left[M'_{C_t} |\tilde \cf_\infty   \right]  
=   -\int_0^{t\wedge \tilde   \sigma} 
 \int_{(0,\infty )} \left(1-\expp{-\theta \ell}
\right)\;\pi(d\ell)\Big(F([\tilde \rho_u, \ell\delta_0])- F(\tilde
\rho_u)\Big)\; du. 
\]
Eventually, as $N_t
=N'_t- \E_{\mu}\left[M'_{C_t} |\tilde \cf_\infty   \right] $,  this
gives \reff{eq:Mart-tr}.

To conclude, notice that from \reff{eq:Emul-K}, the definition of $\hat K$  and
\reff{eq:K-hat}, we have 
\begin{multline*}
     \int_0^{\tilde   \sigma} 
 \int_{(0,\infty )} \left(1-\expp{-\theta \ell}
\right)\;\pi(d\ell)\val{F([\tilde \rho_u, \ell\delta_0])- F(\tilde
  \rho_u)}\; du\\ 
\begin{aligned}
&\leq   \int_0^{\tilde   \sigma} 
 \int_{(0,\infty )} \left(1-\expp{-\theta \ell}
\right)\;\pi(d\ell)\E_{[\tilde \rho_u, \ell\delta_0]} 
\left[\int_0^{\tau_\ell}  \val {K(\rho_s)} \; ds\right]\; du\\
&= \int_0^{\tilde   \sigma} \hat K(\tilde \rho_u)\; du\\
&=
\E_{\mu}\left[\int_0^{ \sigma} \val{K(\rho_s)} \ind_{\{m_s\neq 0\}}\;
ds|\tilde \cf_\infty   \right] ,
\end{aligned}
\end{multline*}
 which is finite $\P_{\mu}$-a.s. since
$\displaystyle  \E_{\mu}\left[\int_0^{ \sigma} \val{K(\rho_s)} \;
ds   \right] <\infty $. 
\end{proof}

\begin{cor}\label{cor:law_total_mass}
  Let $\mu\in \cm_f(\R_+)$. The law  of the total mass process $(\langle
  \tilde\rho_t,1\rangle,\ t\ge  0)$ under  $\P_{\mu,0}^*$ is the  law of
  the total mass process of $\rho^{(\theta)}$ under $\P_\mu^*$.
\end{cor}

\begin{proof}
  Let $X=(X_t, t\geq 0)$ be under $\rP^*_x$, a Lévy process with Laplace
  transform $\psi$  started at  $x>0$ and stopped  when it  reached $0$.
  Under  $\P_\mu$,  the   total  mass  process  $(\langle  \rho_{t\wedge
    \sigma},  1\rangle,   t  \geq  0)$  is  distributed   as  $X$  under
  $\rP^*_{\langle \mu,  1\rangle} $.  Let  $c> 0$.  From  Lévy processes
  theory, we  know that  the process $\expp{-c  X_t} -  \psi(c) \int_0^t
  \expp{-cX_s}\; ds $, for $t\geq 0$ is a martingale. We deduce from the
  stopping time  Theorem that $M=(M_t, t\geq 0)$  is an $\cf$-martingale
  under  $\P_\mu$,  where $M_t=F(\rho_{t\wedge  \sigma})-\int_0^{t\wedge
    \sigma} K(\rho_s) \; ds $,  with $F, K\in \cb( \cm_f(\R_+))$ defined
  by  $F(\nu)=\expp{   -  c  \langle   \nu,  1\rangle}$  for   $\nu  \in
  \cm_f(\R_+)$ and $K=\psi(c)F$. Notice  $K\geq 0$. We have by dominated
  convergence and monotone convergence.
\[
\expp{-   c\langle   \mu,   1\rangle}   =\lim_{t\rightarrow   \infty   }
\E_\mu[M_t]=\E_\mu[\expp{  - c  \langle \rho_\sigma  ,  1\rangle}] -
\psi(c)  \E_\mu  \left[\int_0^\sigma  \expp{  -  c  \langle  \rho_s  ,
    1\rangle} \; ds \right] .
\]
This implies that, for any $\mu\in
\cm_f(R_+)$,  $\displaystyle  \E_\mu\left[\int_0^\sigma \val{K(\rho_s)}\;
  ds  \right]$  is finite. 
For $\nu\in \cm_f(\R_+)$ with compact support, we have 
\begin{multline*}
    \int_{(0,\infty )} \left(1-\expp{-\theta \ell}
\right)\;\pi(d\ell)\val{F([\nu, \ell\delta_0])- F(\nu)}\\
\begin{aligned}
   &=\int_{(0,\infty )} \left(1-\expp{-\theta \ell}
\right)\;\pi(d\ell)\left(\expp{ - c\langle   \nu,   1\rangle} - \expp{ -
    c\langle   \nu,   1\rangle - c \ell}\right)\\
&=\expp{ - c\langle   \nu,   1\rangle} \int_{(0,\infty )} \left(1-\expp{-\theta \ell}
\right)\left(1  - \expp{ - c \ell}\right)\;\pi(d\ell)\\
&= \expp{ - c\langle   \nu,   1\rangle} \left(\psi(c) -
  \psi^{(\theta)}(c)\right). 
\end{aligned}
\end{multline*}
In particular, we have 
\[
    \int_{(0,\infty )} \left(1-\expp{-\theta \ell}
\right)\;\pi(d\ell)\Big[ F([\tilde \rho_u, \ell\delta_0])- F(\tilde
  \rho_u)\Big] 
= \expp{ - c\langle   \tilde \rho _u,   1\rangle} \left(
  \psi^{(\theta)}(c)-\psi(c) \right).
\]
{F}rom Proposition \ref{prop:mart-tr}, we get that $N=(N_t, t\geq 0)$,
with for $t\geq 0$,  
\[
N_t=\expp{-c   \langle \tilde\rho_{t\wedge \tilde \sigma}  ,  1\rangle} -
\psi^{(\theta)}(c) \int_0^{t\wedge \tilde \sigma} \expp{-c   \langle
  \tilde\rho_s  ,  1\rangle}\; ds , 
\]
is under $\P_\mu$ an $\tilde \cf$-martingale. 

Notice  that $\tilde  \sigma=\inf  \{ s\geq  0;  \langle \tilde\rho_s  ,
1\rangle=0\}$.  Let $X^{(\theta)}=(X^{(\theta)}_t ,  t\geq 0)$  be under
$\rP^*_x$,  a  Lévy  process  with Laplace  transform  $\psi^{(\theta)}$
started at $x>0$ and stopped when it reached $0$.
The two non-negative  càd-làg processes $( \langle
  \tilde\rho_{t\wedge \tilde \sigma}   ,  1\rangle, t\geq 0)$ and 
$X^{(\theta)}$ solves the martingale problem: for any $c\geq 0$, the
process defined for $t\geq 0$ by 
\[
\expp{- c Y_{t\wedge \sigma'}} - \psi^{(\theta)}(c) \int_0^{t\wedge \sigma'} \expp{- c Y_s}\; ds,
\]
where $\sigma'=\inf  \{ s\geq  0;  Y_s\leq 0\}$, is a martingale. {F}rom
Corollary 4.4.4 in  \cite{ek:mp}, we deduce that those two processes have
the same distribution. To finish the proof, notice that the total mass
process of $\rho^{(\theta)}$ under $\P_\mu^*$ is distributed as
$X^{(\theta)}$ under $\rP^*_{\langle \mu,  1\rangle} $.
\end{proof}

\subsection{Identification of the law of $\tilde \rho$}

To begin with, let us mention some useful properties of the process $\tilde\rho$.
\begin{lem}
   We have the following properties for the process $\tilde \rho$. 
\begin{enumerate}
   \item[(i)] $\tilde \rho$ is a càd-làg Markov process.
   \item[(ii)] The sojourn time at $0$ of $\tilde \rho$ is  $0$.
   \item[(iii)] 0 is recurrent for $\tilde \rho$.
\end{enumerate}
\end{lem}

\begin{proof}
   
(i) This is a direct consequence of the strong Markov property of the
    process $(\rho, m)$.

(ii) We have for $r>0$, with the change of variable $t=A_s$, a.s. 
\[
\int_0^r \ind_{\{\tilde \rho _t=0\}} \; dt
=\int_0^r \ind_{\{ \rho _{C_t}=0\}} \; dt
=\int_0^{C_r} \ind_{\{\rho _s=0\}} \; dA_s
=\int_0^{C_r} \ind_{\{\rho _s=0\}} \; ds=0,
\] 
as the sojourn time of $\rho$ at $0$ is $0$ a.s.

(iii) Since $ \tilde  \sigma =A_\sigma$ and $\sigma<+\infty $ a.s., we
    deduce that $0$ is recurrent for $\tilde \rho$ a.s. 
\end{proof}

Since the processes $\tilde \rho$ and $\rho^{(\theta)}$ are both
Markov processes, to show that they have the same law, it is enough to
show that they have the same one-dimensional marginals. We first prove
that result under the excursion measure.

\begin{prop}\label{prop:=n}
For every $\lambda>0$ and every non-negative bounded measurable
function $f$,
$$\N\left[\int_0^{\tilde\sigma}\expp{-\lambda t-\langle \tilde\rho_t,f\rangle}dt\right]=\N\left[\int_0^{\sigma^{(\theta)}}\expp{-\lambda t-\langle \rho^{(\theta)}_t,f\rangle}dt\right].$$
\end{prop}

\begin{proof}
On one hand, we compute, using the definition of the pruned process
$\tilde \rho$,
$$
\N\left[\int_0^{\tilde\sigma}\expp{-\lambda t-\langle
    \tilde\rho_t,f\rangle}dt\right]
=\N\left[\int_0^{A_\sigma}\expp{-\lambda t-\langle \rho_{C_t},f\rangle}dt\right].
$$
We now make the change of variable $t=A_u$ to get
\begin{align*}
\N\left[\int_0^{\tilde\sigma}\expp{-\lambda t-\langle
    \tilde\rho_t,f\rangle}dt\right] & =\N\left[\int_0^\sigma
    \expp{-\lambda A_u}\expp{-\langle\rho_u,f\rangle}dA_u\right]\\
& =\N\left[\int_0^\sigma \expp{-\lambda A_u}\expp{-\langle\rho_u,f\rangle}\ind_{\{m_u=0\}}du\right].
\end{align*}

By a time reversibility argument, see Lemma \ref{lem:reversib}, we obtain
\begin{align*}
\N\left[\int_0^{\tilde\sigma}\expp{-\lambda t-\langle
    \tilde\rho_t,f\rangle}dt\right]
& =\N\left[\int_0^\sigma
    \ind_{\{m_u=0\}}\expp{-\langle\eta_u,f\rangle}\expp{-\lambda(A_\sigma-A_u)}du\right]\\
& =\N\left[\int_0^\sigma
    \ind_{\{m_u=0\}}\expp{-\langle\eta_u,f\rangle}\E_{\rho_u,0}^*\left[\expp{-\lambda
    A_\sigma}\right]du\right]\\
& =\N\left[\int_0^\sigma
    \ind_{\{m_u=0\}}\expp{-\langle\eta_u,f\rangle}\expp{-\langle \rho_u,1\rangle{\psi^{(\theta)}}^{-1}(\lambda)}du\right]
\end{align*}
where we applied Lemma \ref{lem:A_s=0} (i) for the last equality. Now,
by definition of $m$, we have by conditioning,
$$\N\left[\int_0^{\tilde\sigma}\expp{-\lambda t-\langle
    \tilde\rho_t,f\rangle}dt\right]=\N\left[\int_0^\sigma
   \expp{-\theta\langle
    \kappa_u,1\rangle}\expp{-\langle\eta_u,f\rangle}\expp{-\langle
    \rho_u,1\rangle{\psi^{(\theta)}}^{-1}(\lambda)}du\right].$$
Now, the Poisson decomposition of Proposition
\ref{prop:poisson_representation1}  and
    standard computations lead to
\begin{align*}
\N & \left[\int_0^{\tilde\sigma}\expp{-\lambda t-\langle
    \tilde\rho_t,f\rangle}dt\right]\\
& =\int_0^{+\infty}da\expp{-\alpha_0
    a}\exp\left\{-\int_0^adx\int_0^1du\int_{(0,+\infty)}\ell\pi(d\ell)\left[1-\expp{-\ell(\theta+(1-u)f(x)+u{\psi^{(\theta)}}^{-1}(\lambda))}\right]\right\}\\
& =\int_0^{+\infty}da\exp\left\{-\int_0^adx\int_0^1du\;
    \psi'\bigl(\theta+(1-u)f(x)+u{\psi^{(\theta)}}^{-1}(\lambda)\bigr)\right\}\\
& =\int_0^{+\infty}da\exp\left\{-\int_0^adx\;\frac{\lambda-\psi^{(\theta)}\bigl(f(x)\bigr)}{{\psi^{(\theta)}}^{-1}(\lambda)-f(x)}\right\}.
\end{align*}

On the other hand, the formula of Proposition
\ref{prop:poisson_representation}, the Poisson representation of
Proposition \ref{prop:poisson_representation1} and the same
computations as before yields
\begin{align*}
\N & \left[\int_0^{\sigma^{(\theta)}}\expp{-\lambda
    t-\langle\rho^{(\theta)}_t,f\rangle}dt\right]\\
& =\int\mathbb{M}(d\mu\, d\nu)\expp{-\langle
    \mu,f\rangle}\expp{-{\psi^{(\theta)}}^{-1}(\lambda)\langle\nu,1\rangle}\\
& =\int_0^{+\infty}da
    \expp{-\alpha_0a}\exp\left\{-\int_0^adx\int_1du\int_{(0,+\infty)}\ell
    \pi^{(\theta)}(d\ell)\left[1-\expp{-\ell(uf(x)+{\psi^{(\theta)}}^{-1}(\lambda)(1-u))}\right]\right\}\\
& =\int_0^{+\infty}da\exp\left\{-\int_0^adx\frac{\lambda-\psi^{(\theta)}\bigl(f(x)\bigr)}{{\psi^{(\theta)}}^{-1}(\lambda)-f(x)}\right\}.
\end{align*}
As the two quantities are equal, the proof is over.
\end{proof}

Now, we prove the same result under $\P_{\mu,0}^*$, that is:
\begin{prop}\label{prop:Laplace_*}
For every $\lambda>0$,  $f\in \cb_+(\R_+)$ bounded and every finite measure $\mu$,
$$\E_{\mu,0}^*\left[\int_0^{\tilde\sigma}\expp{-\lambda t-\langle \tilde\rho_t,f\rangle}dt\right]=\E_\mu^*\left[\int_0^{\sigma^{(\theta)}}\expp{-\lambda t-\langle \rho^{(\theta)}_t,f\rangle}dt\right].$$
\end{prop}

\begin{proof}
{F}rom the Poisson representation, see Lemma \ref{lem:dlg-decomp}, and
using notations of this Lemma and of \reff{eq:def-Ai} we have 
\begin{align*}
\E_{\mu,0}^*\left[\int_0^{\tilde \sigma}\expp{-\lambda
    t-\langle\tilde\rho_t,f\rangle}dt\right] 
& =\E_{\mu,0}^*\left[\int_0^\sigma \expp{-\lambda A_u-\langle
    \rho_u,f\rangle}dA_u\right]\\
& =\E_{\mu,0}^*\left[\sum_{i\in J}\expp{-\lambda
    A_{\alpha_i}-\langle
    k_{-I_{\alpha_i}},f\rangle}\int_0^{\sigma_i}\expp{-\langle
    \rho_s^i,f_{-I_{\alpha_i}}\rangle-\lambda A^i_s}dA^i_s\right]
\end{align*}
where the function $f_r$ is defined by $f_r(x)=f(H^{(\mu)}_r+x)$ and 
$H^{(\mu)}_r=H(k_r\mu)$ is the maximal
element of the 
closed support of $k_r \mu$ (see \reff{def:H}).
We recall that $-I$ is the local time at 0 of the reflected process 
 $X-I$, and that $\tau_r=\inf\{s; -I_s>r\}$ is the  right continuous inverse of $-I$. {F}rom excursion formula, and using the time
change $-I_s=r$ (or equivalently $\tau_r=s$), we get 
\begin{align}
\nonumber
\E_{\mu,0}^*\left[\int_0^{\tilde \sigma}\expp{-\lambda
    t-\langle\tilde\rho_t,f\rangle}dt\right] 
& =\E_{\mu,0}^*\left[\int_0^{\tau_{\langle\mu,1\rangle}}d(-I_s)\expp{-\langle
    k_{-I_s}\mu,f\rangle-\lambda A_s}G(-I_s)\right]\\
& =\E_{\mu,0}^*\left[\int_0^{\langle \mu,1\rangle} dr\expp{-\langle
    k_r\mu,f\rangle-\lambda A_{\tau_r}}G(r)\right]
\label{eq:EGr1}
\end{align}
where the function $G(r)$ is given by
\[
G(r)=\N\left[\int_0^\sigma\expp{-\langle\rho_s,f_r\rangle-\lambda
    A_s}dA_s\right]=  
\N\left[\int_0^{\tilde\sigma}\expp{-\lambda t-\langle
    \tilde\rho_t,f_r\rangle}dt\right].
\]
The same kind of computation gives
\begin{equation}
   \label{eq:EGr2}
\E_\mu^*\left[\int_0^{\sigma^{(\theta)}}\expp{-\lambda t-\langle
    \rho_t^{(\theta)},f\rangle}dt\right]=\E\left[\int_0^{\langle
    \mu,1\rangle}dr\expp{-\langle
    k_r\mu,f\rangle-\lambda\tau_r^{(\theta)}}G^{(\theta)}(r)\right]
\end{equation}
where the function $G^{(\theta)}$ is defined by
\[
G^{(\theta)}(r)=\N\left[\int_0^{\sigma^{(\theta)}}\expp{-\lambda
    s-\langle\rho_s^{(\theta)},f_r\rangle}ds\right]
\]
and $\tau^{(\theta)}$ is the right-continuous inverse of the infimum
    process $-I^{(\theta)}$ of the L\'evy process with Laplace
    exponent $\psi^{(\theta)}$. 

Proposition \ref{prop:=n} says that the functions $G$ and
$G^{(\theta)}$ are equal. Moreover, as the total mass processes have
the same law (see Corollary \ref{cor:law_total_mass}), we know that
the proposition is true for $f$ constant. And, for $f$ constant, the
functions $G$ and
$G^{(\theta)}$ are also constant. Therefore, we have for $f$ constant
equal to $c\geq 0$,
\[
\E_{\mu,0}^*\left[\int_0^{\langle \mu,1\rangle}dr\expp{-c(\langle
    \mu,1\rangle-r)}\expp{-\lambda                     A_{\tau_r}}\right]
=\E\left[\int_0^{\langle                  \mu,1\rangle}dr\expp{-c(\langle
    \mu,1\rangle-r)}\expp{-\lambda  \tau_r^{(\theta)}}\right].
\]
As this is  true for any $c\geq 0$, uniqueness  of the Laplace transform
gives the equality
\[
\E_{\mu,0}^*\left[\expp{-\lambda
    A_{\tau_r}}\right]=\E\left[\expp{-\lambda\tau_r^{(\theta)}}\right]\qquad
    dr-\mbox{a.e.}
\]
In fact this equality holds for every $r$ by right-continuity.

Eventually as $G=G^{(\theta)}$, we have thanks to \reff{eq:EGr1} and
\reff{eq:EGr2}, that, for every bounded non-negative measurable
function $f$,
$$\int_0^{\langle \mu,1\rangle} dr\expp{-\langle
    k_r\mu,f\rangle}\E_{\mu,0}^*\left[\expp{-\lambda A_{\tau_r}}\right]G(r)=\int_0^{\langle
    \mu,1\rangle}dr\expp{-\langle
    k_r\mu,f\rangle}\E\left[\expp{-\lambda\tau_r^{(\theta)}}\right]G^{(\theta)}(r)$$
which ends the proof.
\end{proof}

\begin{cor}
\label{cor:=*} 
  The  process $\tilde  \rho$ under  $\P^*_{\mu,0}$ is  distributed as
   $\rho^{(\theta)}$ under $\P^*_\mu$.
\end{cor}

\begin{proof}
  Let $f\in  \cb_+(\R_+)$ bounded. Proposition  \ref{prop:Laplace_*} can
  be re-written as
$$\int_0^{+\infty}\expp{-\lambda t}\E_{\mu,0}^*\left[\expp{-\langle
    \tilde\rho_t,f\rangle}\ind_{\{t\le \tilde\sigma\}}\right]dt =\int_0^{+\infty}\expp{-\lambda t}\E_{\mu}^*\left[\expp{-\langle
    \rho^{(\theta)}_t,f\rangle}\ind_{\{t\le\sigma^{(\theta)}\}}\right]dt.$$
By uniqueness of the Laplace transform, we deduce that, for almost every
    $t>0$,
$$\E_{\mu,0}^*\left[\expp{-\langle
    \tilde\rho_t,f\rangle}\ind_{\{t\le
    \tilde\sigma\}}\right]=\E_{\mu}^*\left[\expp{-\langle
    \rho^{(\theta)}_t,f\rangle}\ind_{\{t\le\sigma^{(\theta)}\}}\right].$$
In fact this  equality holds for every $r$  by right-continuity.  As the
Laplace functionals characterize the law  of a random measure, we deduce
that,  for  fixed  $t>0$, the  law of  $\tilde\rho_t$ under
$\P_{\mu,0}^*$  is the  same  as the  law  of $\rho^{(\theta)}_t$  under
$\P_\mu^*$.

The Markov property  then give the equality in law for the  càd-làg processes
$\tilde\rho$ and $\rho^{(\theta)}$.
\end{proof}

\begin{proof}[Proof of Theorem \ref{thm:law_pruned}]
   $0$ is recurrent for the Markov càd-làg processes $\tilde \rho$ and
   $\rho^{(\theta)}$.  This two processes have no sojourn at $0$, and
   when killed on the first hitting time of $0$, they have the
   same law, thanks to Lemma \ref{cor:=*}. {F}rom Theorem 4.2 of
   \cite{b:emp}, Section 5, we deduce that $\tilde \rho$ under
   $\P_{\mu,0}$ is distributed as $\rho^{(\theta)}$ under $\P_\mu$. 
\end{proof}


\section{Property of the excursion of the pruned exploration  process}
\label{sec:prop_pruned}

We know, (cf
\cite{b:pl}, Section VII) that the right continuous inverse,
$(\tau_r,r\ge 0)$, of $-I$
is a subordinator with Laplace exponent
$\psi^{-1}$. This subordinator has no drift as \reff{eq:psi/l} implies
$\lim_{\lambda\rightarrow\infty }\lambda^{-1}{\psi^{-1}(\lambda)}=0$.
We denote by  $\pi_* $ its Lévy measure: for $\lambda\geq 0$
\[
\psi^{-1}(\lambda)= \int _{(0,\infty )} \pi_* (dl)
(1-\expp{-\lambda l}).
\]
Recall $\N$  is the excursion  measure of the exploration  process above
$0$.  If  $\sigma$  denotes  the  duration  of  the  excursion,  we  have
$\N[1-\expp{-\lambda  \sigma}]=\psi^{-1}(\lambda)$.  Hence, under  $\N$,
$\sigma$  is   distributed  according  to  the   measure  $\pi_* $.  By
decomposing the measure $\N$ w.r.t.  the distribution of $\sigma$, we get
that  $\N[d\ce]=\int  _{(0,\infty  )}  \pi_* (dr)  \N_r[d\ce]$,  where
$(\N_r, r \in (0,\infty ))$ is a measurable family of probability
measures on the
set of excursions such that $\N_r[\sigma=r]=1$ for $\pi^*$-a.e. $r>0$.

\begin{lem}\label{lem:loi_cond}
  Conditionally on the length of the excursion, the law of the excursion
  of the pruned  exploration process is the law of  the excursion of the
  exploration process.
\end{lem}

\begin{proof}  {F}rom  the previous  Section,  we  get  that the  pruned
  exploration process $(\tilde\rho_t, t\geq 0)$ is distributed according
  to the  law of the  exploration process, $\rho^{(\theta)}$, of  a Lévy
  process,       $X^{(\theta)}$,       with       Laplace       exponent
  $\psi^{(\theta)}=\psi(\theta+\cdot)-\psi(\theta)$.   In particular the
  law of the  pruned exploration process under the  excursion measure is
  the  law  of  the  exploration  process  $\rho^{(\theta)}$  under  the
  excursion measure.

  Let  $\sigma^{(\theta)}$ denote  the length  of the  excursion  of the
  exploration  process $\rho^{(\theta)}$  under  the excursion  measure.
  The following result is known, but since we couldn't give a reference,
  we shall give a proof at the end of this Section.

\begin{lem}
\label{lem:X_pruned}
For  any non-negative measurable 
  function, $G$,  on the space of excursions, we have 
\[
\N \left[ \expp{ \psi(\theta)\sigma^{(\theta)}} [1- \expp{-G(\rho^{(\theta)})}]
\right]
=\N\left[ 1- \expp{-G(\rho)} \right].
\]   
\end{lem}

In particular the distribution  of $\rho^{(\theta)}$ under the excursion
measure is absolutely continuous  w.r.t. to distribution of $\rho$ under
the   excursion   measure,   with   density   given   by   $\expp{-\sigma
  \psi(\theta)}$.   We deduce that  $\displaystyle \pi^{(\theta)}_*(dr)=
\expp{-r \psi(\theta)} \pi_* (dr)$, where $\pi^{(\theta)}_*$ is the Lévy
measure       corresponding      to      the       Laplace      exponent
$(\psi^{(\theta)})^{-1}$. And we have $\pi_*(dr)$-a.e., conditionally on
the  length  of the  excursion  being  equal to  $r$,  the  law of  the
excursion of the pruned exploration  process is the law of the excursion
of the exploration process.

\end{proof}

Recall  $\tilde\sigma=\int_0^\sigma  \ind_{\{m^{(\theta)}_s=0\}}\; ds $
denotes the length of the excursion of the pruned exploration process. We
can compute  the joint  law of $(\tilde\sigma,\sigma)$.   This will
determine  uniquely  the  law  of $\tilde\sigma$  conditionally  on
$\sigma=r$.

\begin{prop}\label{prop:s-sq}
   For all non-negative $\gamma, \kappa, \theta$, the value  $v$ defined
   by $\displaystyle v=\N\left[1-\expp{-\psi(\gamma) \sigma -\kappa
   \tilde\sigma}\right] $ is the unique non-negative solution of
   the equation 
\[
\psi(v+\theta)= \kappa +\psi(\gamma+\theta).
\]
\end{prop}

\begin{proof}

   Using the special Markov property, Theorem \ref{th:SMP}, with
   $\phi(\cs)=\psi(\gamma)\sigma$, we have
\begin{align*}
v
= \N\left[1-\expp{ -\kappa \tilde\sigma-\psi(\gamma)\sigma}\right]
&= \N\left[1-\expp{ -(\kappa+\psi(\gamma) )
    \tilde\sigma-\psi(\gamma) 
    \int_0^\sigma \ind_{\{m_s\neq  0\}}\; ds}\right]\\
&= \N\left[1-\expp{ -(\kappa+\psi(\gamma))
    \tilde \sigma- \tilde \sigma \int_{(0,+\infty )} \pi(d\ell)
    (1-\expp{-\theta 
      \ell}) \E_\ell^* [1-\exp{( -\psi(\gamma) \sigma)}] } \right].
\end{align*}
Notice that  $\sigma$ under $\P^*_\ell$ is distributed  as $\tau_\ell$, the
first  time  for which    the  infimum  of  $X$,  started   at  $0$,  reaches
$-\ell$. Since  $\tau_\ell$ is distributed as a
subordinator with Laplace exponent $\psi^{-1}$ at time
$\ell$,  we have 
\[
\E_\ell^* [1-\expp{ -\psi(\gamma) \sigma}] 
=\E\left[1 - \expp{-\psi(\gamma) \tau_\ell} \right] = 1 -\expp{-\ell \gamma}.
\]
and 
\[
\int_{(0,+\infty )} \pi(d\ell)
    (1-\expp{-\theta 
      \ell}) \E_\ell^* [1-\expp{ -\psi(\gamma) \sigma}] =
    \int_{(0,+\infty )} \pi(d\ell) 
    (1-\expp{-\theta       \ell})   (1-\expp{-\gamma      \ell})  =
    \psi^{(\theta)} (\gamma) - \psi(\gamma).
\]
We get 
\[
v= \N\left[1-\expp{ -\tilde \sigma (\kappa+\psi^{(\theta)}
    (\gamma))} \right] = {\psi^{(\theta)}}^{-1} (\kappa+\psi^{(\theta)}
    (\gamma)).
\]
Using Corollary \ref{cor:Npruned}  and definition \reff{eq:def_psi-q} of
$\psi^{(\theta)}$,       we      have       $\psi(v+\theta)=      \kappa
+\psi(\gamma+\theta)$.  Since $\psi$  is increasing and continuous, this
equation has only one solution.
\end{proof}

\begin{proof}[Proof of Lemma \ref{lem:X_pruned}]

Since an excursion of the exploration process above $0$ can be recovered
from an excursion of the process $X$ above its minimum. We shall prove
the Lemma  in the latter case. 

   Let  $\theta>0$.  We set  $X^{(\theta)}=(X^{(\theta)}_t,  t\geq 0)$  the
Lévy process with Laplace exponent $\psi^{(\theta)}$. 
Notice that $(\expp{-\theta X_t - t\psi(\theta)},  t\geq 0) $ is a martingale
w.r.t. the natural filtration generated  by $X$, $(\ch_t, t\geq 0)$.  We
define a new probability by
\[
d\P^{(\theta)}_{|\ch_t}= \expp{-\theta X_t - t\psi(\theta)}
d\P_{|\ch_t}.
\]
The law of  $(X_u, u\in [0,t])$ under  $\P^{(\theta)}$ is the law of 
$(X_u^{(\theta)}, u\in [0,t])$. Therefore, we have for any non-negative
measurable function on the path space
\begin{equation}
   \label{eq:cgt-proba}
{\E\left[F(X_{\leq  t}^{(\theta)}) \expp{\theta X_t^{(\theta)} +t
    \psi(\theta)} \right]=\E[F(X_{\leq  t})].}
\end{equation}
We define $-I^{(\theta)}_t =-\inf_{u\in [0,t]}
X^{(\theta)}_u$, and $\tau^{(\theta)}$ its  right-continuous
inverse.  In particular, it is a subordinator of Laplace exponent
${\psi^{(\theta)}}^{-1}$. Since ${\psi^{(\theta)}}^{-1} (\lambda)=
\psi^{-1}(\lambda+\psi(\theta) ) - \theta$, we have 
\[
\E\left[ \expp{-\lambda \tau^{(\theta)}_r  }\right]
=\expp{- r [\psi^{-1}(\lambda+\psi(\theta) ) - \theta]}.
\]
Furthermore, this  equality holds for  $\lambda\geq -\psi(\theta)$. With
$\lambda=-\psi(\theta)$,   we    get   $\E\left[   \expp{   \psi(\theta)
    \tau^{(\theta)}_r } \right] =\expp{\theta r}$.

{F}rom  \reff{eq:cgt-proba}, we  get that  the process  $(Q_t,  t\geq 0)$,
where  $Q_t=\expp{\theta   X_t^{(\theta)}  +t  \psi(\theta)}$  is a
martingale. Since   $\displaystyle  M_{\tau_r^{(\theta)}}=
\expp{-\theta r  + \psi(\theta)  \tau^{(\theta)}_r} $ is integrable and 
$\E[M_{\tau_r^{(\theta)}}]=1$, we deduce from \reff{eq:cgt-proba} that 
\begin{equation}
   \label{eq:form-cgt}
\E\left[F(X_{\leq  \tau^{(\theta)}_r}^{(\theta)})  \expp{-\theta  r  +
    \psi(\theta)  \tau^{(\theta)}_r}  \right]=\E[F(X_{\leq  \tau_r})].
\end{equation}
Let     $\ce_i=(X_{t+\alpha_i}     -I_{\alpha_i},    t\in     [\alpha_i,
\alpha_i+\sigma_i])$,  $i\in I$,  be  the excursions  of  $X$ above  its
minimum,  up  to  time  $\tau_r$.  With  $F$  such  that  $\displaystyle
F(X_{\leq \tau_r}) = \expp{ -\sum_{i\in I} G(\ce_i)}$, we get 
\[
\E[
F(X_{\leq  \tau_r})\expp {-\lambda \tau_r} ]
=\expp{-r \N[1-\expp{-G(\ce)-\lambda \sigma}]}.
\]
We deduce from \reff{eq:form-cgt} that 
\[
\expp{-\theta r} \expp{ - r\N [1-\expp{- G(\ce^{(\theta)})+\psi(\theta) 
    \sigma^{(\theta)}}]} 
= \expp{- r \N[1-\expp{- G(\ce)}]},
\]
where $\ce^{(\theta)}$ is an excursion of $X^{(\theta)}$ above its
minimum, that is 
\[
\N [1-\expp{- G(\ce^{(\theta)})+\psi(\theta) 
    \sigma^{(\theta)}}] = \N[1-\expp{- G(\ce)}] -\theta.
\]
Subtracting $\N [1-\expp{\psi(\theta) 
    \sigma^{(\theta)}}] =  -\theta$, in the above equality, we get 
\[
 \N \left[ \expp{ \psi(\theta)\sigma^{(\theta)}} [1- \expp{-G(\ce^{(\theta)})}]
\right]
=\N\left[ 1- \expp{-G(\ce)} \right].
\]

\end{proof}


\section{Link between Lévy snake and fragmentation processes at nodes}
\label{sec:LS-F}

We define  the  fragmentation process.   Let  $\cs=(\rho,M)$ be  a
Lévy Poisson snake.  Recall definition of $m^{(\theta)}$ at the end of
Section  \ref{sec:LPS}.   For  fixed  $\theta>0$, let  us  consider  the
following  equivalence  relation  $\mathcal{R}_\theta$ on  $[0,\sigma]$,
defined under $\N$ or $\N_\sigma$ (see definition in Section
\ref{sec:prop_pruned}) by:
\begin{equation}
   \label{eq:def-R}
s\mathcal{R}_\theta                                               t\iff
m_s^{(\theta)}\bigl([H_{s,t},H_s]\bigr)=m_t^{(\theta)}
\bigl([H_{s,t},H_t]\bigr)=0,
\end{equation}
where  $\displaystyle  H_{s,t}=\inf_{u\in[s,t]}H_u$  (recall  definition
\reff{eq:def_H}). Intuitively, two points $s$ and $t$ belongs to the
same class of equivalence (i.e. the same fragment) at time $\theta$, if there is no cut on their lineage down
to their most recent common ancestor (that is
$m^{(\theta)}_s$ put no mass on $[H_{s,t}, H_s]$ nor 
$m^{(\theta)}_t$ on $[H_{s,t}, H_t]$). Notice cutting occurs on 
branching points, that is at node of the CRT. Each node of the CRT
correspond to a jump of the underlying Lévy process $X$. The
cutting times are, conditionally on the CRT, independent exponential
random times, with parameter  equal to the jump of the
corresponding node.

Let us index the different  equivalent classes in the following way: For
any $s\le  \sigma$, let  us define $H_s^0=0$  and recursively  for $k\in
\N$,
\[
H_s^{k+1}=\inf\Bigl\{u\ge 0\bigm| m_s^\theta\bigl((H_s^k,u]\bigr)>0\Bigr\},
\]
with the usual convention $\inf\emptyset=+\infty$. We  set
 \[
K_s=\sup\{j\in\N,\ H_s^j<+\infty\}.
\]
\begin{rem}
   \label{rem:ks=infini}
Notice  that we have $K_s=\infty $ if 
$M_s(\cdot, [0,\theta])$ has infinitely many atoms. By construction of
$M$ using Poisson point   measures,  this happens  $\N[d\cs]  \;
ds$-a.e.,  if and only if the intensity measure  $\rho_s
+\eta_s$ is infinite. Since $\N[d\cs] $-a.e., $\rho$ and $\eta$ are
finite measure valued process, we get that 
$ \N[d\cs] $-a.e., $K_s<\infty $. 
\end{rem}

Let us remark that $s\mathcal{R}_\theta t$ implies $ K_s=K_t$.
We denote, for any $j\in\N$, $(R^{j,k},k\in J_j)$ the family of
equivalent classes with positive Lebesgue measure such that $K_s=j$. For
$j\in\N$, $k\in J_j$ we set
\[
A_t^{j,k}  =\int_0^t\ind_{\{s\in R^{j,k}\}}ds\quad\text{and}\quad
C_t^{j,k}  =\inf\{u\ge 0,\ A_u^{j,k}>t\},
\]
with   the  convention  $\inf\emptyset=\sigma$.    And  we   define  the
corresponding Lévy snake,  $\tilde \cs^{j,k}=(\tilde \rho^{j,k},\tilde
M^{j,k})$ by: 
for every $f\in \cb_+(\R_+)$, $\varphi\in \cb_+(\R_+\times \R_+)$,
$t\geq 0$, 
\begin{align*}
\bigl\langle \tilde \rho_t^{j,k},f\bigr\rangle  &
=\int_{(H_{C_0^{j,k}},+\infty)}f(x-H_{C_0^{j,k}})\rho_{C_t^{j,k}}(dx)\\
\bigl\langle \tilde M_t^{j,k},\varphi\bigr\rangle  & =
\int_{(H_{C_0^{j,k}},+\infty)\times(\theta,+\infty)}
\varphi(x-H_{C_0^{j,k}},v-\theta)M_{C_t^{j,k}}(dx,dv) .
\end{align*}
Let $\tilde \sigma^{j,k}=A^{j,k}_\infty $ be the length of the excursion
$\tilde\cs^{j,k}$. 
Since $K_s<\infty $ $\N[d\cs] ds$-a.e. (Remark
\ref{rem:ks=infini}), the family $(\tilde \sigma^{j,k}j\in \N, k\in
J_j)$ gives all the equivalent
classes with positive Lebesgue measure.

\begin{rem}
\label{rem:L-q}
   In view of the next Section we introduce the set $\cl^{(\theta)}=(\tilde
\rho^{(j,k)}, j\in \N, k\in J_j)$ of fragments of Lévy snake as well as the 
the set $\cl^{(\theta-)}$ defined similarly but for the equivalence
relation where ${\mathcal R}_\theta$ in \reff{eq:def-R} is replaced by
${\mathcal R}_{\theta-}$ defined as
\begin{equation}
   \label{eq:def-R-}
s{\mathcal R}_{\theta-}  t\iff   M_s\bigl( [H_{s,t},H_s]\times (0,\theta)  \bigr)
=M_s\bigl( [H_{s,t},H_t]\times (0,\theta)  \bigr)=0.
\end{equation}
Notice that $m^{(\theta)}_s(\cdot)=M_s\bigl(\cdot, (0,\theta]\bigr)$. So
the two equivalence relations are equal $\N$-a.e. for fixed
$\theta$, but may differ if $M$ has an atom in  $\{\theta\}\times \R_+$. 
\end{rem}

Let   us   denote   by   $\Lambda   ^\theta=(\Lambda   _1^\theta,\Lambda
_2^\theta,\ldots)$  the sequence  of positive  Lebesgue measures  of the
equivalent classes of  $\mathcal{R}_\theta$, $(\tilde \sigma^{j,k}, j\in
\N, k\in J_j)$, ranked in decreasing order.
Notice this sequence
is at most countable. If it is finite, we complete the sequence with
zeros, so that $\N$-a.s. and $\N_\sigma$-a.s.
\[
           \Lambda ^\theta             \in\mathcal{S}^\downarrow
=\bigl\{(x_1,x_2,\ldots),\  x_1\ge  x_2\ge \cdots  \ge  0,\ \sum  x_i\le
\infty \bigr\}.
\]
For $\pi^*(d\sigma)$-a.e.  $\sigma>0$, let $ \rP_\sigma$  denote the law
of  $(\Lambda ^\theta,  \theta\geq 0)$  under $\N_\sigma$.   (The
law,  $\N_r$, of  $\cs$ conditionally  on the  length of  the excursion,
$\sigma$, being equal to $r$  has been defined in the previous Section.)
By    convention   $\rP_0$   is    the   Dirac    mass   at    $(0,   0,
\ldots)\in\mathcal{S}^\downarrow $.

\begin{theo}\label{theo:frag_property}
For   $\pi_*(dr)$-almost   any   $r$,   under   $\rP_r$,   the   process
$\Lambda =(\Lambda ^\theta,\theta\ge   0)$    is   a   $\mathcal{S}^\downarrow$-valued
fragmentation  process.  More precisely,  the law  under $\rP_r$  of the
process    $(\Lambda ^{\theta+\theta'},\theta'\ge    0)$    conditionally    on
$\Lambda ^\theta=(\Lambda _1,\Lambda _2,\ldots)$  is given  by the  decreasing  reordering of
independents processes of respective law $\rP_{\Lambda _1},\rP_{\Lambda _2},\ldots$.
\end{theo}

\begin{rem}
   We get a self-similar fragmentation when
   $\psi(\lambda)=\lambda^\alpha$, see Corollary
   \ref{cor:alpha-frag}. This particular case was already studied in
   \cite{m:sfdfstsn}. 
\end{rem}

\begin{rem}
  We  may get rid  of the  ``$\pi_*(dr)$-almost any  $r$'' and  have the
  theorem  for any  positive $r$  if we  have a  regular version  of the
  family  of conditional  probability laws  $(\N_r,r> 0)$.  This  is for
  instance the  case when  the Lévy  process is stable  (for which  it is
  possible  to construct  the measure  $\N_r$ from  $\N_1$ by  a scaling
  property)  or  when we  may  construct  this  family via  a  Vervaat's
  transform of the Lévy bridge (see \cite{m:oacfalpnpj}).
\end{rem}

The proof of the Proposition is a consequence of   Lemma
\ref{lem:frag}, and the fact that $\N(\cdot)=\int
_{(0,+\infty)}\pi_*(dr)\N_r(\cdot)$ which implies that the result of
Lemma \ref{lem:frag} holds  $\N_r$-a.s.  for
$\pi_*(dr)$ almost every $r$.

\begin{lem}
\label{lem:frag}
Under $\N$, the law of the family $(\tilde\cs^{j,k}, j\in \N, k\in J_j)$,
conditionally on $(\tilde\sigma^{j,k}, j\in \N, k\in J_j)$, is the law of
independent Lévy Poisson snakes, and the conditional law of
$\tilde\cs^{j,k}$ is $\N_{\tilde\sigma^{j,k}}$. 
\end{lem}

\begin{proof}
For  $j=0$,  notice that  $J_0$  has only  one  element,  say $0$.   And
$\tilde\cs^{0,0}$  is just the  Lévy snake,  $\tilde \cs$,  defined in
Section     \ref{sec:Markov_special}.      Of     course,    we     have
$\tilde\sigma^{0,0}=\tilde  \sigma$.  {F}rom  the special  Markov property
(Theorem  \ref{th:SMP}) and  Proposition  \ref{prop:M_trans}, we  deduce
that conditionally  on $\tilde\sigma^{0,0}$, $\tilde  \cs^{0,0}$ and the
family  $(\cs^i,i\in  I)$  of   excursions  of  $\cs$  out  of  $\{s\geq
0;m_s^\theta=0\}$ (as  defined in Section  \ref{sec:Markov_special}) are
independent.

{F}rom Corollary \ref{cor:Npruned}  and Lemma \ref{lem:loi_cond} for the
exploration   process  and   Proposition   \ref{prop:M_trans}  for   the
underlying  Poisson process,  we deduce  that, conditionally  on $\tilde
\sigma^{0,0}$,   $\tilde\cs^{0,0}$   is   distributed   according   to
$\N_{\tilde\sigma^{0,0}}$.

Furthermore,  from the special  Markov property  (Theorem \ref{th:SMP}),
the  conditional  law   of  $\cs^i$  is  given  by   $\rN$,  defined  in
\reff{eq:def_rN}.  Now  we give a  Poisson decomposition of  the measure
$\rN$.

For  $\cs'=(\rho',  M')$ distributed  according  to  $\rN$, we  consider
$(\alpha'_{l},\beta_{l})_{l\in  I'}$  the  excursion  intervals  of  the
Lévy Poisson snake, $\cs'$, out of $\{H'_s=0\}$. For $l\in I'$, we set
${\cs'}^{l}=({\rho'}^{l}, {M'}^{l})$ where for $s\geq 0$,
\begin{align*}
{\rho'}^{l}_s(dr) & ={\rho'_{(s+\alpha'_{l})\wedge
    \beta'_{l}}}(dr)\ind_{(0,+\infty)}(r),\\ 
{ M'}^{l}_s(dr,dv) & ={M'_{(s+\alpha'_{l})\wedge
    \beta'_{l}}}(dr,dv)\ind_{(0,+\infty)}(r). 
\end{align*}
Let  us   remark  that  in  the  above   definition  ${\rho'}^{l}$  and
${M'}^{l}$ don't have mass at $\{0\}$ and $\{0\}\times \R_+$.

As a direct consequence of the Poisson decomposition of $\P_\ell^*$ (see
Lemma \ref{lem:dlg-decomp}), we get  the following Lemma.

\begin{lem}
 Under $\rN$, the point measure
$\displaystyle \sum_{i'\in  I'}\delta_{{\cs'}^{i'}}$ is a  Poisson point
measure with intensity  $C_\theta\N(d\cs)$ where $ C_\theta=\int_{(0,\infty
)} (1-\expp{-\theta\ell})\ell\pi(d\ell)= \psi'(\theta) -\psi'(0)$.
\end{lem}

By  this Poisson  representation, each  process $\cs^i$  is  composed of
i.i.d.  excursions of  law $\N$.  Thus we get,  conditionally on $\tilde
\sigma^{0,0}$,  a family  $(\cs^{1,k}, k\in  J_1)$ of  i.i.d. excursions
distributed  as the  atoms of  a  Poisson point  measure with  intensity
$\tilde  \sigma^{0,0}  C_\theta  \N$.   Now,  we can  repeat  the  above
arguments  for   each  excursion   $\cs^{1,k}$,  $k\in  J_1$:   so  that
conditionally on $\tilde \sigma^{0,0}$, we can
\begin{itemize}
   \item check that $\tilde \cs^{1,k}$ is built from $\cs^{1,k}$ as $\tilde
   \cs $ from $\cs$ in Section \ref{sec:Markov_special},
 \item  get  a  family   $(\cs^{2,k',k},  k'\in  J^{k}_2)$,  which  are,
   conditionally on $\tilde \sigma^{1,k}$, distributed as the atoms of a
   Poisson point  measure with intensity $  \tilde \sigma^{1,k} C_\theta
   \N$.  and are independent of $\tilde \cs^{1,k}$.
\end{itemize}
If  we  set  $J_2=\cup_{k\in  J_1}  J^{k}_2\times \{k\}$,  we  get  that
conditionally on  $\tilde \sigma^{0,0}$, and  $(\tilde \sigma^{1,k},k\in
J_1)$,
\begin{itemize}
   \item  the  excursions  $\tilde\cs^{0,0}$ and  $(\tilde\cs^{1,k},k\in
J_1)$, are independent,
\item $ \tilde \cs ^{i,k}$ is distributed as $\N_{\tilde \sigma^{j,k}}$,
  for $j\in \{0,1\}$, $k\in J_j$,
\item  $(\cs^{2,k'}, k'\in  J_2)$,  are 
  distributed as  the atoms  of a Poisson  point measure  with intensity
  $\sum_{k\in J_1} \tilde \sigma^{1,k} C_\theta \N$, and are independent
  of $\tilde\cs^{0,0}$ and $ (\tilde\cs^{1,k},k\in J_1)$.
\end{itemize}

Eventually, the result follows by induction.

\end{proof}

Now we check there is no loss of mass during the fragmentation.

\begin{prop}
 \label{prop:tot-length}
For $\pi_*(dr)$ almost every $r$, $\rP_r$-a.s., for every $\theta\ge 0$,
$\displaystyle \sum_{i=1}^{+\infty}\Lambda _i^\theta=r$.
\end{prop}

\begin{proof}
Let $\theta> 0$. We use the notations of the proof of Theorem
\ref{theo:frag_property} and of Lemma \ref{lem:frag}. For $n\in \N$, we
have $\N$-a.e. 
\[
\sigma=\sum_{k=0}^n \sum_{j\in J_k} \tilde \sigma^{j,k}+ \int_0^\sigma
\ind_{\{K_s \geq n+1\}} \; ds.
\]
By monotone convergence, we deduce from Remark \ref{rem:ks=infini}, that
we get as $n\rightarrow +\infty $ that $\N$-a.e. 
\[
\sigma=\sum_{k=0}^\infty  \sum_{j\in J_k} \tilde \sigma^{j,k}.
\]
As  the decreasing reordering  of $(\tilde  \sigma^{j,k}, j\in  \N, k\in
J_j)$   is   $\Lambda ^\theta$,   we   get   that   $\N$-a.e.    $\displaystyle
\sum_{i=1}^{+\infty}\Lambda _i^\theta=\sigma$.      As      the      sequence
$(\sum_{i=1}^\infty  \Lambda _i ^\theta,\theta\geq  0)$ is  non  increasing, we
deduce that the previous equality holds for any $\theta>0$, $\N$-a.e.

Here again the result for $P_r$ is deduced from the one  under $\N$.

\end{proof}


\section{Dislocation measures}
\label{sec:disloc}

 Let $\lambda(\theta)$ be the mass of a tagged fragment at time $\theta$
of    the   fragmentation   process    $\Lambda $   defined    in   Theorem
\ref{theo:frag_property} (typically  the fragment or the equivalent
class which  contains 0).  A
dislocation of  this fragment occurs when $\lambda(\theta)$  has a jump.
Let   $\ct_0$  be   the  set   of  time   jumps  for   $\lambda$.   Recall
$\cs^\downarrow$ denote  the set of  non-negative non-increasing sequence
$(x_i,  i\in  \N^*)$  such  that  $\sum_{i\geq  1}  x_i<\infty  $.   For
$\theta'\in  \ct_0$,  let  $x(\theta')  = (x_i(\theta'),  i\in  \N^*)  \in
\cs^\downarrow$,  the   masses  of   the  fragments  resulting   of  the
dislocation at time $\theta'$.  Following  the Remark after Theorem 3 in
\cite{b:ssf},     we      call     the     random      point     measure
\[
\delta(d\theta,dx)=\sum_{\theta'\in        \ct_0}        \delta_{\theta',
x(\theta')}(d\theta,dx)
\]
 the dislocation process of the fragmentation (or dislocation process of
the $\psi$-CRT  fragmentation at nodes).   Of course, since there  is no
erosion,  that   is  the  total   length  is  constant   cf.  Proposition
\ref{prop:tot-length}, $\lambda(\theta'-)=\sum_{i\geq 1} x_i(\theta')$.

For self-similar fragmentation with  with index $\gamma$ and no erosion,
there   exists   a   measure   $\nu_1$  on   $   \cs^\downarrow_1=\{x\in
\cs^\downarrow;  \sum_{i\geq  1}   x_i=  1\}$,  called  the  dislocation
measure,  such that  the dislocation  process  is a  point process  with
intensity                 $\ind_{\{                \lambda(\theta-)>0\}}
\nu_{\lambda(\theta-)}(dx)d\theta$,  where the  measures  $(\nu_r, r>0)$
are defined by
\begin{equation}
   \label{eq:sel-sim}
\int_{\cs^\downarrow_r}  F(x)\nu_r(dx)
=r^\gamma\int_{\cs^\downarrow_1} F(rx)  \nu_1(dx),
\end{equation}
and  the equality  hold  for any  non-negative  measurable function  on
$\cs^\downarrow$. We refer to \cite{b:ssf} for the proof of this result and
to  \cite{j:cspm} for  the definition  of  intensity of  a random  point
measure.

In order to give the corresponding dislocation measures for the $\psi$-CRT
fragmentation at nodes, we need  to consider $(\Delta S_t, t\geq 0)$ the
jumps of  a subordinator  $S$ with Laplace  exponent $\psi^{-1}$.  Let 
$\mu$  the  measure  on   $\R_+\times  \cs  ^\downarrow$  such that   for  any
non-negative measurable function, $F$, on $\R_+\times \cs ^\downarrow$,
\begin{equation}
   \label{eq:def-mu}
\int_{\R_+\times   \cs  ^\downarrow}   F(r,x)   \mu(dr,dx)=\int  \pi(dv)
  \E[F(S_v, (\Delta S_t, t\leq v))], 
\end{equation}
where $(\Delta S_t, t\leq v)$ has to be understood as the family of
jumps of the subordinator up to time $v$ ranked in decreasing size.

Intuitively, $\mu$ is the  law of $S_T$ and the jumps of  $S$ up to time
$T$, where $T$ and $S$ are independent, and $T$ is distributed according
to  the  infinite measure  $  \pi$. Recall  $\pi_*$  is  the ``law''  of
$\sigma$  under $\N$  (this  is  the Lévy  measure  associated to  the
Laplace exponent $\psi^{-1}$).

\begin{theo}
\label{th:calcul_nu_r}
   The dislocation process of  the $\psi$-CRT fragmentation at nodes, is
   under   $\N$   a   point    process   with   intensity   $   \ind_{\{
   \lambda(\theta-)>0\}}   \nu_{\lambda(\theta-)}(dx)   d\theta$,  where
   $\lambda(\theta-)=  \sum_{i\geq 1}  x_i(\theta)$ is  the mass  of the
   fragment just before $\theta$.  And the family of dislocation measure
   $(\nu_r,   r>0)$   on  $\cs^\downarrow$   is   the   result  of   the
   disintegration of $r \mu(dr,dx)$ w.r.t. $\pi_*(dr)$:
 \[
r \mu(dr,dx)=\nu_r(dx)\pi_* (dr).
\]
\end{theo}
Notice that \reff{eq:def-mu} implies that
$\pi^*(dr)$-a.e. $\nu_r(dx)$-a.e. $\sum_{i\in \N^*} x_i=r$, where
$x=(x_i, i\in \N^*)$. The dislocation measure $\nu_r$ describe the
dislocation of a fragment of size $r$.

\begin{rem} 
   Either from Lemma \ref{lem:loi_cond} or directly, it is easy
   to check that  the dislocation measure of the 
   fragmentation at nodes 
    associated to $\psi^{(\theta)}$ (see \reff{eq:def_psi-q}),
   $(\nu_r^{(\theta)},
   r>0)$, is equal to $(\nu_r,
   r>0)$, $\pi_*(dr)$-a.e. 
\end{rem}

The next Sections are devoted to the proof of the Theorem. In Section
\ref{sec:orf}, we give an other representation of the fragmentation
following ideas  in \cite{as:psf,ap:sac} developed for
$\psi(\lambda)=\lambda^2$. In Section \ref{sec:pp}, we explain how to
compute the intensity of the dislocation process. And we perform the
computation in Section \ref{sec:comp-int}. This will end the proof of
the Theorem.

For the $\lambda^\alpha$-CRT (with $\alpha\in (1,2)$), thanks to scaling
properties, the  corresponding fragmentation is self  similar with index
$1/\alpha$, and we can recover the result of \cite{m:sfdfstsn}.

\begin{cor}
   \label{cor:alpha-frag}
{F}or the $\lambda^\alpha$-CRT fragmentation at nodes, the fragmentation
is self-similar, with index $1/\alpha$, that is \reff{eq:sel-sim} holds
with $\gamma=1/\alpha$. 
And the dislocation measure $\nu_1$
on $\cs^\downarrow_1$ is s.t. 
\[
\int F(x) \nu_1(dx)= \frac{\alpha(\alpha-1)
  \Gamma([\alpha-1]/\alpha)}{\Gamma(2-\alpha)} \E[S_1 \;F((\Delta
S_t/ S_1, t\leq 1))],
\]
holds    for   any    non-negative   measurable    function,    $F$,   on
$\cs^\downarrow_1$, where  $(\Delta S_t,  t\geq 0)$ are  the jumps  of a
stable   subordinator   $S=(S_t,   t\geq   0)$   of   Laplace   exponent
$\psi^{-1}(\lambda)=\lambda^{1/\alpha}$, ranked by decreasing size.
\end{cor}

\begin{proof}
For  $\psi(\lambda)=\lambda^\alpha$,  we  get  $\pi(dr)=\alpha(\alpha-1)
\Gamma(2-\alpha)^{-1}  r^{-1-\alpha}   dr$  as  well   as  $\pi_*  (dr)=
\left[\alpha    \Gamma    ([\alpha-1]    /\alpha)    \right]^{-1}    r^{
-(1+\alpha)/\alpha}  dr$.   In particular,  we  have  for a  non-negative
measurable function, $F$, defined on $\R_+\times \cs^\downarrow_1$,
\begin{align*}
   \int F(r,x) \;r\mu(dr,dx)
&=\E\left[\int  \pi(dv)
 \;S_v F(S_v, (\Delta S_t, t\leq v))\right]\\
&=\frac{\alpha(\alpha-1)}{\Gamma(2-\alpha)} \E\left[ \int 
  \frac{dv}{v^{1+\alpha}}\; S_vF(S_v, (\Delta
S_t, t\leq v))\right ]\\
&=\frac{\alpha(\alpha-1)}{\Gamma(2-\alpha)} \E\left[ \int 
  \frac{dv}{v}\; S_1F(v^\alpha S_1, v^\alpha S_1 ( \Delta
S_t/S_1, t\leq 1))\right ]\\
&= \frac{\alpha-1}{\Gamma(2-\alpha)}
\int \E[S_1 \;F(y, y(\Delta
S_t/S_1, t\leq 1))] \frac{dy}{y},
\end{align*}
where we used the  scaling property of $S$, that is 
$(\Delta S_t, t\leq  r)$ is distributed as  $(r^\alpha \Delta S_t,
t\leq  1)$, for the third equality, and the change of variable
$y=v^\alpha S_1$ for the fourth equality. {F}rom Theorem
\ref{th:calcul_nu_r}, we have  that
\[
\int \inv{ \alpha\Gamma([\alpha-1]/\alpha)}
\frac{dr}{r^{(1+\alpha)/\alpha} } \nu_r(dx)\; F(r,x)=
\int \frac{\alpha-1}{\Gamma(2-\alpha)}
\E[S_1 \;F(y, (y\Delta
S_t, t\leq 1))] \frac{dy}{y}.
\]
This implies that for a.a. $r>0$, 
\[
\int \nu_r(dx) \; F(x)= \frac{\alpha(\alpha-1)
  \Gamma([\alpha-1]/\alpha)}{\Gamma(2-\alpha)} r^{1/\alpha} \E[ S_1 F(r
  (\Delta S_t/S_1, t\leq 1))],
\]
and  thus  $\displaystyle   \int  \nu_r(dx)  \;  F(x)=r^{1/\alpha}  \int
\nu_1(dx)\; F(rx)$, with
\[
\int \nu_1(dx)\; F(x)
=\frac{\alpha(\alpha-1)
  \Gamma([\alpha-1]/\alpha)}{\Gamma(2-\alpha)}  \E[ S_1 F(
  (\Delta S_t/S_1, t\leq 1))].
\]
\end{proof}

\subsection{An other representation of the fragmentation}
\label{sec:orf}

Following the ideas  in \cite{as:psf,ap:sac},  we give  an other
representation  of  the   fragmentation  process  described  in  Section
\ref{sec:LS-F}, using a Poisson point  measure under the epigraph of the
height  process.   

We consider a fragmentation process,  as time $\theta$ increases, of the
CRT,   by  cutting   at  nodes   (set  of   points  $(s,a)$   such  that
$\kappa_s(\{a\})>0$, where $\kappa$  is defined in \reff{eq:def-kappa}).
More precisely, we consider, conditionally on the CRT or equivalently on
the exploration process $\rho$,  a Poisson point process, $ \cq(d\theta,
ds,  da)$   under  the  epigraph  of  $H$,   with  intensity  $d\theta\;
q_\rho(ds,da)$, where
\begin{equation}
  \label{eq:intensity_d} 
 q_\rho(ds,da)= \frac{ds\; \kappa_s(da)}{d_{s,a}-g_{s,a}},
\end{equation}
with  $d_{s,a}=\sup\{u\geq  s,  \min\{H_v,  v\in  [s,u]\}  \geq  a\}  $
and  $g_{s,a}=\inf\{u\leq  s, \min\{H_v, v  \in [u,s]\} \geq  a\}$. (The
set $[g_{s,a}, d_{s,a}]\subset [0,\sigma]$ represent the individuals who
have a common ancestor with the individual $s$ after or at generation $a$.)

Notice that  from this  representation, the cutting  times of  the nodes
are, conditionally on the  CRT, independent exponential random time, and
their parameter is equal to the mass of the node (defined as the mass of
$\kappa$ or equivalently  as the value of the  jump of $X$ corresponding
to the given  node).

We say two points $s, s'\in  [0,\sigma]$ belongs to the same fragment at
time $\theta$,  if there is no cut  on their lineage down  to their most
recent common ancestor $H_{s,s'}$: that is for $v=s$ and $v=s'$,
\[ 
\int\ind_{[H_{s,s'}, H_v]}(a)\ind_{[g_{v,a}, d_{v,a}]}(u) 
 \cq ([0,\theta] ,du,da)=0 .
\]
This define an equivalence relation,  and we call fragment an equivalent
class.  Let $\Lambda ^\theta$  be the  sequences  of Lebesgue  measures of  the
corresponding equivalent classes ranked in decreasing order.

It  is clear  that conditionally  on  the CRT,  the process  $(\Lambda ^\theta,
\theta\geq  0)$ as the  same distribution  as the  fragmentation process
defined  in  Section  \ref{sec:LS-F}.   Roughly  speaking,  in  Section
\ref{sec:LPS}   (which   leads   to   the   fragmentation   of   Section
\ref{sec:LS-F}) we  mark the node as  they appear: that is,  for a given
level   $a$,  the   node   $\{s;  \kappa_s(\{a\})>0\}$   is  marked   at
$g_{s,a}$. Whereas in this Section  the same node is marked uniformly on
$[g_{s,a}, d_{s,a}]$. In both case,  the cutting times of the nodes are,
conditionally on the CRT, independent exponential random time, and their
parameter is equal to the mass  of the node (defined as the common value
of   $\kappa_u(\{a\})$   for   $u\in   \{s;   \kappa_s(\{a\})>0\}$,   or
equivalently  as the  value  of the  jump  of $X$  corresponding to  the
given  node).

Now, we define the fragments of the
Lévy snake corresponding to the cutting of $\rho$ according to
the measure $q_\rho$. For $(s,a)$ chosen according to
the measure $q_\rho(ds, da)$, we can define the following Lévy snake
fragments $(\rho^i, i\in \tilde I)$ of $\rho$ by considering 
\begin{itemize}
   \item the open intervals of excursion after $s$ of $H$ above level $a$:
   $((\alpha_i,\beta_i), i\in \tilde I_+)$, which are such that
   $\alpha_i>s$, $H_{\alpha_i}=H_{\beta_i}=a$, and for $s'\in
   (\alpha_i, \beta_i)$ we have $H_{s'}>a$ and
   $H_{s',s}=a$ (recall definition  \reff{eq:def_H});
   \item the open intervals of excursion before  $s$ of $H$ above level $a$:
   $((\alpha_i,\beta_i), i\in \tilde I_-)$, which are such that
   $\beta_i<s$, $H_{\alpha_i}=H_{\beta_i}=a$, and for $s'\in
   (\alpha_i, \beta_i)$ we have $H_{s'}>a$ and
   $H_{s',s}=a$;
   \item the excursion, $i_s$,  of $H$ above level $a$ that straddle $s$:
   $(\alpha_{i_s}, \beta_{i_s})$, which is such that $\alpha_{i_s}<s<
   \beta_{i_s} $, $H_{\alpha_{i_s}}=H_{\beta_{i_s}}=a$, and for $s'\in
   (\alpha_{i_s}, \beta_{i_s})$ we have $H_{s'}>a$ and
   $H_{s',s}=a$;
   \item the excursion, $i_0$,  of $H$ under level $a$: $\{s\in [0,\sigma];
   H_{s',s}<a\}=[0, \alpha_{i_0})\cup(\beta_{i_0}, \sigma]$. 
\end{itemize}
For $i\in \tilde I_+ \cup \tilde I_-\cup\{i_s\}$, we set
$\rho^i=(\rho^i_s, s\geq 0)$ where 
\[
\int f(r) \rho^i_s(dr)=\int f(r-a)\ind_{\{r>a\}} \rho_{(\alpha_i+s)\wedge
  \beta_i}  (dr)
\]
for  $f\in  \cb_+(\R)$.  For  $i_0$,  we set  $\rho^{i_0}=(\rho^{i_0}_s,
s\geq   0)$  where  $\rho^{i_0}_s=   \rho_s$  if   $s<\alpha_{i_0}$  and
$\rho^{i_0}_s=  \rho_{s-\beta_{i_0}+\alpha_{i_0}}$  if  $s>\beta_{i_0}$.
Eventually,  we  set  $\tilde  I=\tilde I_+  \cup  \tilde  I_-\cup\{i_s,
i_0\}$. And $(\rho^i, i\in \tilde I)$ correspond to the fragments of the
Lévy snake corresponding to the cutting of $\rho$ according to one point
chosen with the measure $q_\rho$.  We shall denote $\tilde \nu_\rho$ the
distribution of $(\rho^i, i\in \tilde I)$ under $\N$.

In Section  \ref{sec:comp-int}, we shall  use $\sigma^i$, the  length of
fragment  $\rho^i$.  For $i\in  \tilde  I_-  \cup  \tilde I_+$,  we  have
$\sigma^i=\beta_i-\alpha_i$.            We           also           have
$\sigma^{i_s}=\sigma^{i_s}_-+\sigma^{i_s}_+$
(resp.        $\sigma^{i_0}=\sigma^{i_0}_-+\sigma^{i_0}_+$),       where
$\sigma^{i_s}_-=s-\alpha_{i_s}  $  (resp. $\sigma^{i_0}_-=\alpha_{i_0}$)
is the length of the fragment before $s$ and
$\sigma^{i_s}_+=\beta_{i_s}-s   $  (resp. $\sigma^{i_0}_+=\sigma-\beta_{i_0}$)
is the length of the fragment after $s$. 
Notice that $\N$-a.e. $\sigma=\sum_{i\in \tilde I} \sigma^i$.

\subsection{The dislocation process is a point process}
\label{sec:pp}

Let  $\ct$ the  set of  time  jumps of  the Poisson  process $\cq$.  For
$\theta\in \ct$, consider $\cl^{(\theta-)}=(\rho_i, i\in I^{(\theta-)})$
and $\cl^{(\theta)}=(\rho_i, i\in  I^{(\theta)})$ the families of Lévy
snakes defined in Remark \ref{rem:L-q}. The length, ranked in decreasing
order,   of  those  families   of  Lévy   snakes  correspond
respectively  to  the
fragmentation   process  just   before   time  $\theta$   and  at   time
$\theta$. Notice that for $\theta\in \ct$ the families $\cl^{(\theta-)}$
and  $\cl^{(\theta)}$ agree  but  for only one  snake  $ \rho^{i_\theta}  \in
\cl^{(\theta-)}$  which  fragments  in  a  family  $(\rho^i,i\in  \tilde
I^{(\theta)})\subset \cl^{(\theta)}$. Thus we have that
\[
\cl^{(\theta)}=\left(\cl^{(\theta-)}\backslash \{
  \rho^{i_\theta}\}\right) \bigcup (\rho^i,i\in  \tilde
I^{(\theta)}).
\]
{F}rom  the
representation of  the previous Section, this fragmentation  is given by
cutting  the Lévy  snake according  to the  measure $q_\rho$:  that is
 the measure $\tilde \nu_\rho$ defined at the end of
Section \ref{sec:orf}.
{F}rom Lemma \ref{lem:frag} and the construction of the Lévy Poisson
Snake, we deduce that 
\[
\sum_{\theta\in \ct} \delta_{\theta, \cl^{(\theta-)}, (\rho^i,i\in  \tilde
I^{(\theta)})}
\]
is a point process with intensity $d\theta\; \delta_{\cl^{(\theta-)}}
\sum_{\rho\in \cl^{(\theta-)}} \tilde  \nu _\rho$.

Notice the evolution of a tagged fragment of the Lévy snake has the
same distribution as the evolution of the fragment  of the Lévy snake which
contains $0$, say $\rho^{0,(\theta)}\in \cl^{(\theta)}$. (This is known
as the re-rooting property of the CRT.)
Then,  we get that 
\[
\sum_{\theta\in \ct_0} \delta_{\theta,  (\rho^i,i\in
  \tilde I^{(\theta)})}
= \sum_{\theta\in \ct} \delta_{\theta,  (\rho^i,i\in
  \tilde I^{(\theta)})}\ind_{\{0\text{ belongs to } (\rho^i,i\in
  \tilde I^{(\theta)})\}} ,
\]
where $\ct_0$ is   the  set   of  time   fragmentation of the fragment
which contains 0, is a point process with intensity $d\theta\;  \tilde \nu
_{\rho^{0,(\theta-)}}$.

Now,    in    the   dislocation    process    of   the    fragmentation,
$\delta(d\theta,dx)=\sum_{\theta'\in       \ct_0}       \delta_{\theta',
x(\theta')}(d\theta,dx)$,  the  sequences  $x(\theta')$ are  the  length
ranked in  decreasing order, $(\sigma^i, i\in  \tilde I^{(\theta)})$, of
the  Lévy   snakes  $(\rho^i,i\in  \tilde   I^{(\theta)})$.   Using  a
projection argument, one can check that the dislocation process is
a point process with intensity $d\theta\; \nu _{\sigma^{0,(\theta-)}}$,
where $\sigma^{0,(\theta-)}$ is the  length of $ \rho^{0,(\theta-)}$ and
$ \nu  _{\sigma^{0,(\theta-)}}$ is  the distribution of  the decreasing
lengths  of  Lévy  snakes  under $\tilde  \nu  _{\rho^{0,(\theta-)}}$,
integrated w.r.t.  to the law of  $\rho^{0,(\theta-)}$ conditionally on
$\sigma^{0,(\theta-)}$. More precisely we have $\pi_*(dr)$-a.e.
\[
\int_{\cs^\downarrow}  F(x) \nu_r(dx)=\N_r\left[\int F((\sigma^i, i\in
\tilde I)) \tilde \nu _\rho (d(\rho^i, i\in \tilde I))\right],
\]
for    any   non-negative    measurable   function    $F$    defined   on
$\cs^\downarrow$,  where   $(\sigma^i,  i\in  \tilde  I)$   as  to  be
understood  as the  family of  length, of  the fragments  $(\rho^i, i\in
\tilde  I)$,  ranked  in  decreasing  size.  

This prove that the dislocation process is a point process. And we will
now explicit the family of dislocation measures $(\nu_r, r>0)$.

As  computations are  more tractable  under $\N$  than under  $\N_r$, we
shall  compute  for  $\lambda\geq  0$, and  any  non-negative  measurable
function, $F$,  defined on $\cs^\downarrow$
\[
\int_{ \R_+\times \cs^\downarrow}  \expp{-\lambda r} F(x)
\pi_*(dr)\nu_r(dx) .
\]
{F}rom the definition of $\tilde \nu_\rho$, and using the notation at the end of
Section  \ref{sec:orf}, we get that this last quantity is equal to 
\begin{equation}
   \label{eq:N-q}
\N\left[\expp{-\lambda\sigma} \int  q_\rho( ds,da)  F((\sigma^i, i\in
  \tilde I)) \right],
\end{equation}
where $(\sigma^i, i\in  \tilde I)$ as to be understood  as the family of
  length ranked  in
  decreasing size.

\subsection{Computation of dislocation measures}
\label{sec:comp-int}

In order to compute  quantities like \reff{eq:N-q}, we shall consider
for $p>0$, $p'>0$ and $h\in \cb_+(\cm((0,+\infty )))$
\[
A=\N\left[\expp{-\lambda\sigma} \int  q_\rho( ds,da)  (\sigma^{i_s}+
    \sum_{i\in \tilde     I_- \cup \tilde I_+} \sigma^i ) 
\expp{ - p \sigma^{i_0} - p'\sigma^{i_s}} h\big(\sum_{i\in \tilde     I_-
    \cup \tilde I_+} \delta _{\sigma^i}\big) \right].
\]
As        $\displaystyle       q_\rho(       ds,da)        =       \frac
{\kappa_s(da)}{d_{s,a}-g_{s,a}}$                and                since
$d_{s,a}-g_{s,a}=\sigma^{i_s}+  \sum_{i\in \tilde  I_- \cup  \tilde I_+}
\sigma^i $, we get
\[
A=  \N\left[\int_0^\sigma   ds  \int  {\kappa_s(da)}\expp{  -(p+\lambda)
   \sigma^{i_0} -(p'+\lambda) \sigma^{i_s}} h\big(\sum_{i\in \tilde I_-
   \cup  \tilde  I_+} \delta_{\sigma^i}\big)\expp{-\lambda  \sum_{i\in
   \tilde I_- \cup \tilde I_+} {\sigma^i}} \right] .
\]
We set $\displaystyle h_{(\lambda)} \big(\sum_{i\in \tilde I_- \cup \tilde I_+}
\delta_{\sigma^i}\big)=h\big(\sum_{i\in \tilde I_- \cup \tilde I_+}
\delta_{\sigma^i}\big)\expp{-\lambda \sum_{i\in \tilde I_- \cup \tilde I_+}
{\sigma^i}} $. Now, we can replace
\begin{align*}
   B
&=\expp{ -(p+\lambda) \sigma^{i_0} -(p'+\lambda)  \sigma^{i_s}}
  h_{(\lambda)}\big(\sum_{i\in \tilde     I_- \cup \tilde I_+}
   \delta_{\sigma^i}\big)\\
&=\expp{ -(p+\lambda) \sigma^{i_0}_- -(p'+\lambda)  \sigma^{i_s}_-
   -(p+\lambda) \sigma^{i_0}_+ -(p'+\lambda)  \sigma^{i_s}_+  }
  h_{(\lambda)}\big(\sum_{i\in \tilde     I_- } 
   \delta_{\sigma^i}+ \sum_{i\in \tilde     I_+}
   \delta_{\sigma^i}\big) 
\end{align*}
by its optional projection $B'$:
\begin{multline*}
   B'=\expp{ -(p+\lambda) \sigma^{i_0}_- -(p'+\lambda)  \sigma^{i_s}_-}\\
\E^*_{\rho_s} \Big[\expp{ -(p+\lambda) \int_0^\sigma \ind_{\{H_{0,u}<a\}} \;
  du  -(p'+\lambda) \int_0^\sigma \ind_{\{H_{0,u}>a\}} \; du }
  h_{(\lambda)}\big(\mu' + \sum_{j\in \tilde I_+} \delta_{\sigma^j}\big)  \Big]
_{\Big| \mu'= \sum_{i\in \tilde I_-} \delta_{\sigma^i}} .
\end{multline*}
Using notations introduced above Lemma \ref{lem:dlg-decomp}, we
have 
\begin{multline*}
   B'=\expp{ -(p+\lambda) \sigma^{i_0}_- -(p'+\lambda)  \sigma^{i_s}_-}\\
\E^*_{\rho_s} \Big[\expp{ -(p+\lambda) \sum_{k\in I} \sigma_k
  \ind_{\{h_k<a\}}   -(p'+\lambda)\sum_{k\in I} \sigma_k
  \ind_{\{h_k> a\}}   }
  h_{(\lambda)}\big(\mu' + \sum_{k\in I; h_k=a} \delta_{\sigma_k}\big) \Big]
_{\Big| \mu'= \sum_{i\in \tilde I_-} \delta_{\sigma^i}} . 
\end{multline*}
Then we deduce from Lemma \ref{lem:dlg-decomp}, that 
\begin{align*}
   B'
&=\expp{ -(p+\lambda) \sigma^{i_0}_- -(p'+\lambda)  \sigma^{i_s}_-}\\
&\hspace{1cm}\expp{ -\rho_s([0,a)) \N[1-\expp{-(p+\lambda)\sigma}]
  -\rho_s((a,+\infty )) \N[1-\expp{- (p'+\lambda)\sigma}]}
  \E[h_{(\lambda)}(\mu'+ \cp)]_{\Big| \mu'= \sum_{i\in \tilde I_-}
    \delta_{\sigma^i}} \\
&=\expp{ -(p+\lambda) \sigma^{i_0}_- -(p'+\lambda)  \sigma^{i_s}_-}
\expp{ -\rho_s([0,a)) \psi^{-1}(p+\lambda)
  -\rho_s((a,+\infty )) \psi^{-1}(p'+\lambda)}
  \E[h_{(\lambda)}(\mu'+ \cp)]_{\Big| \mu'= \sum_{i\in \tilde I_-}
    \delta_{\sigma^i}} ,
\end{align*}
where  $\cp$ is  under  $ \P$  a  Poisson point  measure with  intensity
$\rho_s(\{a\})    \N[d\sigma]=\rho_s(\{a\})    \pi_*(dr)$.    By    time
reversibility (see Corollary 3.1.6 in \cite{dlg:rtlpsbp}), we get
\begin{align*}
A
&=   \N\Bigg[\int_0^\sigma ds \int
   {\kappa_s(da)} \expp{ -(p+\lambda) \sigma^{i_0}_- -(p'+\lambda)
     \sigma^{i_s}_-}\\ 
&\hspace{3cm} \expp{ -\rho_s([0,a)) \psi^{-1}(p+\lambda)
  -\rho_s((a,+\infty )) \psi^{-1}(p'+\lambda)}
  \E[{h_{(\lambda)}(\mu'+ \cp)}]_{\Big| \mu'= \sum_{i\in \tilde I_-}
    \delta_{\sigma^i}} \Bigg] \\
&=   \N\Bigg[\int_0^\sigma ds \int
   {\kappa_s(da)} \expp{ -(p+\lambda) \sigma^{i_0}_+ -(p'+\lambda)
     \sigma^{i_s}_+}\\ 
&\hspace{3cm} \expp{ -\eta_s([0,a)) \psi^{-1}(p+\lambda)
  -\eta_s((a,+\infty )) \psi^{-1}(p'+\lambda)}
  \E[{h_{(\lambda)}(\mu'+ \cp')}]_{\Big| \mu'= \sum_{i\in \tilde I_+}
    \delta_{\sigma^i}} \Bigg] ,
\end{align*}   
where $\cp'$ is under $ \P$ a Poisson point measure with intensity
$\eta_s(\{a\}) \pi_*(d\sigma)$.
Using the same computation as above, we eventually get 
\[
A=\N\left[\int_0^\sigma ds \int
   {\kappa_s(da)}
\expp{ -\kappa_s([0,a)) \psi^{-1}(p+\lambda)
  -\kappa_s((a,+\infty ))\psi^{-1}(p+\lambda)}
  \E[{h_{(\lambda)}(\cp'')}]\right ] ,
\]
where $\cp''$ is under $ \P$ a Poisson point measure with intensity
$\kappa_s(\{a\}) \pi_*(d\sigma)$.
We write $g_\lambda(\kappa_s(\{a\}))$ for $\E[{h_{(\lambda)}(\cp'')}]$. 
Thanks to the Poisson
representation of Proposition \ref{prop:poisson_representation1}, we get
\begin{align*}
   A
&= \E\left[ \int_0^\infty  da\expp{- \alpha_0 a  }
\sum_{x_i\leq  a} \ell_i 
\expp{ -\sum_{x_j<x_i} \ell_j \psi^{-1}(p+\lambda)
     - \sum_{a\geq x_j>x_i} \ell_j \psi^{-1}(p'+\lambda)} g_\lambda
   (\ell_i)  \right] \\
&=\int_0^\infty  da\expp{- \alpha_0 a  } \E\left[
\sum_{x_i\leq  a} \ell_i 
\expp{ -x_i \int \ell\pi(d\ell) \;   [1-\expp{-\ell\psi^{-1}(p+\lambda)}]  
     - (a-x_i)  \int \ell\pi(d\ell) \;  
     [1-\expp{-\ell\psi^{-1}(p'+\lambda)}]   } 
     g_\lambda   (\ell_i) \right] \\
&=\int_0^\infty  da \; \E\left[ 
\sum_{x_i\leq  a} \ell_i 
\expp{ -x_i \psi'(\psi^{-1}(p+\lambda))  
     - (a-x_i)  \psi'(\psi^{-1}(p'+\lambda)) }
g_\lambda   (\ell_i) \right] \\
&= \int_0^\infty  da
\int_{(0,\infty )} \ell  \pi(d\ell) \int dt\;  \ind_{[0,a]}(t) \; \ell 
\expp{ -t\psi'(\psi^{-1}(p+\lambda)) - (a-t)
  \psi'(\psi^{-1}(p'+\lambda))} g_\lambda(\ell)\\
&= \frac{\int_{(0,\infty )} \ell^2g_\lambda(\ell)\;  \pi(d\ell)}
{\psi'(\psi^{-1}(p+\lambda))
  \psi'(\psi^{-1}(p'+\lambda))},
\end{align*}
where we used \reff{eq:def-psi1} for the fourth equality.

On the other side, let $(\Delta S_t, t\geq 0)$ be the jumps of a
subordinator $S=(S_t, t\geq 0)$ with Laplace exponent
$\psi^{-1}$ and Lévy measure $\pi_* $. Standard computations yield
for $r>0$, 
\begin{align*}
 G(r)
&=  \E\left[\expp{-\lambda S_r} \sum_{t\leq r, s\leq  r, \; t\neq s}
  \Delta S_t \Delta S_s 
   \expp{ - p \Delta S_t - p'\Delta S_s } h\big( \sum_{u\leq
   r, u\not \in \{s,t\}} \delta_{\Delta S_u}\big)\right]\\
&=  \E\left[\sum_{t\leq r, s\leq  r, \; t\neq s} \Delta S_t \Delta S_s
   \expp{ - (p+\lambda) \Delta S_t -(p'+\lambda) \Delta S_s} h_{(\lambda)}\big(
     \sum_{u\leq    r, u\not \in \{s,t\}} \delta_{\Delta S_u}\big) \right]\\
  & = r^2  \left[\int \pi_*  (d\ell) \; \ell \expp{- (p+\lambda)\ell}\right]
\left[\int \pi_*  (d\ell) \; \ell \expp{- (p'+\lambda)\ell}\right]
\E\left[{ h_{(\lambda)}\big(     \sum_{u\leq    r} \delta_{\Delta S_u}\big)  }
 \right] \\
  & = r^2  {\psi^{-1}}' (p+\lambda){\psi^{-1}}'  (p'+\lambda) g_\lambda(r),
\end{align*}
as $\sum_{u\leq    r} \delta_{\Delta S_u}$ is a Poisson measure with
intensity $r\pi_*(d v)$. Notice that ${\psi^{-1}}'=1/\psi'\circ \psi^{-1}$ 
to conclude that 
\[
\int \pi(dr ) G(r)= \frac{\int_{(0,\infty )} r^2g_\lambda(r)\;  \pi(dr)
  }{\psi'(\psi^{-1}(p+\lambda)) 
  \psi'(\psi^{-1}(p'+\lambda))}.
\]
Therefore, we deduce that for any $p>0$, $p'>0$ and $h\in
\cb_+(\cm((0,+\infty )))$, we have 
\begin{multline*}
\N\left[\expp{-\lambda\sigma} \int  q_\rho( ds,da)  (\sigma^{i_s}+
    \sum_{i\in \tilde     I_- \cup \tilde I_+} \sigma^i ) 
\expp{ - p \sigma^{i_0} - p'\sigma^{i_s}} h\big(\sum_{i\in \tilde     I_-
    \cup \tilde I_+} \delta _{\sigma^i}\big) \right]\\   
= \int \pi(dr )\E\left[\expp{-\lambda S_r} \sum_{t\leq r, s\leq  r, \; t\neq s}
  \Delta S_t \Delta S_s 
   \expp{ - p \Delta S_t - p'\Delta S_s} h\big( \sum_{u\leq
   r, u\not \in \{s,t\}} \delta_{\Delta S_u}\big) \right] .
\end{multline*}
Recall $\tilde I= \tilde     I_- \cup \tilde I_+ \cup\{i_0, i_s\}$. 
{F}rom monotone class Theorem, we deduce that  for any $h\in
\cb_+(\R_+\times \R_+\times \cm((0,+\infty )))$,
\begin{multline*}
 \int \pi_* (dr)\expp{-\lambda r}   \N_r \left[
  \int q(ds,da)  h\big(\sigma^{i_0}, \sigma^{i_s}, \sum_{i\in \tilde
  I-\cup \tilde I_+ } \delta_{\sigma^i}\big)   \right]\\
= \int \pi(dr )\E\left[ S_r \expp{-\lambda S_r} \sum_{t\leq r, s\leq
   r, \; t\neq s} \frac{\Delta S_t}{S_r} \frac{\Delta S_s}{S_r-\Delta
   S_t}  h\big( \Delta S_t, \Delta S_s, \sum_{u\leq r, u\not \in \{s,
  t\} } \delta_{\Delta S_u} \big)
   \right].
\end{multline*}
For a measurable  non-negative function $F$ defined on $\cs^\downarrow$,
we deduce that 
\[
 \int \pi_* (dr)\expp{-\lambda r}   \N_r \left[
  \int q(ds,da)  F\big((\sigma^i, i\in \tilde I)\big)  \right]
= \int \pi(dr )\E\left[ S_r \expp{-\lambda S_r} F\big(( \Delta S_u,
  u\leq r )\big) 
   \right],
\]
where $(\sigma^i, i\in  \tilde I)$ and $( \Delta S_u, u\leq  r )$ are to
be  understood as the  family of  length or  jumps ranked  in decreasing
size.  {F}rom the end of Section \ref{sec:pp}, we deduce that
\[
   \int_{ \R_+\times \cs^\downarrow}  \expp{-\lambda r} F(x)
\pi_*(dr)\nu_r(dx) 
= \int \pi(dv )\E\left[ S_v \expp{-\lambda S_v} F\big(( \Delta S_u,
  u\leq v )\big) 
   \right].
\]
{F}rom definition \reff{eq:def-mu} of $\mu$,  we deduce that 
\[
   \int_{ \R_+\times \cs^\downarrow}  \expp{-\lambda r} F(x)
\pi_*(dr)\nu_r(dx) 
= \int \expp{-\lambda r}  F(x)\;  r \mu (dr, dx).
\]
This ends the proof of theorem \ref{th:calcul_nu_r}.

\bibliographystyle{abbrv}
\bibliography{/home/delmas/cermics/Bibliographie/delmas}

\end{document}